\documentclass{article}

\usepackage{latexsym}

\newcommand{\bremark}{\begin{remark}}
\newcommand{\eremark}{\end{remark}}

\newcommand{\blem}{\begin{lemma}}
\newcommand{\elem}{\end{lemma}}
\newcommand{\bth}{\begin{theorem}}
\newcommand{\eth}{\end{theorem}}
\newcommand{\benu}{\begin{enumerate}}
\newcommand{\eenu}{\end{enumerate}}
\newcommand{\bdes}{\begin{description}}
\newcommand{\edes}{\end{description}}
\newcommand{\bdf}{\begin{definition}}
\newcommand{\edf}{\end{definition}}
\newcommand{\bcor}{\begin{cor}}
\newcommand{\ecor}{\end{cor}}
\newcommand{\bprp}{\begin{proposition}}
\newcommand{\eprp}{\end{proposition}}
\newcommand{\bmlem}{\begin{mlemma}}
\newcommand{\emlem}{\end{mlemma}}
\newcommand{\bclm}{\begin{claim}}
\newcommand{\eclm}{\end{claim}}
\newcommand{\bprf}{{\bf Proof}.\hspace{2mm}}

\newcommand{\eprf}{\hspace*{\fill} $\Box$}

\newcommand{\ovl}{\overline}
\newcommand{\beqn}{\begin{equation}}
\newcommand{\eeqn}{\end{equation}}
\newcommand{\beqnarr}{\begin{eqnarray}}
\newcommand{\eeqnarr}{\end{eqnarray}}
\newcommand{\beqnarrs}{\begin{eqnarray*}}
\newcommand{\eeqnarrs}{\end{eqnarray*}}

\newcommand{\smlskp}{\\ \smallskip \\}
\newcommand{\spand}{\,\&\,}

\newtheorem{theorem}{Theorem}[section]
\newtheorem{definition}[theorem]{Definition}
\newtheorem{proposition}[theorem]{Proposition}
\newtheorem{lemma}[theorem]{Lemma}
\newtheorem{cor}[theorem]{Corollary}
\newtheorem{mlemma}[theorem]{Main Lemma}
\newtheorem{claim}[theorem]{Claim}

\newtheorem{remark}[theorem]{Remark}

\newcommand{\alp}{\alpha}

\newcommand{\veps}{\varepsilon}
\newcommand{\del}{\delta}
\newcommand{\Del}{\Delta}
\newcommand{\ome}{\omega}
\newcommand{\Ome}{\Omega}
\newcommand{\bet}{\beta}
\newcommand{\gam}{\gamma}
\newcommand{\Gam}{\Gamma}
\newcommand{\kap}{\kappa}
\newcommand{\sig}{\sigma}
\newcommand{\Sig}{\Sigma}
\newcommand{\tht}{\theta}

\newcommand{\vphi}{\varphi}

\newcommand{\fal}{\forall}
\newcommand{\exi}{\exists}
\newcommand{\rarw }{\rightarrow}
\newcommand{\Rarw }{\Rightarrow}
\newcommand{\larw}{\leftarrow}
\newcommand{\Larw}{\Leftarrow}
\newcommand{\lrarw}{\leftrightarrow}
\newcommand{\Lrarw}{\Leftrightarrow}
\newcommand{\uarw}{\uparrow}
\newcommand{\darw}{\!\downarrow}

\newcommand{\calb}{{\cal B}}
\newcommand{\calc}{{\cal C}}
\newcommand{\cald}{{\cal D}}

\newcommand{\calg}{{\cal G}}

\newcommand{\calt}{{\cal T}}
\newcommand{\calu}{{\cal U}}
\newcommand{\calv}{{\cal V}}

\newcommand{\calw}{{\cal W}}
\newcommand{\calx}{{\cal X}}
\newcommand{\caly}{{\cal Y}}
\newcommand{\calz}{{\cal Z}}

\newcommand{\incl}{\subseteq}

\title{
Wellfoundedness proofs by means of non-monotonic inductive definitions II: \\
first order operators
}
\author{Toshiyasu Arai
\\
Graduate School of Science,
Chiba University
\\
1-33, Yayoi-cho, Inage-ku,
Chiba, 263-8522, JAPAN
}
\date{Nov. 23, 2005}
\begin{document}
\maketitle
\begin{abstract}
In this paper, we give two proofs of the wellfoundedness of recursive notation 
systems for $\Pi_{N}$-reflecting ordinals. 
One is based on $\Pi_{N-1}^{0}$-inductive definitions, and the other is based on distinguished classes.
\end{abstract}

\section{Introduction}\label{sec:intro}

This is a sequel to our \cite{Wienpi3d}.
In \cite{Wienpi3d} we proved the wellfoundedness of recursive notation 
systems for reflecting ordinals up to $\Pi_{3}$-reflection by relevant inductive definitions.

Let KP$\Pi_{N}\, (3<N<\omega)$ denote a set theory whose intended models 
are $L_{\pi}$ with $\Pi_{N}$-reflecting ordinals $\pi$.
It is easy to see that KP$\Pi_{N}$ proves that any $\Pi^{0}_{N-1}$-inductive definition
eventually reaches to a closed point. 
We have designed a recursive notation system $Od(\Pi_{N})$ of ordinals
so that the order type of its countable fragment is the proof-theoretic 
ordinal of KP$\Pi_{N}$. An element of the notation system is called an {\it ordinal diagram\/}.
One half of this result was accomplished by cut-elimination in \cite{ptpiN}.
The other was to show that KP$\Pi_{N}$ proves the wellfoundedness of
$Od(\Pi_{N})$ up to each countable ordinal in it

In this paper, we give two proofs of the wellfoundedness.
One is based on $\Pi_{N-1}^{0}$-inductive definitions, and the other is based on distinguished classes.
Proof theoretic study for $\Pi_{N}$-reflecting ordinals via ordinal diagrams $Od(\Pi_{N})$ will be reported in a forthcoming paper \cite{ptpiN}.

$Ord$ denotes the class of ordinals. $\Pi^{0}_{n}$ denotes the arithmetical hierarchy on $\ome$, while $\Pi_{n}$ the L\'evy hierarchy on sets.

\bdf\label{df:operation}{\rm (Richter-Aczel \cite{Richter-Aczel74})}
\benu
\item {\rm By an} operator {\rm we mean any function} $\Gam$ {\rm on} ${\cal P}(\ome)$.
\item {\rm An operator} $\Gam$ {\rm determines a transfinite sequence} $(\Gam^{x}: x\in Ord)$ {\rm of subsets of} $\ome${\rm , where} $\Gam^{x}=\bigcup\{\Gam(\Gam^{y}):y<x\}$.  {\rm The} closure ordinal $|\Gam|$ {\rm of} $\Gam$ {\rm is the least ordinal} $x$ {\rm such that} $\Gam^{x+1}=\Gam^{x}$. {\rm The set} defined by $\Gam$ {\rm is} $\Gam^{\infty}=\Gam^{|\Gam|}$.
\item {\rm An operator} $\Gam$ {\rm is said to be} $\Pi^{0}_{n}$ {\rm if} $\{(n,X)\in\ome\times{\cal P}(\ome): n\in\Gam(X)\}$ {\rm is defined by a} $\Pi^{0}_{n}${\rm -formula.}
\item {\rm For operators}  $\Gam_{0}, \Gam_{1}$ {\rm let}
\[
\\n\in[\Gam_{0},\Gam_{1}](X) :\Lrarw n\in\Gam_{0}(X) \, \vee\, [\Gam_{0}(X)\incl X\spand n\in\Gam_{1}(X)].
\]
\item {\rm Let} $\Phi_{0},\Phi_{1}$ {\rm be  classes of formulae.} {\rm An operator} $\Gam$ {\rm is said to be} $[\Phi_{0},\Phi_{1}]$ {\rm if} $\Gam=[\Gam_{0},\Gam_{1}]$ {\rm for some} $\Gam_{i}\in\Phi_{i}\, (i<2)$.
\item {\rm For} $n\in\Gam^{\infty}$ {\rm its} norm $| n|_{\Gam}$ {\rm (relative to} $\Gam${\rm ) is the least ordinal} $x$ {\rm such that} $n\in\Gam^{x+1}(\lrarw n\in\Gam(\Gam^{x}))$.
\eenu
\edf

Let L(PA) denote the language of Peano arithmetic. Variables in L(PA) are denoted by $n,m,\ldots$.

\bdf\label{df:Fixtheory}
{\rm Let} $\Phi$ {\rm be one of the classes} $\Pi^{0}_{n}, [\Pi^{0}_{n},\Pi^{0}_{m}]\, (n,m\in\ome)${\rm . Then} $\Phi\mbox{{\rm -Fix}}$ {\rm denote a first order two-sorted thoery defined as follows. Its} language 
{\rm L(Fix) is obtained from L(PA) by adding variables} $x,y,\ldots$ {\rm for} ordinals{\rm , and binary predicates} $x=y, x<y$ {\rm (less than relation on ordinals) and} binary predicate $n\in\Gam^{x}$ {\rm for each} $\Gam\in\Phi$.

 Axioms {\rm of the theory} $\Phi\mbox{{\rm -Fix}}$ {\rm are classified into four groups:}
\benu
\item
 Axioms {\rm of PA in the language L(Fix) and equality axioms for either sort.}
\item 
{\rm The} defining axioms {\rm for} $n\in\Gam^{x}${\rm :} $n\in\Gam^{x}\lrarw \exi y<x[n\in\Gam(\Gam^{y})]$.
\item
{\rm Axioms for the well ordering} $<$ {\rm on ordinals:}
$<$ {\rm is a linear ordering and} transfinite induction schema {\rm for} any {\rm formula} $F\in\mbox{{\rm L(Fix)}}${\rm :}
\[\fal x[\fal y<xF(y)\rarw F(x)]\rarw \fal xF(x).\]
\item 
Closure axiom: $\Gam(\Gam^{\infty})\incl\Gam^{\infty}$ {\rm for} $\Gam^{\infty}=\{n:\exi x(n\in\Gam^{x})\}$.
\eenu
\edf
{\bf General Conventions}.  Let $(X,<)$ be a quasiordering. Let $F$ be a function $F:X\ni\alp\mapsto F(\alp)\incl X$. For subsets $Y,Z\subset X$ of $X$ and elements $\alp,\bet\in X$, put
\benu
\item $\alp\leq\bet\Lrarw \alp<\bet \mbox{ or } \alp=\bet$.
\item $Y| \alp=\{\bet\in Y:\bet<\alp\}$.
\item $Y<Z  :\Lrarw  \exi\bet \in Z\fal \alp \in Y(\alp<\bet)$.
\item $Z\leq Y  :\Lrarw  \fal \bet \in Z\exi \alp \in Y(\bet\leq\alp)$.
\item $Y<\bet  :\Lrarw  Y<\{\bet\}$; $\alp< Z  :\Lrarw \{\alp\}< Z$; \\
$\bet\leq Y  :\Lrarw  \{\bet\}\leq Y$; $Z\leq \alp :\Lrarw Z\leq\{\alp\}$.
\item $F(Y)=\bigcup\{F(\alp):\alp\in Y\}$.
\item $[\alp,\bet]=\{\gam\in X:\alp\leq\gam\leq\bet\}$. Open intervals $(\alp,\bet)$ and half-open intervals $[\alp,\bet), (\alp,\bet]$ are defined similarly.
\item When $(X,<)$ is a linear ordering with its least element $0$ and $Y$ is a finite subset of $X$, $\max Y$ denotes the maximum of elements $\alp\in Y$ with respect to the ordering $<$. If $Y=\emptyset$, then, by convention, set $\max\emptyset:=0$.
\item Let $X^{<\ome}$ denote the set of finite sequences on $X$. Then $<_{lex}$ denotes the {\it lexicographic ordering\/} on $X^{<\ome}$ with $\fal s\in X^{n}\fal t\in X^{n+1}(s<_{lex}t)$.
If $<$ is a linear [well] ordering, then so is $<_{lex}$, resp.
\eenu

In each system of ordinal diagrams the {\it Veblen function\/} $\vphi_\alp(\bet)=\vphi\alp\bet$ is built-in as a constructor so that $\vphi0\bet=\ome^\bet$. Natural numbers are defined from this as usual and denoted by $i,j,k,l,m,n$

Now let us mention the content of this paper. 

In Section \ref{subsec:odsuper} we recall briefly the system $Od$ of ordinal diagrams (abbreviated by o.d.'s) in \cite{Wienpi3d}. The system $Od$ is a super-system of all systems of o.d.'s considered in this paper.

Let $Od^{\prime}$ be a subsystem of $Od$ which is closed under subdiagrams.
The Section \ref{sec:CX} is divied into two subsections. In the first subsection \ref{subsec:CX}, following Setzer \cite{Setzer98} and Buchholz \cite{BuchholzBSL}
\footnote{In the previous article \cite{Wienpi3d} , I attributed techniques on distinguished classes 
totally to Buchholz \cite{BuchholzBSL}, 
some of which were actually developed by Setzer \cite{Setzer98}.
I have learned the fact from the review \cite{SetzerMR} of \cite{Wienpi3d} written by  Setzer.}, 
we define sets $\calc^{\alp}(X)\incl Od^{\prime}$ for $\alp\in Od^{\prime}, X\incl Od^{\prime}$. 
In the second subsection \ref{subsec:proopr} we examine operators related to the sets. 
Almost all of these materials are reproduced from \cite{Wienpi3d}, and proofs are omitted.

In Section \ref{sec:wfidN} we introduce a system $Od(\Pi_{N})$ of o.d.'s for each positive integer $N\geq 4$.

The section is divided into seven subsections. 
The first four subsections are intended to give a set-theoretic interpretation of o.d.'s which is suggested by our wellfoundedness proof in \cite{odpiN} and in Section \ref{sec:5awf}.

The first subsection \ref{sec:prl5a} begins with defining iterations of Mahlo operations. These are intended to resolve 
a $\Pi_{N}$-reflecting universe. 
Since the resolving of $\Pi_{N}$-reflecting universes by iterations of $\Pi_{2}$-recursively Mahlo operations is 
so complicated, we first explain the simplest case in subsection \ref{subsec:M32}. 
Namely $\Pi_{3}$-reflecting universe $M_{3}^{2}$ on $\Pi_{3}$-reflecting ordinals. 
The subsection consists in three subsubsections. 
In the first one \ref{subsubsec:M32MC} we define Mahlo classes to resolve $M_{3}^{2}$ by iterations of 
$\Pi_{2}$-recursively Mahlo operations. 
In the second one \ref{subsubsec:odM32} we define the subsystem $Od(M_{3}^{2})\subset Od$ of ordinal diagrams for 
$M_{3}^{2}$. 
In the third one \ref{subsubsec:odM32wf} we prove the wellfoundedness of the countable fragment $Od(M_{3}^{2})|\Ome$ 
by means of a $[\Pi^{0}_{2},\Pi^{0}_{2}]$-inductive definition.

In subsection \ref{sec:prl5a.1} we define Mahlo classes for $\Pi_{N}$-reflection, and 
in subsection \ref{sec:prl5a.2} we associate Mahlo classes to o.d.'s in $Od(\Pi_{N})$.

In subsection \ref{subsec:od5api} we define the system $Od(\Pi_{N})$ of o.d.'s. 

In subsection \ref{subsec:od5fine} a finer analysis of relations $\prec_{i}$ on $Od(\Pi_{N})$ is given.
The analysis inspires us about the ramification procedure in subsections \ref{sec:prl5a.1} and \ref{sec:prl5a.2}, and
enables us to prove the wellfoundedness in Sections \ref{subsec:wfpiNid} and \ref{sec:5awf}.

In subsection \ref{subsec:daggareast} 
we introduce decompositions $\alp(s)$ of ordinal diagrams $\alp$, where 
$s$ denotes a function in ${}^{[i,k)}2\, (2\leq i\leq k\leq N-3)$.
In the next section \ref{subsec:wfpiNid} we define a suitable $\Pi^{0}_{N-1}$-operator $\Gam_{N}$ 
through the decompositions.

In Section \ref{subsec:wfpiNid} we prove the wellfoundedness of $Od(\Pi_{N})|\Ome$ by means of a $\Pi_{N-1}$-inductive 
definition as an extension of \cite{Wienpi3d}.

In subsection \ref{subsec:calgi}
we define operators $\calg_{i}\, (1\leq i<N-1)$ on $Od(\Pi_{N})$ recursively.
Using these operators, a $\Pi^{0}_{N-1}$-operator $\Gam_{N}$ is defined.
In subsection \ref{subsec:adequate} we show the adequacy of  the operator $\Gam_{N}$.
In subsection \ref{subsec:lem:id5wf21}
we conclude the proof.

In Section \ref{sec:5awf} we show that $\mbox{KP}\Pi_{N}$ does the same job, i.e., that for each $\alp<\Ome$ in $Od(\Pi_{N})$ $\mbox{KP}\Pi_{N}$ proves that 
$(Od(\Pi_{N})|\alp,<)$ is a well ordering. 
The wellfoundedness proof is based on the {\em distinguished class\/} (in German: 
Ausgezeichnete Klasse) and is an extension of ones in \cite{odMahlo}, \cite{odpi3}.
The proof is essentially the same given in \cite{odpiN}.
 
In the first subsection \ref{subsec:5awf.1} distinguished classes are defined and elementary facts on these classes are 
established. 
In the second subsection \ref{subsec:5awf.2} we introduce several classes of Mahlo universes and establish key facts 
on these classes. 
This is a crux in showing $Od(\Pi_{N})$ to be wellfounded without assuming the existence of the maximal distinguished class
 $\calw_{D}$, which is $\Sig^{1}_{2}$ on $\ome$ and hence a proper class in $\mbox{KP}\Pi_{N}$. 
These classes imitate the ramification procedure described in subsections \ref{sec:prl5a.1} and \ref{sec:prl5a.2}. 
In the third subsection \ref{subsec:5awf.3} we conclude a proof of wellfoundedness of $Od(\Pi_{N})$.

We rely on the previous \cite{Wienpi3d}, and state some lemmata without proofs since we gave proofs of these in \cite{Wienpi3d}

\section{The system $Od$}\label{subsec:odsuper}

In this section let us recall briefly the system $Od$ of ordinal diagrams (abbreviated by o.d.'s) in \cite{Wienpi3d}.
 The system $Od$ is a super-system of all systems of o.d.'s considered in this paper, and each of them is obtained by posing restrictions on the construction $(\sig,\alp,q)\mapsto d_{\sig}^{q}\alp$ in Definition \ref{df:opi}.\ref{caldQ}.

Let $0,\Ome,\pi,+,\vphi, {}^{+}$ and $d$ be distinct symbols. Each o.d. in the system $Od$ is a finite sequence of these symbols. $\vphi$ is the binary Veblen function. $\Ome$ denotes the first recursively regular ordinal $\ome^{CK}_{1}$. 
$\sig^{+}$ denotes the next recursively regular ordinal to $\sig$. 

$\ell\alp$ denotes the number of occurrences of symbols in the o.d. $\alp$. Let $sd(\alp)$ denote the {\it set of proper subdiagrams\/} (subterms) of $\alp$. Thus $\alp\not\in sd(\alp)$. Also put $sd^+(\alp)=sd(\alp)\cup\{\alp\}$.

The set $Od$ is classified into subsets $R=\{\pi\}\dot{\cup}\cald^{Q}\dot{\cup}SR, SC, P$ according to the intended meanings of o.d.'s. $P$ denotes the set of additive principal numbers, $SC$ the set of strongly critical numbers and $R$ the set of recursively regular ordinals.
Ordinal diagrams are denoted $\alp,\bet,\gam,\ldots$, while $\sig,\tau,\ldots$ denote o.d.'s in the set $R$.

Let us reproduce generating clauses of o.d.'s not in \cite{odMahlo}.
 
\bdf\label{df:opi} $Od$.
 \begin{enumerate}
 \item\label{df:opi.6} $\sig \in\cald^{Q} \spand 0<k<\ome \Rarw \sig^{+k} \in SR$.\footnote{In proof-theoretic studies, i.e., cut-elimination in \cite{ptMahlo}, \cite{ptpi3} and \cite{ptpiN} the construction $\sig\mapsto\sig^{+}$ is not needed. The construction helps us to define the system $Od(\Pi_{N})$ smoothly.}
 \item\label{df:opi.7} {\rm Let} $\alp \in Od \, \& \,  \sig\in\{\pi\}\cup SR=\{\Ome\}\cup\{\kap^{+k}:\kap\in\cald^{Q}, k>0\}${\rm . Put} $\eta:=d_{\sig}\alp$ {\rm and define}
  \[b(\eta)=\alp, Q(\eta)=\emptyset, c(\eta)=\{\alp\}\cup Q(\eta)=\{\alp\}.\]
  {\rm Assume that the following condition is fulfilled:}

\beqn\label{eq:Odmu}
\calb_{>\sig}(\{\sig\}\cup c(\eta))<b(\eta)=\alp 
\eeqn

 {\rm Then} $\eta=d_{\sig}\alp\in SC$.
 \item\label{caldQ} {\rm Let} $\alp \in Od \, \& \,  \sig\in\{\pi\}\cup\cald^{Q} \spand q=\ovl{j\kap\tau\nu}\incl Od$, {\rm where} $q=\ovl{j\kap\tau\nu}$ {\rm denotes a non-empty sequence of} quadruples $j_{m}\kap_m\tau_m\nu_m$ {\rm of length} $l+1 \, (l\geq 0)${\rm . Put} $\eta:=d^q_\sig\alp\in\cald^{Q}_{\sig}$ {\rm and define}
  \[b(\eta)=\alp, Q(\eta)=q=\{j_{m}, \kap_m, \tau_m, \nu_m:m\leq l \}, c(\eta)=\{\alp\}\cup Q(\eta).\]
  {\rm Assume that the condition (\ref{eq:Odmu}) is fulfilled.}

{\rm Then} $\eta=d^q_\sig\alp\in\cald^{Q}_{\sig}$.
 \end{enumerate}
\edf

$\cald_{\sig}$ denotes the set of diagrams of the form $d_{\sig}^{q}\alp$.
$\alp\prec\bet$ denotes the transitive closure of the relation $\{(\alp,\bet):\alp\in\cald_{\bet}\}$, and $\rho\prec\sig\Lrarw \tau\preceq\sig$ for $\rho\in\cald_{\tau}$. 
For any $\alp\in Od$, set $\alp<\pi^{+}=\infty$ where $\pi^{+}=\infty$ denotes an extra symbol not in $Od$.

Finite subsets of subdiagrams of o.d.'s $\alp$, $K^{d}\alp\incl\cald\cap sd^{+}(\alp)$, $K\alp\incl SC\cap sd^{+}(\alp)$, $\calb_{\sig}(\alp), \calb_{>\sig}(\alp)\incl sd(\alp)$ and
$K_{\sig}\alp\incl sd^{+}(\alp)\cap\cald$ are defined as in \cite{odMahlo}, \cite{Wienpi3d}.
Specifically for $\alp\in\cald_{\tau}$
\[K_{\sig}\alp=\left\{
        \begin{array}{ll}
         K_{\sig}(\{\tau\}\cup c(\alp)), & \mbox{{\rm if }} \sig<\tau \\
         K_{\sig}\tau, & \mbox{{\rm if }} \tau<\sig \, \& \, \tau\not\preceq\sig \\
         \{\alp\}, & \mbox{{\rm if }} \alp\prec\sig
        \end{array}
        \right.\]

\blem\label{lem:3.2}
\benu
\item $\alp\prec\bet\Rarw \ell\bet<\ell\alp\,\&\, \alp<\bet$.
\label{lem:3.2.2}
\item If $\bet\in K_\sig\alp$, then $\bet\prec\sig$, $\bet$ is a subdiagram of $\alp$ and 
$\sig$ is a proper subdiagram of $\alp$.
\label{lem:3.2.1}
\item $\alp<\sig\,\&\, K_\sig\alp<\bet\prec\sig \Rarw \alp<\bet$.
\label{lem:3.2.5}
\item $K_\sig\alp\leq\alp$.
\label{lem:3.2.9}
\item $\alp\in\cald_{\sig} \spand \kap<\sig \Rarw K_{\kap}\alp<\alp$, and 
$\alp\preceq\tau\in\cald_\sig \,\&\, \kap<\sig\Rarw K_{\kap}\tau<\alp$.
\label{lem:3.2.10}
\eenu
\elem

As in \cite{ptMahlo} we see the following lemma.

\blem\label{lem:Od3}
For $\alp,\bet,\tau\in Od$, $\alp\leq\bet<\tau \Rarw \calb_{\tau}(\alp)\leq\calb_{\tau}(\bet)$.
\elem

\blem\label{lem:Npi11exist}
For $\alp,\bet\in Od$, $\alp\prec\bet\in\cald \Rarw b(\bet)<b(\alp)$.
\elem

\section{Sets $\calc^{\alp}(X)$ and operators}\label{sec:CX}

Let $Od^{\prime}$ be a subsystem of $Od$ which is closed under subdiagrams: $\alp\in Od^{\prime} \Rarw sd(\alp)\incl Od^{\prime}$. $X,Y,\ldots$ ranges over subsets of $Od^{\prime}$.
In this section we define sets $\calc^{\alp}(X)\incl Od^{\prime}$ for $\alp\in Od^{\prime}, X\incl Od^{\prime}$, and examine operators related to the sets.

Almost everything in this section is reproduced from \cite{Wienpi3d}, and omitted proofs of lemmata can be found there.

\subsection{The sets $\calc^{\alp}(X)$}\label{subsec:CX}

Following Setzer \cite{Setzer98} and Buchholz \cite{BuchholzBSL}
 we define sets $\calc^{\alp}(X)\incl Od^{\prime}$ for $\alp\in Od^{\prime}, X\incl Od^{\prime}$ as follows.

\bdf\label{df:CX}
{\rm For} $\alp\in Od^{\prime}, X\incl Od^{\prime}${\rm , let} 
\beqnarr
\calc^{\alp}(X) & := & 
\mbox{{\rm closure of }} \{0,\Ome,\pi\}\cup(X|\alp) \mbox{ {\rm under} } +,\vphi, \kap\mapsto\kap^{+} 
\nonumber \\
&& \mbox{ {\rm and} } (\sig,\alp,q)\mapsto d_{\sig}^{q}\alp \mbox{ {\rm for} } \sig>\alp \mbox{ {\rm in }} Od^{\prime}
\label{eq:CX}
\eeqnarr
\edf

\blem\label{lem:CX1}
$X|\alp=Y|\alp \Rarw \calc^{\alp}(X)=\calc^{\alp}(Y)$ and $X\mapsto\calc^{\alp}(X)$ is monotonic.
\elem

\bdf\label{df:CXcond}
 {\rm Consider the following conditions for} $X\incl Od^{\prime}$ {\rm :}
\bdes
\item[(A)] $\fal\alp\in X[\alp\in\calc^{\alp}(X)]$.
\item[(K)] $\fal\alp\in X\fal\sig[K_{\sig}\alp\incl X]$.
\item[(KC)] $\fal\alp\fal\bet\fal\sig[\alp\in\calc^{\bet}(X)\spand\sig\leq\bet \Rarw K_{\sig}\alp\incl X]$.
\edes
\edf

\blem\label{lem:KC} Assume $X$ enjoys the condition {\bf (K)}. Then $X$  enjoys the condition {\bf (KC)}, too.
\elem

\blem\label{lem:CX2}
Assume $X\incl Od^{\prime}$ enjoys the condition {\bf (A)}.
\benu
\item\label{lem:CX2.3} $\alp\leq\bet \Rarw \calc^{\bet}(X)\incl\calc^{\alp}(X)$.
\item\label{lem:CX2.4} $\alp<\bet<\alp^{+} \Rarw \calc^{\bet}(X)=\calc^{\alp}(X)$, where
$\alp^{+}:=\min\{\sig\in R\cup\{\infty\}:\alp<\sig\}$.
\eenu
\elem

\blem\label{lem:CX3} Assume $\gam\in\calc^{\alp}(X)$, $\alp<\bet$ and $\fal\kap\leq\bet[K_{\kap}\gam<\alp]$.
\benu
\item\label{lem:CX3.1} Assume $\mbox{{\rm LIH:} } \fal\del[\ell\del\leq\ell\gam \spand \del\in\calc^{\alp}(X)|\alp \Rarw \del\in\calc^{\bet}(X)\supseteq X|\alp]$. Then $\gam\in\calc^{\bet}(X)$.
\item\label{lem:CX3.2} $\calc^{\alp}(X)|\alp \incl X \Rarw \gam\in\calc^{\bet}(X)$.
\eenu
\elem

\bdf\label{df:persistent}
{\rm An operator} $\Gam$ {\rm on} $Od^{\prime}${\rm , i.e.,} $\Gam:{\cal P}(Od^{\prime})\rarw{\cal P}(Od^{\prime})$ {\rm is said to be} persistent {\rm if} 
$\fal\alp\in Od^{\prime}\fal X, Y\incl Od^{\prime}[X|\alp=Y|\alp \Rarw \Gam(X)|(\alp+1)=\Gam(Y)|(\alp+1)]$.
\edf

\bdf\label{df:G()}
\benu
\item\label{df:G().1}
$\calg(X):=\{\alp: \alp\in\calc^{\alp}(X)\spand \calc^{\alp}(X)|\alp\incl X\}$.
\item\label{df:G().1.5}
$R^{\prime}:=\{\alp\in Od^{\prime}: \cald_{\alp}\cap Od^{\prime}\neq\emptyset\}\incl R$.
\item\label{df:G().2}
$\alp\in \Gam_{2}(X) :\Lrarw \alp<\pi \,\wedge\, \alp\not\in R^{\prime}\,\wedge\, \alp\in\calg(X)$.
\eenu
\edf

\blem\label{lem:CXG}
The operators $\calg$ and $\Gam_{2}$ are $\Pi^{0}_{1}$ and persistent.
\elem

\blem\label{lem:CX4}
Assume $\alp\in\calc^{\alp}(X)$ and $\alp\preceq\sig$.
\benu
\item\label{lem:CX4.0} $\sig\in\calc^{\alp}(X)$.
\item\label{lem:CX4.1} If $\alp\in\calg(X)$, then $\sig\in\calc^{\bet}(X)$ for any $\bet$ with $\alp\leq\bet\leq\sig$.
\item\label{lem:CX4.2} If $\alp\in\calg(X)$ and $\sig\in\cald_{\tau}$, then $\sig\in\calc^{\bet}(X)$ for any $\bet$ with $\alp\leq\bet<\tau$.
\eenu
\elem

\blem\label{lem:CX4aro}
Assume $\bigcup\{K_{\sig}\nu: \sig\leq \kap\}\incl X|\rho$. Then 
$\nu\in\calc^{\kap}(X)$.
\elem
\bprf
We show, by induction on $\ell\gam$,
\[
\fal\gam\in sd^{+}(\nu)[\bigcup\{K_{\sig}\gam: \sig\leq\kap\}\incl X|\rho \Rarw \gam\in\calc^{\kap}(X)]
\]
for the set of subdiagram $sd^{+}(\nu)$ of $\nu$. 

If $\gam\not\in\cald$, then IH yields $\gam\in\calc^{\kap}(X)$. 
Suppose $\gam\in\cald_{\tau} \spand \bigcup\{K_{\sig}\gam: \sig\leq\kap\}\incl X|\rho$ with some $\{\tau\}\cup c(\gam)$.

If $\tau>\kap$, then IH yields $\{\tau\}\cup c(\gam)\incl \calc^{\kap}(X)$, and hence we are done.

Suppose $\tau\leq\kap$.
By $\bigcup\{K_{\sig}\gam: \sig\leq\kap\}\incl X|\rho$ we have 
$\{\gam\}=K_{\tau}\gam\incl X|\rho\incl X|\kap\incl\calc^{\kap}(X)$.
\eprf

\blem\label{lem:CX6}
Assume $X$ enjoys the conditions {\bf (A)} and {\bf (K)}. Further assume $\alp\in\calg(X)$ and $\alp\prec\sig$. Then either $\exi\del\in X[\alp\leq\del\prec\sig]$ or $\sig\in\calg(X)$.
\elem
\bprf
Suppose $\neg\exi\del\in X[\alp\leq\del\prec\sig]$. We show $\sig\in\calg(X)$. $\sig\in\calc^{\sig}(X)$ follows from Lemma \ref{lem:CX4}.\ref{lem:CX4.1}. It remains to show
$\gam\in\calc^{\sig}(X)|\sig \Rarw \gam\in X$. 
\\
{\bf Case1} $\gam<\alp$: By Lemma \ref{lem:CX2}.\ref{lem:CX2.3} we have $\gam\in\calc^{\sig}(X)|\alp\incl\calc^{\alp}(X)|\alp\incl X$.
\\
{\bf Case2} $\gam=\alp$: By $\gam=\alp\prec\sig$ we can assume $\gam\in\cald_{\kap}$ for a $\kap>\sig$ with $\{\kap\}\cup c(\gam)\incl\calc^{\sig}(X)$. Then we would have $\kap\preceq\sig<\kap$. This is not the case.
\\
{\bf Case3} $\gam>\alp$: Then $\alp<\gam<\sig$. Lemma \ref{lem:3.2}.\ref{lem:3.2.5} yields $\alp\leq K_{\sig}\gam$. 

On the other side, $X$ enjoys the condition {\bf (KC)} by Lemma \ref{lem:KC}. Thus $K_{\sig}\gam\incl X$. However $\del\prec\sig$ for any $\del\in K_{\sig}\gam$ by Lemma \ref{lem:3.2}.\ref{lem:3.2.1}. 
Therefore we would have $\exi\del\in X[\alp\leq\del\prec\sig]$.
\eprf

\subsection{Families of sets in wellfoundedness proofs}\label{subsec:proopr}
Let $\tht[X]$ denote a (definable) property on subsets $X$ of $Od^{\prime}$.
In this subsection some conditions on families of sets are extracted from wellfoundedness proofs, and we derive properties of sets under the conditions.

\bdf\label{df:calgtht}
\benu
\item $\calw_{\tht}:=\bigcup\{X: \tht[X]\}$.
\item
$\calg_{\tht}:=\calg(\calw_{\tht})|\pi$.
\eenu
\edf

\bdf\label{df:wftg}
\benu
\item $Prg[X,Y] :\Lrarw \fal\alp(X|\alp\subseteq Y \spand \alp\in X \rarw \alp\in Y)$.

\item {\rm For a definable class} $\calx$, $TI[\calx]$ {\rm denotes the schema:}\\
$TI[\calx] :\Lrarw Prg[\calx,\caly]\rarw \calx\incl\caly \mbox{ {\rm holds for} any definable class } \caly$.
\item
{\rm For} $X\incl Od^{\prime}$ {\rm let} $WX$ {\rm denote the} wellfounded part {\rm of} $X$. $WX$ {\rm is defined from the monotonic} $\Pi^{0}_{1}${\rm -operator (in} $X${\rm )} $\Gam_{A}(Y):=\{\alp\in X : \fal\bet\in X|\alp(\bet\in Y)\}$.
\eenu
\edf

Consider the following conditions on the property $\tht$:

\bdes
\item[($\tht$.0)] $\calw_{\tht}<\pi$.
\item[($\tht$.1)] $\fal \alp\in X\exi Y\incl X\{\tht[X] \Rarw \tht[Y] \spand \alp\in\calg(Y) \spand X|\alp=Y|\alp\}$.
\item[($\tht$.2)] For any $X$ with $\tht[X]$, $TI[X]$ and $\alp\in X \Rarw \calw_{\tht}|\alp=X|\alp$.

\item[($\tht$.3)] $\Gam_{2}(\calw_{\tht})\incl\calw_{\tht}$, i.e.,
$\fal\alp<\pi[\alp\not\in R^{\prime}\spand\alp\in\calg_{\tht} \Rarw \alp\in\calw_{\tht}]$.
\item[($\tht$.4)] $\fal\alp[\alp\in SR\spand\alp\in\calg_{\tht} \Rarw \alp\in\calw_{\tht}]$.
\edes

\subsubsection{Elementary properties of the family of sets}\label{subsubsect:elemprop}

\blem\label{lem:id7wf8}
Assume that $\tht$ enjoys hypotheses {\bf ($\tht$.i)} for $i\leq 1$.
For any $X\in\{X:\tht[X]\}\cup\{\calw_{\tht}\}$,
\benu
\item\label{lem:wel.3}
$\fal\alp\in X[\alp\in\calc^{\alp}(X)]$, and hence $X$ enjoys the condition {\bf (A)} in Definition \ref{df:CXcond}.
\item $\fal\tau[\alp\in X\Rarw K_{\tau}\alp\incl X]$.
\label{lem:id7wf8.3}
\item $\fal\bet\fal\tau[\alp\in\calc^{\bet}(X) \Rarw K_{\tau}\alp\incl\calc^{\bet}(X)]$.
\label{lem:id7wf8.3.5}
\eenu

Hence $X$ enjoys the conditions {\bf (K)} and {\bf (KC)} in Definition \ref{df:CXcond}.
\elem

\blem\label{lem:CX6+}
Assume that $\tht$ enjoys hypotheses {\bf ($\tht$.i)} for $i\leq 1$.
For any $X\in\{X:\tht[X]\}\cup\{\calw_{\tht}\}$,
if $\alp\in\calg(X)$ and $\alp\prec\sig$, then either $\exi\del\in X[\alp\leq\del\prec\sig]$ or $\sig\in\calg(X)$.
\elem
\bprf
This follows from Lemma \ref{lem:CX6} with Lemmata \ref{lem:id7wf8}.\ref{lem:wel.3} and \ref{lem:id7wf8}.\ref{lem:id7wf8.3}.
\eprf

\subsubsection{Hypotheses on wellfoundedness}\label{subsubsec:hyp3-5}

\blem\label{cor:A1}
Assume that $\tht$ enjoys hypotheses {\bf ($\tht$.i)} for $i\leq 2$.
\benu
\item\label{cor:A1.3}
$TI[\calw_{\tht}]$.
\item\label{cor:A1.3.1}
For any $X\in\{X:\tht[X]\}\cup\{\calw_{\tht}\}$,
$X\incl\calg(X)$.
Hence $\alp\in\calc^{\bet}(X)|\bet \spand \bet\in X \Rarw \bet\in X$.
\item\label{cor:A1.4}
Let $X\in\{X:\tht[X]\}\cup\{\calw_{\tht}\}$ and assume {\bf ($\tht$.3)} $\Gam_{2}(X)\incl X$.
Then
$\fal\alp<\pi[K\alp\incl X \Rarw \alp\in  X]$.
\eenu
\elem

\blem\label{lem:WWOme}
Assume that $\tht$ enjoys hypotheses {\bf ($\tht$.i)} for $i\leq 3$.
Then 
\[WOd^{\prime}|\Ome=\calw_{\tht}|\Ome.\]
\elem

\bdf\label{df:3awp.1}
{\rm Set}
\[\calw_{\pi}=\calc^{\pi}(\calw_{\tht}).\]
\edf

\blem\label{lem:id3wf19-1}
Assume that $\tht$ enjoys hypotheses {\bf ($\tht$.i)} for $i\leq 4$.
\benu
\item $\calw_{\tht}\incl\calw_{\pi}$ and $\calw_{\pi}|\pi=\calw_{\tht}|\pi$.
\label{lem:id3wf19-1.1}
\item $\pi\in\calw_{\pi}$.
\label{lem:id3wf19-1.2}
\item
$\alp\in\calw_{\pi} \Lrarw K^{d}\alp\incl\calw_{\pi} \Lrarw K^{d}\alp\incl\calw_{\tht}$.
\label{lem:id3wf19-1.3}
\item $\alp\in\calw_{\pi}\Rarw K_{\sig}\alp\incl\calw_{\tht}\spand \alp\in\calc^{\bet}(\calw_{\tht})$ for any $\sig$ and $\bet$.
\label{lem:id3wf19-1.3.5}
\item For each $n\in\ome,$ $TI[\calw_{\pi}|\ome_{n}(\pi+1)]$ with $\ome_{0}(\alp)=\alp, \ome_{n+1}(\alp)=\ome^{\ome_{n}(\alp)}$.
\label{lem:id3wf19-1.4}
\eenu
\elem

\bdf\label{df:id4wfA}{\rm For o.d.'s} $\alp\in Od^{\prime}$ {\rm and finite sequences of quadruples} $q\incl Od^{\prime}${\rm , define:}
\benu
\item\label{df:id4wfA.1}
\[
 A(\alp,q) :\Lrarw 
 \fal\sig\in\calw_{\pi}\fal\alp_1\in \cald_{\sig}[\alp=b(\alp_1) \spand q=Q(\alp_{1}) \Rarw \alp_1\in\calw_{\tht}].
\]
\item \label{df:id4wfA.2}
\[
\mbox{{\rm MIH}}(\alp) :\Lrarw
 \fal\bet\in\calw_{\pi}|\alp\fal q\incl \calw_{\pi} A(\bet,q).
\]
\item\label{df:id4wfA.3}
\[
\mbox{{\rm SIH}}(\alp,q) :\Lrarw
 \fal q_{0}\incl \calw_{\pi}[q_{0}<_{lex}q \Rarw A(\alp,q_{0})],
\]
{\rm where sequences of quadruples} $q=(j_{m}, \kap_m, \tau_m, \nu_m:m\leq l)$ {\rm are arranged in the ordering:} $q=(j_{l}, \kap_{l}, \tau_{l}, \nu_{l},\ldots, j_{0}, \kap_{0}, \tau_{0}, \nu_{0})$.
\eenu
\edf

\blem\label{lem:id3wf20}
Assume that $\tht$ enjoys hypotheses {\bf ($\tht$.i)} for $i\leq 4$.
Suppose $\mbox{{\rm MIH}}(\alp)$.
 For any o.d. $\bet\in Od^{\prime}$ 
\[
\bet\in\calc^{\sig}(\calw_{\tht}) \spand \calb_{>\sig}(\bet)<\alp 
 \Rarw \bet\in\calw_{\pi}.
\]
\elem

\bth\label{th:id5wf21}
Assume that $\tht$ enjoys hypotheses {\bf ($\tht$.i)} for $i\leq 4$.
Assume $\{\alp\}\cup q\incl \calw_{\pi}$, $\mbox{{\rm MIH}}(\alp)$, and $\mbox{{\rm SIH}}(\alp,q)$ in Definition \ref{df:id4wfA}.
Then
\[
 \fal\sig\in\calw_{\pi}\fal\alp_1\in \cald_{\sig}[\alp=b(\alp_1) \spand q=Q(\alp_{1}) \Rarw \alp_1\in\calg_{\tht}].
\]
\eth

\subsection{Hypotheses on operators}\label{subsec:ensureoperators}

Let $\Gam$ denote a first order operator on a subsystem $Od^{\prime}$ of $Od$.
Let $\tht[X] :\Lrarw \exi x\in Ord[X=\Gam^{x}]$. Thus

\bdf\label{df:calg1}
\benu
\item $\calw:=\bigcup\{\Gam^{x}:x\in Ord\}$.
\item
$\calg:=\calg(\calw)|\pi$.
\item
$|\alp|:=|\alp|_{\Gam}$ {\rm for} $\alp\in\calw$.
\eenu
\edf

Consider the following conditions on $\Gam$ for any $X\in\{\Gam^{x}:x\in Ord\}\cup\{\calw\}$:
\bdes
\item[($\Gam$.0)] $\Gam(X)<\pi$.
\item[($\Gam$.1)] $\Gam(X)\incl\calg(X)$.
\item[($\Gam$.2)] $\alp,\bet\in\calw \spand \alp<\bet \Rarw | \alp|<|\bet|$.
\item[($\Gam$.3)] $\Gam_{2}(X)\incl\Gam(X)$, i.e., 
$\fal\alp<\pi[\alp\not\in R^{\prime}\spand\alp\in\calg(X) \Rarw \alp\in\Gam(X)]$.
\item[($\Gam$.4)] $\fal\alp[\alp\in SR\spand\alp\in\calg \Rarw \alp\in\calw]$.
\item[($\Gam$.5)] $\Gam(\calw)\incl\calw$.
\edes

\blem\label{lem:welGam}
Assume $\Gam$ enjoys the hypotheses {\bf ($\Gam$.i)} for $i\leq 2$.
Then $\tht[X] :\Lrarw \exi x\in Ord[X=\Gam^{x}]$ enjoys {\bf ($\tht$.i)} for $i\leq 2$, and $X\incl\calg(X)$ for any $X\in\{\Gam^{x}:x\in Ord\}$.

Furthermore if $\Gam$ enjoys {\bf ($\Gam$.i)} for $3\leq i\leq 5$, then {\bf ($\tht$.i)} for $3\leq i\leq 4$ holds.
\elem

In general for any persistent operator $\Gam$ satisfying {\bf ($\Gam$.2)}, $\calw\incl\Gam(\calw)$ holds.

\blem\label{lem:A1.2}

If $\Gam$ is persistent and enjoys {\bf ($\Gam$.2)}, then $\calw\incl\Gam(\calw)$.
\elem

\section{$\Pi_{N}$-reflection}
\label{sec:wfidN}

In this section we introduce a recursive notation system $Od(\Pi_{N})$ of ordinals for each positive integer $N\geq 4$, 
which we studied first in \cite{hndodpiN}. 

Let $\mbox{KP}\Pi_{N}$ denote the set theory for $\Pi_{N}$-reflecting universes. $\mbox{KP}\Pi_{N}$ is obtained from 
the Kripke-Platek set theory with the Axiom of Infinity by adding the axiom: for any $\Pi_{N}$ formula $A(u)$, 
$A(u)\rarw \exi z(u\in z\spand A^{z}(u))$, where $A^{z}$ denotes 
the result of restricting any unbounded quantifiers $Q x\, (Q\in\{\exi,\fal\})$ in $A$ to $Q x\in z$.

We show that for each $\alp<\Ome$ in $Od(\Pi_{N})$, both $\mbox{KP}\Pi_{N}$ and $\Pi_{N-1}\mbox{-Fix}$ prove that the initial segment of $Od(\Pi_{N})$ determined by $\alp$ is a well ordering, where $\Ome$ denotes the first recursively regular ordinal $\ome^{CK}_{1}$ and $\pi$ the first $\Pi_{N}$-reflecting ordinal. Each $\alp\in Od(\Pi_{N})$ is less than the next epsilon number $\veps_{\pi+1}$ to $\pi$.

\subsection{Prelude}\label{sec:prl5a}

The main constructor in $Od(\Pi_{N})$ is to form an o.d. $d_{\sig}^{q}\alp<\sig$ from a symbol $d$ and o.d.'s $\sig,q,\alp$,
 where $\sig$ denotes a recursively regular ordinal and $q$ a finite sequence of quadruples of o.d.'s. 
By definition we set $d_{\sig}^{q}\alp<\sig$. Let $\gam\prec_{2}\sig:\Lrarw \gam\prec\sig$,
 and $\preceq_{2}$ its reflexive closure. 
Then the set $\{\tau:\sig\prec_{2}\tau\}$ is finite and linearly ordered by $\prec_{2}$ for each $\sig$, 
namely $\{\sig:\sig\preceq_{2}\pi\}$ is a tree with the root $\pi$. 
In the diagram $d_{\sig}^{q}\alp$ $q$ includes some data telling us how the diagram $d_{\sig}^{q}\alp$ is constructed from 
its {\it predecessors\/} $\{\tau:d_{\sig}^{q}\alp\prec_{2}\tau\}=\{\tau:\sig\preceq_{2}\tau\}$. 
Here involves subtle and complicated requirements to which $d_{\sig}^{q}\alp$ have to obey, cf. Definition \ref{df:piN}. 
These were obtained solely from finitary analysis of finite proof figures for $\Pi_{N}$-reflection, cf. \cite{ptpiN}: 
its generation has not referred to any set-theoretic considerations. 
Despite the lack of meaning it now turns out that our wellfoundedness proof of $Od(\Pi_{N})$ in \cite{odpiN} and in 
Section \ref{sec:5awf}, 
which is formalizable in KP$\Pi_{N}$, suggests a set-theoretic interpretation. Let us explain this.
 
In a wellfoundedness proof for $Od(\Pi_{N})$ using the maximal distinguished class $\calw_{D}$ introduced 
in \cite{Buchholz75}, the main task is to show the tree $\{\sig:\sig\preceq_{2}\pi\}$ to be wellfounded. 
When we assume the existence of a $\Sig^{1}_{2}$-class, i.e., $\calw_{D}$ as a {\it set\/}, then we \cite{hndodpiN} 
can show that $Od(\Pi_{N})$ is wellfounded. 
Nevertheless $\calw_{D}$ is a proper {\it class\/} in $\mbox{KP}\Pi_{N}$. 
Therefore we have to show for each o.d. $\eta$ there exists a {\it set\/}, say $P$, in which we can imitate constructions 
in \cite{hndodpiN} up to the given $\eta$. 
Namely the maximal distinguished class defined {\it on\/} $P$ denoted $\calw^{P}$ has to enjoy the same closure properties
 as $\calw_{D}$ up to $\eta$. 
Such a set $P$ is said to be $\eta${\it -Mahlo\/}. 
Then the existence of $\eta$-Mahlo sets guarantees the wellfoundedness of the chain $\{\tau:\tau\preceq_{2}\eta\}$ 
with respect to $\prec_{2}$. 
Thus a crux in showing $Od(\Pi_{N})$ to be wellfounded without assuming the existence of a $\Sig^{1}_{2}$-class $\calw_{D}$
 is to show the existence of $\eta$-Mahlo sets for each $\eta$.
We have learnt in \cite{odpi3} that if a set is $\Pi_{2}$-reflecting on $\gam$-Mahlo universes for any
 $\gam\prec_{2}\eta$, then the set is $\eta$-Mahlo.

Let L denote a $\Pi_{N}$-refecting universe: $(\mbox{L};\in)\models \mbox{KP}\Pi_{N}$. Transitive sets in $\mbox{L}\cup\{\mbox{L}\}$ are denoted $P,Q,\ldots$, and $\mbox{L}^{t}$ denotes the set of transitive sets in L. For a transitive set $P$ let $ord(P)$ denote the set of ordinals in $P$. Also let $ord(P)^{+}$ denote the supremum of $\Del_{1}$-wellfounded relations on $P$. A class $\calx\incl \mbox{L}^{t}\cup\{\mbox{L}\}$ is said to be a $\Pi_{i}${\it -class\/} if there exists a set-theoretic $\Pi_{i}$-formula $F$ such that for any $P\in \mbox{L}^{t}\cup\{\mbox{L}\}$, $P\in\calx\Lrarw P\models F :\Lrarw (P;\in)\models F$. 
By a $\Pi^{1}_{0}${\it -class} we mean a $\Pi_{i}$-class for some $i\geq 2$.
A sequence $\{\calx_{\xi}\}_{\xi<\alp}\, (\calx_{\xi}\incl \mbox{L}^{t}\cup\{\mbox{L}\})$ is said to be a $\Pi_{i}${\it -sequence\/} if classes $\calx_{\xi}$ are $\Pi_{i}$-classes uniformly, i.e., there exists a set-theoretic $\Pi_{i}$-formula $F(\xi)$ such that for any $P\in \mbox{L}^{t}\cup\{\mbox{L}\}$ and any $\xi<\min\{\alp,ord(P)^{+}\}$, $P\in\calx_{\xi}\Lrarw P\models F(\xi)$. 
By a $\Pi^{1}_{0}${\it -sequence} we mean a $\Pi_{i}$-sequence for some $i\geq 2$.

A $\Pi_{i}$-recursively Mahlo operation is defined through a universal $\Pi_{i}$-formula $\Pi_{i}(a)$, cf. \cite{Richter-Aczel74} as follows
\[P\in M_{i}(\calx) :\Lrarw \fal b\in P[P\models\Pi_{i}(b)\rarw \exi Q\in \calx\cap P(Q\models\Pi_{i}(b))].\]

We consider two kinds of its iterations, both are defined by transfinite recursion on ordinals $\bet$ (or, in general, along a well founded relation), and their mixture.
First an inner iteration is defined from a sequence $\{\calx_{\xi}\}_{\xi<\alp}$: let
\[M_{1,i}^{\bet}(\{\calx_{\xi}\}_{\xi<\alp}):=
\bigcap\{M_{i}(M_{i}(M_{1,i}^{\nu}(\{\calx_{\xi}\}_{\xi<\alp})\cap\calx_{\xi})\cap\calx_{\del}):
\del\leq\xi<\alp, \nu<\bet\}.\]
Second an outer iterartion is defined from a class $\calx$:
\[M_{2,i}^{\bet}(\calx):=
\calx\cap\bigcap\{M_{i}(M_{2,i}^{\nu}(\calx)):\nu<\bet\}.\]
For the case $\calx=\mbox{L}^{t}\cup\{\mbox{L}\}$ put 
\[M_{i}^{\bet}(\mbox{L}):=M_{2,i}^{\bet}(\mbox{L}^{t}\cup\{\mbox{L}\})=
\bigcap\{M_{i}(M_{i}^{\nu}(\mbox{L})):\nu<\bet\}.\]

Finally their mixture is defined as follows:
\beqnarrs
&& M_{i}^{\bet}(\calx;\{\calx_{\xi}\}_{\xi<\alp}):= \\
&&\calx\cap
\bigcap\{M_{i}(M_{i}(M_{i}^{\nu}(\calx;\{\calx_{\xi}\}_{\xi<\alp})\cap\calx_{\xi})\cap\calx_{\del}):
\del\leq\xi<\alp, \nu<\bet\}.
\eeqnarrs

Obviously $M_{i}^{\bet}(\mbox{L}^{t}\cup\{\mbox{L}\};\{\calx_{\xi}\}_{\xi<\alp})=M_{1,i}^{\bet}(\{\calx_{\xi}\}_{\xi<\alp})$.

Observe that $M_{i}(\calx)$ is a $\Pi_{i+1}$-class if $\calx$ is a $\Pi^{1}_{0}$-class. $\{M_{1,i}^{\nu}(\{\calx_{\xi}\}_{\xi<\alp})\}_{\nu<\bet}$ is a $\Pi_{i+1}$-sequence for any $\Pi^{1}_{0}$-sequence $\{\calx_{\xi}\}_{\xi<\alp}$, and $\{M_{2,i}^{\nu}(\calx):\nu<\bet\}$ is a $\Pi_{i+1}$-sequence for any $\Pi_{i+1}$-class $\calx$. Finally $\{M_{i}^{\nu}(\calx;\{\calx_{\xi}\}_{\xi<\alp})\}_{\nu<\bet}$ is a $\Pi_{i+1}$-sequence for any $\Pi^{1}_{0}$-sequence $\{\calx_{\xi}\}_{\xi<\alp}$ and any $\Pi_{i+1}$-class $\calx$.

Therefore for any $\Pi_{i+1}$-class $\calx$ and any $\Pi_{i+1}$-sequence $\{\calx_{\xi}\}_{\xi<\alp}$
\beqnarr
&& P\in\calx\cap \bigcap\{M_{i+1}(\calx_{\xi}):\xi<\alp\} \Rarw \nonumber \\
&& P\in M_{i}^{\bet}(\calx;\{\calx_{\xi}\}_{\xi<\alp})\cap 
\bigcap\{M_{i}(M_{i}^{\bet}(\calx;\{\calx_{\xi}\}_{\xi<\alp})\cap\calx_{\xi}):\xi<\alp\}
\label{eq:5aprl1}
\eeqnarr
for any 'ordinal' $\bet<ord(P)^{+}$. This is seen by induction on $\bet$ using the following fact:
\beqn\label{eq:fact}
\mbox{For a } \Pi^{1}_{0}\mbox{-class } \calx \mbox{ and a } \Pi_{i}\mbox{-class } \caly, 
\caly\cap M_{i}(\calx)\incl M_{i}(\caly\cap\calx)
\eeqn

In general $\bet$ can be replaced by a $\Del_{1}$-definable well founded relation on each $P$. For example $\mbox{L}_{\kap}\in M_{i+1}(\mbox{L}^{t})\Rarw \mbox{L}_{\kap}\in M_{i}^{<\kap^{+}}(\mbox{L})=\bigcap\{M_{i}^{\bet}(\mbox{L}):\bet<\kap^{+}\}$ for the next admissible $\kap^{+}$ to $\kap$.

 Let $\calx\prec_{i}\caly$ denote 
\beqnarrs
\calx\prec_{i}\caly & :\Lrarw & \caly\incl M_{i}(\calx)\mbox{, i.e., } \fal P\in\caly(P\in M_{i}(\calx)), \\
\mbox{ and }
\calx\preceq_{i}\caly & :\Lrarw & \calx=\caly \mbox{ or } \calx\prec_{i}\caly.
\eeqnarrs
For example (\ref{eq:5aprl1}) is written as 
$M_{i}^{\nu}(\calx;\{\calx_{\xi}\}_{\xi<\alp})\cap\calx_{\xi}\prec_{i}
M_{i}^{\bet}(\calx;\{\calx_{\xi}\}_{\xi<\alp})\prec_{i}
\calx\cap \bigcap\{M_{i+1}(\calx_{\xi}):\xi<\alp\}$ for $\nu<\bet$ and $\xi<\alp$, 
and (\ref{eq:fact}) is as $\caly\cap\calx\prec_{i}\caly\cap M_{i}(\calx)$. Observe that the relation $\prec_{i}$ is transitive and wellfounded since the $\in$-relation is wellfounded, and $\prec_{i+1}\incl\prec_{i}$.

$\Pi_{3}$-reflecting universes $\mbox{L}$ are so simple to analyse: $\mbox{L}$ can be resolved or approximated by $M_{2}^{\alp}(\mbox{L})$, cf. \cite{odpi3}. Nevertheless $\Pi_{N}$-reflecting universes $\mbox{L}$ for $N\geq 4$ involves a complicated ramification process to resolve by using iterations of $\Pi_{2}$-recursively Mahlo operations $M_{2}^{\alp}(\mbox{L})$. Such a resolving is needed to define $\eta$-Mahloness. In fact the following ramification is inspired from a finer analysis of o.d.'s, which has been constructed purely from combinatorial considerations. 
Since the resolving of $\Pi_{N}$-reflecting universes $\mbox{L}$ by iterations of $\Pi_{2}$-recursively Mahlo operations is so complicated, we first explain the simplest case. Namely $\Pi_{3}$-reflection on $\Pi_{3}$-reflecting ordinals.

\subsection{$\Pi_{3}$-reflection on $\Pi_{3}$-reflecting ordinals}\label{subsec:M32}

\subsubsection{Mahlo classes for $M_{3}^{2}$}\label{subsubsec:M32MC}

Let $M_{3}^{2}=M_{3}(M_{3})$ denote the class of $\Pi_{3}$-reflecting universes on $\Pi_{3}$-reflecting universes $M_{3}$. And let $\mbox{L}\in M_{3}^{2}$.

Define Mahlo classes $M_{2}(2;\bet), M_{2}((2,1);(\bet,\bet_{0}))$ on $\mbox{L}^{t}\cup\{\mbox{L}\}$ as follows.
Let $\pi:=ord(\mbox{L})$.
\beqnarrs
M_{2}(2;\bet) & = & \bigcap\{M_{2}(M_{2}(2;\nu)\cap M_{3}):\nu<\bet\} \\
M_{2}((2,1);(\bet,\bet_{0})) & = & M_{2}(2;\bet) \cap \bigcap\{M_{2}(M_{2}((2,1);(\bet,\nu))):\nu<\bet_{0}\}
\eeqnarrs
and let
\[\calt = \{M_{2}(2;\bet)\cap M_{3}, M_{2}((2,1);(\bet,\bet_{0})):\bet<\veps_{\pi+1},\bet_{0}<\pi\}.\]

Note that $M_{3}=M_{2}(2;0)\cap M_{3}$ and $M_{2}^{\bet}(\mbox{L})=M_{2}(2,1);(0,\bet))$.
For $\calx\in\calt$ define a pair $h(\calx)=\langle h_{0}(\calx),h_{1}(\calx)\rangle$ by
\beqn\label{eq:hcalx}
h(M_{2}(2;\bet)\cap M_{3}):=\langle\bet,\pi\rangle \mbox{ and } h(M_{2}((2,1);(\bet,\bet_{0})):=\langle\bet,\bet_{0}\rangle
\eeqn
Then $h(\calx)<_{lex}h(\caly) \Rarw \calx\prec_{2}\caly$ for $\calx, \caly\in\calt$. 
Namely we see the following facts as in Section \ref{sec:prl5a}:
\beqnarr
&& \nu<\bet_{0} \Rarw M_{2}((2,1);(\bet,\nu))\prec_{2}M_{2}((2,1);(\bet,\bet_{0})) \label{eq:M32.6} \\
&& \nu<\bet \Rarw M_{2}(2;\nu)\cap M_{3}\prec_{2}M_{2}(2;\bet) \label{eq:M32.4} \\
&& M_{2}((2,1);(\bet,\gam))\prec_{2}M_{2}(2;\bet)\cap M_{3} \label{eq:M32.3} \\
&& \nu<\bet \Rarw M_{2}(2;\nu)\cap M_{3}\prec_{2}M_{2}((2,1);(\bet,\bet_{0})) \label{eq:M32.5} \\ 
&& \fal\calx\in\calt[\calx\prec_{2} M_{3}^{2}] \nonumber
\eeqnarr

\subsubsection{Ordinal diagrams $Od(M_{3}^{2})$ for $M_{3}^{2}$}\label{subsubsec:odM32}

We define the subsystem $Od(M_{3}^{2})\subset Od$ of ordinal diagrams.

For $\rho\in\cald^{Q}_{\sig}\cap Od(M_{3}^{2})$, $Q(\rho)$ is a pair $(rg_{2}(\rho),st_{2}(\rho))$ of o.d.'s. Namely $Q(\rho)=(2,rg_{2}(\rho),\sig,st_{2}(\rho))$.
Let $\alp\prec_{2}\bet:\Lrarw\alp\prec\bet$.

\bdf\label{df:odM32} $Od(M_{3}^{2})$.\\
{\rm The system} $Od(M_{3}^{2})$ {\rm of o.d.'s is obtained from} $Od$ {\rm by restricting the construction} $(\sig,q,\alp)\mapsto d_{\sig}^{q}\alp$ {\rm in Definition \ref{df:opi}.\ref{caldQ} as follows: first set} $\pi^{+}:=\veps_{\pi+1}:=\infty$ {\rm and}
\[
M_{3}:=\{\pi\}\cup\{\alp\in\cald^{Q}\cap Od(M_{3}^{2}): rg_{2}(\alp)=\pi\}.
\]
{\rm Assume} $\alp \in Od(M_{3}^{2}) \, \& \,  \sig\in\{\pi\}\cup\cald^{Q} \spand q=(\kap,\nu)\incl Od(M_{3}^{2})$ {\rm such that} 
 \benu
 \item 
 \beqn\label{eq:rgM32}
 \sig\preceq\kap\in M_{3} 
 \eeqn
 \item 
 \beqn\label{eq:stboundM32pi}
 \sig=\pi \Rarw \nu\leq\alp
 \eeqn
 {\rm and}
 \beqn\label{eq:stboundM32}
 \nu<\kap^{+}
 \eeqn
 \eenu
{\rm Put} $\rho:=d^q_\sig\alp\in \cald^{Q}_{\sig}\incl Od$. {\rm For this} $\rho$ {\rm define}
\[pd_{2}(\rho)=\sig, \, st_{2}(\rho)=\nu, \, rg_{2}(\rho)=\kap.\]
{\rm Then}
$\rho \in Od(M_{3}^{2})$ {\rm if the following conditions are fulfilled besides (\ref{eq:Odmu}) in Definition \ref{df:opi}.\ref{caldQ}:}
\bdes
\item[Case1] $\rho\in M_{3}${\rm , i.e.,} $\kap=rg_{2}(\rho)=\pi${\rm : Then} 
\[\fal\eta\in M_{3}|\pi[\rho\prec\eta \Rarw st_{2}(\rho)<st_{2}(\eta)].
\]
\item[Case2] $\rho\not\in M_{3}${\rm , i.e.,} $\kap=rg_{2}(\rho)<\pi${\rm : Then} 
\[
rg_{2}(\rho)=\min\{\kap<\pi: \rho\prec\kap\in M_{3}\}
\]
{\rm i.e.,} $\fal\tau[\rho\prec\tau\prec\kap \Rarw rg_{2}(\tau)=\kap]$ {\rm and}
\[\fal\tau[\rho\prec\tau\prec\kap \Rarw st_{2}(\rho)<st_{2}(\tau)].
\]
\item [$(\cald.2)$]
   \beqn \label{cnd:KstM32}
   \fal\tau\leq rg_{2}(\rho)(K_{\tau}st_{2}(\rho)<\rho)
   \eeqn
\edes
\edf

Now for each o.d. $\rho\in\cald^{Q}\cap Od(M_{3}^{2})$ we associate a Mahlo class $\calx(\rho)\in\calt$ and a pair $h(\rho)$, cf. (\ref{eq:hcalx}) as follows:
\bdes
\item[Case1] $\rho\in M_{3}$: $\calx(\rho)=M_{2}(2;\bet)\cap M_{3}$ for $\bet=st_{2}(\rho)$, and
$h(\rho):=h(\calx(\rho))=\langle\bet,\pi\rangle$.
\item[Case2] $\rho\not\in M_{3}$: $\calx(\rho)=M_{2}((2,1);(\bet,\gam))$ for $\gam=st_{2}(\rho)$ and $\bet=st_{2}(rg_{2}(\rho))$, and
$h(\rho):=h(\calx(\rho))=\langle\bet,\gam\rangle$.
\edes

Then we see

\blem\label{lem:lexM32}
If $\del\prec_{2}\eta$, then $h(\del)<_{lex}h(\eta)$ and hence $\calx(\del)\prec_{2}\calx(\eta)$.
\elem
\bprf
Suppose $\eta=pd_{2}(\del)$.

(\ref{eq:M32.3}) corresponds to the case: $rg_{2}(\del)=\eta\in M_{3}$ with $st_{2}(\eta)=\bet$, 
$st_{2}(\del)=\gam<\eta^{+}$.

(\ref{eq:M32.6}) corresponds to the case: $rg_{2}(\del)=rg_{2}(\eta)<\pi$ with $\nu=st_{2}(\del)<st_{2}(\eta)=\bet_{0}$.

(\ref{eq:M32.4}) corresponds to the case: $\del,\eta\in M_{3}$ and $\nu=st_{2}(\del)<st_{2}(\eta)=\bet$.

Finally consider the case when $\del\in M_{3}$ and $\eta\not\in M_{3}$. Then $rg_{2}(\eta)=\min\{\sig: \del\prec\sig\in M_{3}\}$. (\ref{eq:M32.5}) corresponds to this case with $\nu=st_{2}(\del)<st_{2}(rg_{2}(\eta))=\bet$ (and $st_{2}(\eta)=\bet_{0}$).
\eprf

We could prove the wellfoundedness of $Od(M_{3}^{2})|\Ome$ using Lemma \ref{lem:lexM32} and distinguished classes as in Section \ref{sec:5awf}.

\subsubsection{Wellfoundedness proof for $Od(M_{3}^{2})$}\label{subsubsec:odM32wf}

We show the

\bth\label{th:id5wfeachpiM32}
For each $\alp<d_{\Ome}\veps_{\pi+1}$,i.e., each $\alp\in Od(M_{3}^{2})|\Ome$, $[\Pi^{0}_{2},\Pi^{0}_{2}]\mbox{{\rm -Fix}}$ proves that $(Od(M_{3}^{2})|\alp,<)$ is a well ordering.
\eth

Work in $[\Pi^{0}_{2},\Pi^{0}_{2}]\mbox{-Fix}$.

\bdf
\[
\alp\lhd\bet :\Lrarw \alp,\bet\in\cald^{Q} \spand \alp\prec\bet\prec rg_{2}(\alp).
\]
\edf

Observe that for $\alp,\bet\in\cald^{Q}$ with $\alp\prec\bet$, $\alp\lhd\bet$ iff $\alp\in M_{3}$ or $rg_{2}(\alp)=rg_{2}(\bet)$. Hence $\alp\lhd\bet \Rarw h(\alp)<_{lex}h(\bet)$. Specifically supposing $\gam\prec\alp$ we see:
\benu
\item $\alp\not\in M_{3}$: Then $\gam\lhd\alp$ iff 
$\gam\not\in M_{3}\spand rg_{2}(\gam)=rg_{2}(\alp)\mbox{(, and hence } h_{0}(\gam)=h_{0}(\alp)\mbox{)} \spand h_{1}(\gam)=st_{2}(\gam)<st_{2}(\alp)=h_{1}(\alp)$, cf. (\ref{eq:M32.6}), or 
$\gam\in M_{3} \spand h_{0}(\gam)=st_{2}(\gam)<st_{2}(rg_{2}(\alp))=h_{0}(\alp)$, cf. (\ref{eq:M32.5}).
\item $\alp\in M_{3}$: Then $\gam\lhd\alp$ iff $\gam\in M_{3} \spand h_{0}(\gam)=st_{2}(\gam)<st_{2}(\alp)=h_{0}(\alp)$, cf. (\ref{eq:M32.4}).
\eenu
Moreover if $M_{3}\not\ni\gam\prec\alp\in M_{3}$, then $rg_{2}(\gam)\preceq\alp$ and $h_{0}(\gam)=st_{2}(rg_{2}(\gam))\leq st_{2}(\alp)=h_{0}(\alp) \spand h_{1}(\gam)=st_{2}(\gam)<\pi=h_{1}(\alp)$, cf. (\ref{eq:M32.3}).

\bdf\label{df:oprtM32V}
\[
\alp\in V(X):\Lrarw \fal\gam\lhd\alp[\gam\in\calg(X) \rarw \gam\in X].
\]
\edf

Now let us define an operator $\Gam_{32}$ on $Od(M_{3}^{2})$ from $V(X)$.

\bdf\label{df:oprtM32}
\benu
\item $\alp\in\Gam_{30}(X)$ {\rm iff} $\pi>\alp\not\in M_{3}$, $\alp\in\calg(X)\cap V(X)$ {\rm and}
$ [\alp\in SR \Rarw \fal\gam\in\cald_{\alp}(\gam\in\calg(X) \rarw \gam\in X)]$.
\item $\alp\in\Gam_{32}(X) $ {\rm iff}
\[
\alp\in\Gam_{30}(X) \, \vee\, [\Gam_{30}(X)\incl X\spand \alp\in M_{3}|\pi\spand \alp\in\calg(X)\cap V(X)].
\]
\eenu
\edf

Let us examine the complexity of these operators. Both $V$ and $\Gam_{30}$ are $\Pi^{0}_{2}$, and hence $\Gam_{32}$ is $[\Pi^{0}_{2},\Pi^{0}_{2}]$.

We write $\Gam$ for $\Gam_{32}$, $| \alp|$ for $| n|_{\Gam_{32}}$. 

We see easily that $\Gam=\Gam_{32}$ enjoys the hypotheses {\bf ($\Gam$.0)}, {\bf ($\Gam$.1)} and {\bf ($\Gam$.5)} in Subsection \ref{subsec:ensureoperators}. Furthermore {\bf ($\Gam$.3)} and {\bf ($\Gam$.4)} follow from the facts: if $\alp\not\in R^{\prime}$ or $\alp\in SR$, then $\alp\in V(X)$ for any $X$.

We next show that $\Gam$ enjoys the hypothesis {\bf ($\Gam$.2)}.
\bth\label{th:AM} 
Assume $\alp,\bet\in\calw$. Then
\[\alp<\bet \Rarw x=|\alp|<|\bet|=y.\]
\eth
{\bf Proof} of Theorem \ref{th:AM} for $\Gam=\Gam_{32}$. Assume $\alp,\bet\in\calw$ and $\alp<\bet$. Put $x=|\alp|, y=|\bet|$. We show $x<y$ by induction on the natural sum $x\#y$. Suppose $x\geq y$. Put $X=\Gam^{x}, Y=\Gam^{y}$.  We show $\alp\in Y$. As in \cite{Wienpi3d} we see, using IH, $\alp\in\calg(X)|\bet=\calg(Y)|\bet\spand \alp\in V(Y)$, and we can assume $\alp\prec\bet\in\cald^{Q}$ and $\alp\not\leq Y$ by IH.
There are two cases to consider.
\\
{\bf Csae1} $\exi\sig[\alp\preceq\sig\spand \sig\lhd\bet]$: Then by Lemma \ref{lem:CX6+} and $\alp\not\leq Y$ we have $\sig\in\calg(Y)$. $\bet\in V(Y)$ yields $\alp\leq\sig\in Y$.
\\
{\bf Case2} Otherwise : We have $\sig\not\!\!\lhd\bet$ for any $\sig$ with $\alp\preceq\sig\prec\bet$, and $rg_{2}(\alp)\preceq\bet\spand\alp\not\in M_{3}$. We claim $M_{3}\ni rg_{2}(\alp)=\bet$. Otherwise we would have $rg_{2}(\alp)\lhd\bet$.
Thus $\bet\in M_{3}$, and hence $\alp\in\Gam_{30}(Y)\incl Y$. 
\eprf

Let $V:=V(\calw)$.

\blem\label{lem:GWpiM32}
\[\alp\in\calg\cap V \Rarw \alp\in\calw.\]
\elem
 
The following lemma is seen from Lemma \ref{lem:CX4aro} and (\ref{cnd:KstM32}).

\blem\label{lem:GWQpiM32}
If $\rho\in\calg\cap\cald^{Q}$ and $\kap:=rg_{2}(\rho)$, then $st_{2}(\rho)\in\calc^{\kap}(\calw)$.

In particular, $\rho\in\calg\cap\cald^{Q}\cap M_{3} \Rarw st_{2}(\rho)\in\calc^{\pi}(\calw)=\calw_{\pi}$.
\elem

\blem\label{lem:N4aroM32}
For $\rho\in\cald^{Q}\cap Od(M_{3}^{2})$, $Q(\rho)=\{st_{2}(\rho),rg_{2}(\rho)\}\leq\max\{b(\rho),\pi\}$.
\elem
\bprf
First off, $rg_{2}(\rho)\leq\pi$. It remains to show $st_{2}(\rho)\leq\max\{b(\rho),\pi\}$. By (\ref{eq:stboundM32pi}) and (\ref{eq:stboundM32}) we can assume $\rho\in M_{3}\spand pd_{2}(\rho)<\pi$. Let $\alp_{1}$ denote the diagram such that $\rho\prec\alp_{1}\in\cald_{\pi}$. Then $\alp_{1}\in M_{3}$ and $st_{2}(\rho)<st_{2}(\alp_{1})\leq b(\alp_{1})<b(\rho)$ by Lemma \ref{lem:Npi11exist}.
\eprf

\blem\label{lem:id5wf21}(cf. Definition \ref{df:id4wfA}.\ref{df:id4wfA.1}.) For each $n\in\ome$
\[\fal\alp\in\calw_{\pi}|\ome_{n}(\pi+1)\fal q\incl\calw_{\pi}|\ome_{n}(\pi+1) A(\alp,q).\]
\elem
\bprf
 We have to show for each $n\in\ome$
\[
\fal\alp\in\calw_{\pi}|\ome_{n}(\pi+1)\fal q\incl\calw_{\pi}|\ome_{n}(\pi+1) A(\alp,q).
\] 
By main induction on $\alp\in\calw_{\pi}|\ome_{n}(\pi+1)$ with subsidiary induction on $q\incl\calw_{\pi}|\ome_{n}(\pi+1)$.
Here observe that if $\bet_{1}\in\cald$ with $b(\bet_{1})<\ome_{n}(\pi+1)$, then by Lemma \ref{lem:N4aroM32} we have $Q(\bet_{1})\leq\max\{b(\bet_{1}),\pi\}<\ome_{n}(\pi+1)$.

Let $\alp_{1}\in\cald_{\sig}|\pi$ with $\sig\in\calw_{\pi}$ and $\alp=b(\alp_{1})\spand q=Q(\alp_{1})$.
By Theorem \ref{th:id5wf21} we have $\alp_{1}\in\calg$. We show $\alp_{1}\in\calw$. By Lemma \ref{lem:GWpiM32} it suffices to show $\alp_{1}\in V$.
By $\sig\in\calw_{\pi}$ we have $\sig\in V\cup\{\pi\}$. If $\sig=\pi$, i.e., $pd_{2}(\alp_{1})=\pi$, then $\alp_{1}\in M_{3}\spand st_{2}(\alp_{1})\in\calw_{\pi}$ by Lemma \ref{lem:GWQpiM32}, and we see $\alp_{1}\in V$ from $\calg\ni\gam\lhd\alp_{1}\Rarw \calw_{\pi}\ni st_{2}(\gam)<st_{2}(\alp_{1})$.
Therefore we can assume $\sig\in V$. If $\sig\prec rg_{2}(\alp_{1})$, then $\alp_{1}\lhd\sig\in V$, and hence $\alp_{1}\in\calw$. Assume $rg_{2}(\alp_{1})=\sig$.

Let $\calg\ni\gam\lhd\alp_{1}$. We have to show $\gam\in\calw$. We can assume $rg_{2}(\gam)=\sig$ by $\sig\in V$. Thus $st_{2}(\gam)<st_{2}(\alp)$. By Lemma \ref{lem:GWQpiM32} we have $st_{2}(\gam), st_{2}(\alp)\in\calc^{\sig}(\calw)$. 
We have $\calb_{>\sig}(st_{2}(\gam))<b(\alp_{1})=\alp$, cf. Lemma \ref{lem:5.4}.\ref{lem:5.4.3}. Lemma \ref{lem:id3wf20} with $\mbox{MIH}(\alp)$ yields $st_{2}(\gam)\in\calw_{\pi}$. Therefore $\gam\in\calw$ is seen by induction on $st_{2}(\gam)$.

We are done.
\eprf

Lemma \ref{lem:id5wf21} yields Lemma \ref{th:id4wf22}: $\alp_{1}\in\calw_{\pi}$ for each $\alp_{1}\in Od(M_{3}^{2})$ as in \cite{Wienpi3d}.

\blem\label{th:id4wf22}
For each $\alp_{1}$, $\alp_{1}\in\calw_{\pi}$.
\elem

Consequently Lemma \ref{lem:WWOme} yields Theorem \ref{th:id5wfeachpiM32}.

\subsection{Mahlo classes for $\Pi_{N}$-reflection}\label{sec:prl5a.1}

In what follows $N$ denotes a fixed integer $N\geq 4$.
Let us resolve a $\Pi_{N}$-reflecting universe L.
Now let $M_{i}(\alp;\bet)\, (i\geq 2)$ denote the following class:
\[M_{i}(\alp;\bet):=M_{1,i}^{\bet}(\{\calx_{\xi}\}_{\xi<\alp})
\mbox{ with the } \Pi_{i+2}\mbox{-class } \calx_{\xi}=M_{i+1}^{\xi}(\mbox{L}).\]
Namely
\[
M_{i}(\alp;\bet)=
\bigcap\{M_{i}(M_{i}(M_{i}(\alp;\nu)\cap M_{i+1}^{\xi}(L))\cap M_{i+1}^{\del}(L)): \del\leq\xi<\alp, \nu<\bet\}
\]
Then $M_{i}(\alp;\bet)$ is again a $\Pi_{i+1}$-class and hence from (\ref{eq:5aprl1})

\[
M_{i}(\alp;\bet)\prec_{i}M_{i+1}^{\alp}(\mbox{L})\, (\alp>0)
\spand M_{i}(\alp;\bet)\supset M_{i+1}^{\alp}(\mbox{L}).
\]

Moreover let $\bar{\alp}=(\alp_{n}>\cdots>\alp_{1})\,(n\geq 0)$ denote a decreasing sequence of ordinals and $\bar{\bet}=(\bet_{n},\ldots,\bet_{1})$ a sequence of ordinals of the same length. By induction on the length $n$ of the sequences $\bar{\alp}, \bar{\bet}$ we define classes $M_{i}(\bar{\alp};\bar{\bet})$ as follows. $M_{i}(\langle\rangle;\langle\rangle)=\mbox{L}^{t}\cup\{\mbox{L}\}$ for the empty sequence $\langle\rangle$. Let $\bar{\alp}*(\alp)=(\alp_{n}>\cdots>\alp_{1}>\alp)$ and $\bar{\bet}*(\bet)=(\bet_{n},\ldots,\bet_{1},\bet)$ denote the concatenated sequences for ordinals $\alp,\bet$ with $\alp<\alp_{1}$. Then define 
\[
M_{i}(\bar{\alp}*(\alp);\bar{\bet}*(\bet))=
M_{i}^{\bet}(\caly;\{\calx_{\xi}\}_{\xi<\alp})
\]
with the $\Pi_{i+2}$-class $\calx_{\xi}=M_{i+1}^{\xi}(L)$ and $\caly=M_{i}(\bar{\alp};\bar{\bet})$.

Then as above we see that $M_{i}(\bar{\alp};\bar{\bet})$ are $\Pi_{i+1}$-classes and
\beqnarr
&& M_{i}(\bar{\alp}*(\alp);\bar{\bet}*(\bet))=M_{i}(\bar{\alp};\bar{\bet})\cap  \nonumber \\
&& \bigcap\{M_{i}(M_{i}(M_{i}(\bar{\alp}*(\alp);\bar{\bet}*(\nu))\cap M_{i+1}^{\xi}(L))\cap M_{i+1}^{\del}(L)): 
\nonumber \\
&& \del\leq\xi<\alp, \nu<\bet\} 
\label{eq:Mexp}
\eeqnarr

Thus for any $\bar{\alp}=(\alp_{n}>\cdots>\alp_{1})$, $\alp<\alp_{0}<\alp_{1}$ and any ordinals $\bar{\bet}, \bet_{0}$, we see from (\ref{eq:5aprl1})

\beqn\label{eq:5aprl3a}
M_{i}(\bar{\alp};\bar{\bet})\cap M_{i+1}^{\alp}(\mbox{L})\prec_{i}M_{i+1}^{\alp_{n}}(\mbox{L})\, (n>0)
\eeqn
and
\beqn\label{eq:5aprl3b}
M_{i}(\bar{\alp}*(\alp_{0});\bar{\bet}*(\bet_{0}))\cap M_{i+1}^{\alp}(\mbox{L})\prec_{i}
M_{i}(\bar{\alp};\bar{\bet})\cap M_{i+1}^{\alp_{0}}(\mbox{L})
\eeqn

Now let us depict this ramification process as wellfounded trees $\calt^{i}_{N}\, (2\leq i<N)$. These trees are defined by induction on $N-i$. Each class in $\calt^{i}_{N}$ is a $\Pi_{N}$-class for any $i$. Pick a $\Del_{1}$-well ordering $<$ on $\mbox{L}$. Assume its order type is, for example, $\veps_{\pi+1}$ for the least ordinal $\pi$ not in L. For a $\Pi_{N}$-reflecting universe $\mbox{L}$, the singleton class $\{\mbox{L}\}$ sits on each root of $\calt^{i}_{N}$. In the tree $\calt_{N}^{N-1}$ its sons are $\Pi_{N}$-classes $M_{N-1}^{\alp}(\mbox{L})$ for each 'ordinal' $\alp<\veps_{\pi+1}$ such that $M_{N-1}^{\alp}(\mbox{L})\prec_{N-1}M_{N}(\mbox{L})\ni \mbox{L}$. Further the node $M_{N-1}^{\alp}(\mbox{L})$ has sons $M_{N-1}^{\bet}(\mbox{L})$ for $\bet<\alp$. Each node $M_{N-1}^{\alp}(\mbox{L})$ can be identified with the ordinal $\alp<\veps_{\pi+1}$, and the relation $M_{N-1}^{\bet}(\mbox{L})\prec_{N-1}M_{N-1}^{\alp}(\mbox{L})$ on $\calt_{N}^{N-1}$ with $\bet<\alp$. Note that these are $\Del_{1}$-relations on L.

Suppose a wellfounded tree $\calt^{i+1}_{N}$ has been constructed for $2\leq i<N-1$ so that each class in $\calt^{i+1}_{N}$ is a $\Pi_{N}$-class, and if a class $\calx$ is a son of a class $\caly$, then $\calx\prec_{i+1}\caly$, i.e., the tree ordering is compatible with the relation $\prec_{i+1}$. Classes in the tree $\calt_{N}^{i+1}$ are assumed to be ordered by the relation $\prec_{i+1}$. Moreover suppose that the tree $\calt^{i+1}_{N}$ and the relation $\prec_{i+1}$ on $\calt^{i+1}_{N}$ are coded by $\Del_{1}$-relations on L ({\it Coding Supposition\/}). Then another wellfounded tree $\calt^{i}_{N}$ is defined as follows.

For a branch $\bar{\calx}*(\calx_{1})=(\calx_{n},\ldots,\calx_{2},\calx_{1})$ with $\calx_{1}\prec_{i+1}\calx_{2}\prec_{i+1}\cdots\prec_{i+1}\calx_{n}\prec_{i+1}\mbox{L}\spand\calx_{m}\in\calt^{i+1}_{N} \, (n\geq 1)$ in the tree $\calt^{i+1}_{N}$ and a sequence $\bar{\bet}*(\bet_{1})=(\bet_{n},\ldots,\bet_{1})$ of ordinals, a class $M_{i}(\bar{\calx};\bar{\bet})$ is defined by replacing $\xi<\alp$ in the definition (\ref{eq:Mexp}) of $M_{i}(\bar{\alp};\bar{\bet})$ by the wellfounded relation $\prec_{i+1}$ in $\calt_{N}^{i+1}$:
\beqnarr
&& M_{i}(\bar{\calx}*(\calx_{1});\bar{\bet}*(\bet_{1}))=M_{i}(\bar{\calx};\bar{\bet})\cap \nonumber \\
&& \bigcap\{M_{i}(M_{i}(M_{i}(\bar{\calx}*(\calx_{1});\bar{\bet}*(\nu))\cap \calx_{0})\cap\calx^{0}): \nonumber \\
&& \calx^{0}\preceq_{i+1}\calx_{0}\prec_{i+1}\calx_{1}, \calx^{0},\calx_{0}\in\calt_{N}^{i+1}, \nu<\bet_{1}\} 
\label{eq:5aprl1.5}
\eeqnarr
By the {\it Coding Supposition\/} $M_{i}(\bar{\calx};\bar{\bet})$ is a uniform $\Pi_{i+1}$-class, and hence as in (\ref{eq:5aprl3a}) and (\ref{eq:5aprl3b}) we see for $\calx\prec_{i+1}\calx_{0}\prec_{i+1}\calx_{1}\prec_{i+1}\cdots\prec_{i+1}\calx_{n}\prec_{i+1}\{\mbox{L}\}\, (n\geq 0)$ in the tree $\calt^{i+1}_{N}$ with $\bar{\calx}=(\calx_{n},\ldots,\calx_{1})$ and a sequence $\bar{\bet}=(\bet_{n},\ldots,\bet_{1})$ of ordinals:

\beqn\label{eq:5aprl3ax}
M_{i}(\bar{\calx};\bar{\bet})\cap \calx_{0}\prec_{i}\calx_{n}\, (n>0)
\eeqn
and
\beqn\label{eq:5aprl3bx}
M_{i}(\bar{\calx}*(\calx_{0});\bar{\bet}*(\bet_{0}))\cap \calx\prec_{i}M_{i}(\bar{\calx};\bar{\bet})\cap\calx_{0}
\incl M_{i}(\bar{\calx}*(\calx_{0});\bar{\bet}*(\bet_{0}))
\eeqn

Each class $\calx$ in the tree $\calt^{i}_{N}$ except $\{\mbox{L}\}$ is of the form
\beqn\label{eq:treei}
\calx=M_{i}(M_{i}(\bar{\calx};\bar{\bet})\cap\calx_{0})\cap\calx^{0}
\eeqn
for some branch $\calx^{0}\preceq_{i+1}\calx_{0}\prec_{i+1}\calx_{1}\prec_{i+1}\cdots\prec_{i+1}\calx_{n}\prec_{i+1}\mbox{L}\, (n\geq 0)$ in the tree $\calt^{i+1}_{N}$ with $\bar{\calx}=(\calx_{n},\ldots,\calx_{1})$ and a sequence $\bar{\bet}=(\bet_{n},\ldots,\bet_{1})$ of ordinals. Thus $\calx$ is a $\Pi_{N}$-class.

The root $\mbox{L}$ has sons $M_{i}(\calx)\cap\calx$ for each son $\calx$ of L in the tree $\calt_{N}^{i+1}$, i.e., in (\ref{eq:treei}) $\bar{\calx}=\bar{\bet}=\langle\rangle$ and $\calx_{0}=\calx^{0}=\calx$: $M_{i}(\calx)\cap\calx\prec_{i}\mbox{L}$, cf. (\ref{eq:5aprl1}). 

Now the class $\calx\in\calt^{i}_{N}$ defined in (\ref{eq:treei}) has sons of three kinds.
 
A class of the form 
\[\caly=M_{i}(M_{i}(\bar{\calx};\bar{\bet})\cap\calx_{0})\cap\caly^{0}
\mbox{ with } \calt^{i+1}_{N}\ni\caly^{0}\prec_{i+1}\calx^{0}\]
 is a first son of the father $\calx$.
 
  $\caly\prec_{i}\calx$ is seen from (\ref{eq:fact}) and the fact that $M_{i}(M_{i}(\bar{\calx};\bar{\bet})\cap\calx_{0})$ is a $\Pi_{i+1}$-class.

Second a second son is a class 
\[\calz=
M_{i}(M_{i}(\bar{\calx}*(\calx_{0});\bar{\bet}*(\gam))\cap\calz_{0})\cap\calz_{0}
\]
with an ordinal $\gam$ and a class $\calz_{0}\prec_{i+1}\calx_{0}$ in $\calt^{i+1}_{N}$.

$\calz\prec_{i}\calx$ is seen from $\calz_{0}\prec_{i+1}\calx_{0}$ and (\ref{eq:5aprl3bx}):\\
$\calz=
M_{i}(M_{i}(\bar{\calx}*(\calx_{0});\bar{\bet}*(\gam))\cap\calz_{0})\cap\calz_{0}
\prec_{i}
M_{i}(\bar{\calx};\bar{\bet})\cap\calx_{0}\prec_{i}\calx$.

Finally for the case $n>0$ a third son is a class 
\[\calu=
M_{i}(M_{i}(\bar{\calx}\uarw k;\bar{\bet}\uarw(k+1)*(\gam))\cap\calu_{0})\cap\calu_{0}
\]
for a $k\, (1\leq k\leq n)$, an ordinal $\gam<\bet_{k}$ and a class $\calu_{0}\prec_{i+1}\calx_{k}$ in $\calt^{i+1}_{N}$, where $\bar{\calx}\uarw k=(\calx_{n},\ldots,\calx_{k})$ and $\bar{\bet}\uarw(k+1)*(\gam)=(\bet_{n},\ldots,\bet_{k+1},\gam)$.

$\calu\prec_{i}\calx$ is seen from $\calu_{0}\prec_{i+1}\calx_{k}\spand \gam<\bet_{k}$ and (\ref{eq:5aprl1.5}): 
\beqnarrs
\calu & = & M_{i}(M_{i}(\bar{\calx}\uarw k;\bar{\bet}\uarw(k+1)*(\gam))\cap\calu_{0})\cap\calu_{0} \\
& \prec_{i} & M_{i}(\bar{\calx}\uarw k;\bar{\bet}\uarw(k+1)*(\bet_{k}))
\supset M_{i}(\bar{\calx};\bar{\bet})\prec_{i}\calx.
\eeqnarrs

Note that we have $\calu\prec_{i}\calz\prec_{i}\caly\prec_{i}\calx$ for the sons $\caly, \calz, \calu$ of $\calx$.

These three kinds of classes are sons of $\calx$ and this completes a description of the tree $\calt^{i}_{N}$.  
The tree represents a ramification procedure in decomposing a $\Pi_{N}$-reflecting universe $\mbox{L}$ in terms of 
iterations of $\Pi_{i}$-recursively Mahlo operations.

The class $\calx\in\calt^{i}_{N}$ defined in (\ref{eq:treei}) can be coded by a code $(i;\bar{x};\bar{\bet};x_{0},x^{0})$, 
where $\bar{x}$ are codes for $\bar{\calx}$, and $x_{0}$ [$x^{0}$] for $\calx_{0}$ [for $\calx^{0}$], resp. 
By the {\it Coding Supposition\/}, '$x$ is a code for a class $\calx$ in $\calt^{i+1}_{N}$' is recursive, 
and so is the relation $x\prec_{i+1}y :\Lrarw \calx\prec_{i+1}\caly\, (\calx,\caly\in\calt^{i+1}_{N})$. 
Therefore the tree $\calt^{i}_{N}$ and the relation $\prec_{i}$ on $\calt^{i}_{N}$ are again coded by $\Del_{1}$-relations on L.

In this way we get a wellfounded tree $\calt_{N}=\calt_{N}^{2}$ ordered by the relation $\prec_{2}$. 
For $2\leq i<N-1$ each class $\calx\neq \{\mbox{L}\}$ in $\calt_{N}^{i}$ is of the form described in (\ref{eq:treei})
for some branch 
$\calx^{0}\preceq_{i+1}\calx_{0}\prec_{i+1}\calx_{1}\prec_{i+1}\cdots\prec_{i+1}\calx_{n}\prec_{i+1}\mbox{L}\, (n\geq 0)$ 
in the tree $\calt^{i+1}_{N}$ and a sequence $(\bet_{n},\ldots,\bet_{1})$ of ordinals. 
Therefore we can associate its {\it construction tree\/} with depth $N-i$: 
$\calx$ sits on the root and its sons are $\{(\calx_{m},\bet_{m}): 1\leq m\leq n\}$ and $\calx_{0},\calx^{0}$, 
and each son $\calx_{m},\calx^{0}$ has sons and so on. 
Does the construction tree remind you an ordinal structure with addition and exponentiation?

\subsection{Mahlo classes of ordinal diagrams}\label{sec:prl5a.2}

Now let us turn to o.d.'s in $Od(\Pi_{N})$ and explain what class in $\calt_{N}$ corresponds to a diagram of the form $\eta=d_{\sig}^{q}\alp$.

$q$ in $\eta=d_{\sig}^{q}\alp$ includes some data $st_{i}(\eta), rg_{i}(\eta)$ for $2\leq i<N$. $st_{N-1}(\eta)$ is an o.d. less than $\veps_{\pi+1}$ and $rg_{N-1}(\eta)=\pi$, while $st_{i}(\eta), rg_{i}(\eta)$ for $i<N-1$ may be undefined. If these are defined, then $\kap=rg_{i}(\eta)$ is an o.d. such that $\eta\prec_{i}\kap$, where $\prec_{i}$ is a transitive closure of the relation on o.d.'s $\{(\eta,\kap):\kap=pd_{i}(\eta)\}$ such that $\prec_{i+1}\incl\prec_{i}$. Therefore the diagram $pd_{i}(\eta)$ is a proper subdiagram of $\eta$. $q$ also determines the diagrams $pd_{i}(\eta)$. For any $\eta=d_{\sig}^{q}\alp$ and any $i$ we have $\eta\prec_{i}\pi$. $st_{i}(\eta)$ is an o.d. less than the next admissible $\kap^{+}$ to $\kap=rg_{i}(\eta)$.

Let $M^{\pi}_{i}=\{\mbox{L}\}\in\calt^{i}_{N}$ for $i$ with $2\leq i<N$. Now we associate a $\Pi_{N}$-class $M^{\eta}_{i}\in\calt^{i}_{N}$ for each such diagram $\eta$ and each $i\, (2\leq i<N)$ so that 
\beqn\label{eq:pretree}
\gam\prec_{i}\eta<\pi \Rarw M^{\gam}_{i}\prec_{i}M^{\eta}_{i}
\eeqn
in the tree $\calt^{i}_{N}$. Namely the relation $\prec_{i}$ on o.d.'s is embedded in the relation $\prec_{i}$ on the tree $\calt^{i}_{N}$.
 
 The definition of the class $M^{\eta}_{i}\in\calt^{i}_{N}$ is based on induction on $N-i$. First set $M^{\eta}_{N-1}=M_{N-1}^{\bet}(\mbox{L})\in\calt^{N-1}_{N}$ with $\bet=st_{N-1}(\eta)$. Then (\ref{eq:pretree}) is satisfied since $st_{N-1}(\eta)$ is always defined for diagrams $\eta$ of the form $d_{\sig}^{q}\alp<\pi$, and enjoys $\gam\prec_{N-1}\eta\Rarw st_{N-1}(\gam)<st_{N-1}(\eta)$.

$q$ determines a sequence $\{\eta_{i}^{m}: m<lh_{i}(\eta)\}$ of o.d.'s in $\{\bet<\pi:\eta\preceq\bet\}$ with its length $lh_{i}(\eta)=n+1>0$. The sequence enjoys the following property:
\beqn\label{eq:lem2.9}
\eta\preceq_{i+1}\eta_{i}^{0}\prec_{i+1}\eta_{i}^{1}\prec_{i+1}\cdots\prec_{i+1}\eta_{i}^{n}<\pi
\eeqn
where $\preceq_{i}$ denotes the reflexive closure of $\prec_{i}$.

Moreover $st_{i}(\eta_{i}^{m}), rg_{i}(\eta_{i}^{m})$ have to be defined for $0\leq m<n$ so that\\
 $rg_{i}(\eta_{i}^{m})\preceq_{i+1}\eta_{i}^{m+1}$, and these sequences are defined so that if $\eta=pd_{i}(\gam)$, one of the following holds, cf. Lemma \ref{lem:3.23.1} in Subsection \ref{subsec:od5fine}:
\bdes
\item[Case1]
\[
\eta=pd_{i}(\gam)=pd_{i+1}(\gam) \spand lh_{i}(\gam)=lh_{i}(\eta)\spand \fal m<lh_{i}(\gam)[\gam_{i}^{m}=\eta_{i}^{m}]
\]
\item[Case2]
\[
rg_{i}(\gam)=pd_{i}(\gam)=\eta \spand \gam_{i}^{0}=\gam \spand \fal m<lh_{i}(\eta)=lh_{i}(\gam)-1[\eta_{i}^{m}=\gam_{i}^{1+m}]
\]
\item[Case3]
\beqnarrs
&& \eta=pd_{i}(\gam)\prec_{i}rg_{i}(\gam) \spand \gam_{i}^{0}=\gam \spand  \\
&& \exi m[0<m\leq lh_{i}(\eta)-1 \spand rg_{i}(\eta_{i}^{m-1})=rg_{i}(\gam)\spand st_{i}(\eta_{i}^{m-1})>st_{i}(\gam)\spand \\
&& \fal k<lh_{i}(\eta)-m+1=lh_{i}(\gam)(k>0 \rarw \eta_{i}^{m-1+k}=\gam_{i}^{k})]
\eeqnarrs
\edes

From the sequence $\{\eta_{i}^{m}\}$ we define a class $M_{i}^{\eta}\in\calt_{N}^{i}\, (2\leq i<N-1)$ as follows, 
cf. Definition \ref{df:5etaMh}.\ref{df:5etaMh.3}:
\[
M_{i}^{\eta}=M_{i}(M_{i}(\bar{\calx};\bar{\bet})\cap\calx_{0})\cap\calx^{0}
\]
for $\calx_{m}=M_{i+1}^{\eta_{i}^{m}}\in\calt^{i+1}_{N}\, (0\leq m\leq n)$, $\calx^{0}=M_{i+1}^{\eta}\in\calt^{i+1}_{N}$, 
$\bar{\calx}=(\calx_{1},\ldots,\calx_{n})$ and o.d.'s $\bet_{m}=st_{i}(\eta_{i}^{m-1})\, (0<m\leq n)$.

From $\eta_{i}^{m}<\pi$, (\ref{eq:pretree}) for the case $i+1$ and (\ref{eq:lem2.9}) we have 
$\calx^{0}\preceq_{i+1}\calx_{0}\prec_{i+1}\calx_{1}\prec_{i+1}\cdots\prec_{i+1}\calx_{n}\prec_{i+1}\mbox{L}$ in the tree 
$\calt^{i+1}_{N}$. Thus, cf. (\ref{eq:treei}), $M_{i}^{\eta}\in\calt_{N}^{i}$. 
We verify that (\ref{eq:pretree}) holds for the case $i$. Assume $\eta=pd_{i}(\gam)$, and let $\calx=M_{i}^{\eta}$ and 
\[\calv=M_{i}^{\gam}=M_{i}(M_{i}(\bar{\calv};\bar{\gam})\cap\calv_{0})\cap\calv^{0}\]
for $\calv_{m}=M_{i+1}^{\gam_{i}^{m}}\in\calt^{i+1}_{N}\, (0\leq m<lh_{i}(\gam))$, 
$\calv^{0}=M_{i+1}^{\gam}\in\calt^{i+1}_{N}$ and o.d.'s $\gam_{m}=st_{i}(\gam_{i}^{m-1})\, (0<m<lh_{i}(\gam))$. 
We show $\calv\prec_{i}\calx$.
\\
{\bf Case1} $\eta=pd_{i}(\gam)=pd_{i+1}(\gam)$: Then 
$lh_{i}(\gam)=lh_{i}(\eta)\spand \fal m<lh_{i}(\gam)[\gam_{i}^{m}=\eta_{i}^{m}]$ and $\gam\prec_{i+1}\eta$. 
Hence $\calx_{m}=\calv_{m}\spand \bar{\bet}=\bar{\gam}\spand \calv^{0}\prec_{i+1}\calx^{0}$. 
This means $\calv$ is a first son $\caly$ of $\calx$ in $\calt_{N}^{i}$. Therefore $\calv\prec_{i}\calx$.
\\
{\bf Case2} $rg_{i}(\gam)=pd_{i}(\gam)=\eta$: Then $\gam_{i}^{0}=\gam$, 
$\fal m<lh_{i}(\eta)=lh_{i}(\gam)-1[\eta_{i}^{m}=\gam_{i}^{1+m}]$. By (\ref{eq:lem2.9}) we have 
$\gam=\gam_{i}^{0}\prec_{i+1}\gam_{i}^{1}=\eta_{i}^{0}$ and hence 
$\bar{\calv}=\bar{\calx}*(\calv_{0}) \spand \bar{\gam}=\bar{\bet}*(\gam_{1})\spand \calv^{0}=\calv_{0}\prec_{i+1}\calx_{0}$.
 This means that $\calv$ is a second son $\calz$ of $\calx$ in $\calt_{N}^{i}$. Therefore $\calv\prec_{i}\calx$.
\\
{\bf Case3} $\eta=pd_{i}(\gam)\prec_{i}rg_{i}(\gam)$: Then we have $\gam_{i}^{0}=\gam$, 
$rg_{i}(\eta_{i}^{m-1})=rg_{i}(\gam)$ and $st_{i}(\eta_{i}^{m-1})>st_{i}(\gam)\spand
\fal k<lh_{i}(\eta)-m+1=lh_{i}(\gam)(k>0 \rarw \eta_{i}^{m-1+k}=\gam_{i}^{k})$ for some $m$ with $0<m\leq lh_{i}(\eta)-1=n$.
 In particular $\gam=\gam_{i}^{0}\prec_{i+1}\gam_{i}^{1}=\eta_{i}^{m}$.
 Hence $\bar{\calv}=\bar{\calx}\uarw m \spand \bar{\gam}=\bar{\bet}\uarw(m+1)*(\gam_{m}) \spand 
 \calv^{0}=\calv_{0}\prec_{i+1}\calv_{1}=\calx_{m}$. 
 This means $\calv$ is a third son $\calu$ of $\calx$ in $\calt_{N}^{i}$. Therefore $\calv\prec_{i}\calx$.

This completes a proof of (\ref{eq:pretree}), cf. Lemma \ref{lem:5etaMh}. 
Our proof is based on the fact that the ramification procedure that produces three sons $\caly$, $\calz$ and $\calu$ from 
$\calx$ imitates the decomposition procedure of the relation $\eta=pd_{i}(\gam)$ in terms of sequences 
$\{\gam_{i}^{m}\}$ and $\{\eta_{i}^{m}\}$, and the relation $\prec_{i+1}$ between them, cf. Subsection \ref{subsec:od5fine}. 

In particular if $\gam\prec_{2}\eta$, then $M_{2}^{\gam}\prec_{2}M_{2}^{\eta}$, i.e., every $P\in M_{2}^{\eta}$ is 
$\Pi_{2}$-reflecting on the class $M_{2}^{\gam}$. This means that every $P\in M_{2}^{\eta}$ is $\eta$-Mahlo. 

Next we show the existence of an $\eta$-Mahlo set.
Corresponding to the construction tree of the class $M^{\eta}=M_{2}^{\eta}\in\calt_{N}$ we can associate a tree 
$\{\eta(s): s\in Tree(\eta)\}$ of o.d.'s in $\{\bet<\pi:\eta\preceq\bet\}$ with its depth $N-2$. 
First for the empty sequence $\langle\rangle$ $\eta(\langle\rangle)=\eta$. 
For each nonleaf $s\in Tree(\eta)$ let $\{s_{m}:-1\leq m\leq n\}$ be sons of $s$ in $Tree(\eta)$ with 
$n=lh_{i}(\eta(s))-1\geq 0$ and $s_{m}=s*(m)$. Then $\eta(s_{-1})=\eta(s), \eta(s_{m})=(\eta(s))_{i}^{m}$ for 
$0\leq m\leq n$, where $N-1>i=dp(s)+2$ with the depth $dp(s)$ of $s$, e.g., $dp(\langle\rangle)=0$. 
$s$ is a leaf if $dp(s)=N-3$.

For each $s\in Tree(\eta)$ we associate a $\Pi_{N}$-class $\calx(s;\eta)\in\calt_{N}^{i}$ with $i=dp(s)+2$ as follows. 
For a leaf $s$ put $\calx(s;\eta)=M_{N-1}^{st_{N-1}(\eta(s))}(\mbox{L})\in\calt_{N}^{N-1}$. 
Suppose $s$ is a nonleaf node and let $\{s_{m}:-1\leq m\leq n\}$ be sons of $s$. 
Then put $\calx(s;\eta)=M_{i}^{\eta(s)}\in\calt_{N}^{i}$ with $i=dp(s)+2$. Thus $\calx(\langle\rangle;\eta)=M_{2}^{\eta}$.

Now we show
\beqn\label{eq:exiL}
\fal s\in Tree(\eta)[\mbox{L}\in\calx(s;\eta)]
\eeqn
by tree induction on $s\in Tree(\eta)$ as follows.

 For a leaf $s$ we have $\eta(s)<\pi$, and hence 
 $\mbox{L}\in\calx(s;\eta)=M_{N-1}^{st_{N-1}(\eta(s))}(\mbox{L})\in\calt_{N}^{N-1}$ by induction on o.d.'s 
 $st_{N-1}(\eta(s))<\veps_{\pi+1}$.
 
Suppose $s$ is a nonleaf node with $i=dp(s)+2<N-1$ and let $\{s_{m}:-1\leq m\leq n\}$ be sons of $s$. Then 
\[
\calx(s;\eta)=M_{i}^{\eta(s)}=M_{i}(M_{i}(\bar{\calx};\bar{\bet})\cap\calx(s_{0};\eta))\cap\calx(s_{-1};\eta)
\]
for $\bar{\calx}=(\calx(s_{1};\eta),\ldots,\calx(s_{n};\eta))$ and o.d.'s $\bet_{m}=st_{i}((\eta(s))_{i}^{m-1})\, (0<m\leq n)$.

By IH(=Induction Hypothesis) we have $\mbox{L}\in\bigcap\{\calx(s_{m};\eta): -1\leq m\leq n\}$.
If $n=0$, then $\calx(s;\eta)=M_{i}(\calx(s_{0};\eta))\cap\calx(s_{-1};\eta)$. Since $\calx(s_{m};\eta)$ are 
$\Pi_{N}$-classes, a $\Pi_{N}$-reflecting universe L reflects these classes, i.e., $\mbox{L}\in\calx(s;\eta)$ by IH. 
Next assume $n>0$. Then by IH we have $\mbox{L}\in\calx(s_{n};\eta)\cap\calx(s_{-1};\eta)$. 
Hence (\ref{eq:5aprl3ax}) yields $\mbox{L}\in\calx(s;\eta)$. This completes a proof of (\ref{eq:exiL}). 
In particular we have $\mbox{L}\in\calx(\langle\rangle;\eta)=M_{2}^{\eta}$. 
Once again by reflecting the $\Pi_{N}$-class $M_{2}^{\eta}$ we conclude $\mbox{L}\in M_{2}(M_{2}^{\eta})$. 
Consequently $M_{2}^{\eta}\cap\mbox{L}\neq\emptyset$. This shows the existence of a {\it set\/} in $M_{2}^{\eta}$, 
cf. Theorem \ref{th:pi11exist}.

Let us examine the above proofs of (\ref{eq:pretree}) and (\ref{eq:exiL}). 
First these are based on the fact $\eta_{i}^{m}<\pi$ for any $m<lh_{i}(\eta)$ and any $i$ with $2\leq i<N-1$. 
Our proof of (\ref{eq:pretree}) is based on the fact that the second sons $\calz$ in $\calt_{N}^{i}\, (2\leq i<N-1)$ is 
less than its father $\calx$ with respect to $\prec_{i}$: $\calz\prec_{i}\calx$. The fact is based on (\ref{eq:5aprl3bx}), 
i.e., on (\ref{eq:5aprl3b}). On the other side our proof of (\ref{eq:exiL}) is also based on the fact (\ref{eq:5aprl3ax}), 
i.e., on (\ref{eq:5aprl3a}). These two facts (\ref{eq:5aprl3a}) and (\ref{eq:5aprl3b}) follow from (\ref{eq:5aprl1}), which
in turn, is shown by induction on ordinals $\bar{\bet}$, i.e., by induction on ordinal diagrams 
$\bet_{m}=st_{i}((\eta(s))_{i}^{m-1})\, (0<m\leq n, 2\leq i<N-1)$. Therefore we have to restrict o.d.'s to ones in 
wellfounded parts with respect to o.d.'s $st_{i}(\eta)\, (2\leq i<N-1)$ in advance, cf. 
Definition \ref{df:5wfuv}.\ref{df:5wfuv.V*d} of $V^{*}(X)$.
 Otherwise we would be in a circle, for the aim of (\ref{eq:pretree}) and (\ref{eq:exiL}) is to show that the system 
 $Od(\Pi_{N})$ of o.d.'s is wellfounded.
 
Furthermore our proof of (\ref{eq:exiL}) for leaves $s$ is based on induction on o.d.'s $st_{N-1}(\eta(s))<\veps_{\pi+1}$. 
When we restirct o.d.'s to a suitable subclass, then we can show transfinite induction up to each $\alp<\veps_{\pi+1}$, cf. 
Lemma \ref{lem:id3wf19-1}.\ref{lem:id3wf19-1.4} in Subsection \ref{subsec:proopr}. 
In this way we can show (\ref{eq:exiL}) for {\it each\/} $\eta$ in the subclass.

Note, here, that $\eta(s)<\pi$. If we would have $\eta(s)=\pi$ and put, e.g., 
$\calx(s;\eta)=M_{N-1}^{\pi}=M_{N-1}^{\veps_{\pi+1}}(\mbox{L})$, then we would need to invoke induction up to 
$\veps_{\pi+1}+1$ in showing (\ref{eq:exiL}) for leaves $s$. If we would put $M_{N-1}^{\pi}=\{\mbox{L}\}$, then, again, 
we would need to invoke induction up to $\veps_{\pi+1}+1$ in showing (\ref{eq:pretree}) for the case $N-1$. 
Therefore $\eta(s)<\pi$ is desired, cf. (\ref{eq:5pred}) in Definition \ref{df:piN} in Section \ref{subsec:od5api}.

This ends a set-theoretic explanation of o.d.'s.

\subsection{The system $Od(\Pi_{N})$}\label{subsec:od5api}

In this subsection we define the subsystem $Od(\Pi_{N})\subset Od$ of ordinal diagrams.

For $\rho\in\cald^{Q}\cap Od(\Pi_{N})$, we define o.d.'s $rg_{i}(\rho), st_{i}(\rho), pd_{i}(\rho)$ and a pair $in_{i}(\rho)$ of o.d.'s for $2\leq i<N$ and a set $In(\sig)\incl\{i: 2\leq i<N\}$. $st_{i}(\rho)$ and $rg_{i}(\rho)$ may be undefined. In this case we denote $st_{i}(\rho)\uarw$ and $rg_{i}(\rho)\uarw$. Otherwise we denote $st_{i}(\rho)\darw$ and $rg_{i}(\rho)\darw$.

Using $pd_{i}(\rho)$ we define a relation $\alp\prec_i\bet$ and its reflexive closure $\alp\preceq_i\bet$ as follows.
\bdf
$\alp\prec_i\bet$ {\rm denotes the transitive closure of the relation} $\{(\alp,\bet): pd_i(\alp)=\bet\}${\rm , and}
$\alp\preceq_i\bet$ {\rm its reflexive closure}.
\edf

\bdf $Od(\Pi_{N})$\label{df:piN}
{\rm The system} $Od(\Pi_{N})$ {\rm of o.d.'s is obtained from} $Od$ {\rm by restricting the construction} $(\sig,q,\alp)\mapsto d_{\sig}^{q}\alp$ {\rm in Definition \ref{df:opi}.\ref{caldQ} as follows:}

{\rm Assume} $\alp \in Od(\Pi_{N}) \, \& \,  \sig\in\{\pi\}\cup\cald^{Q} \spand q=\ovl{j\kap\tau\nu}$, {\rm where} $q=\ovl{j\kap\tau\nu}$ {\rm denotes a sequence of} quadruples $j_{m}\kap_{m}\tau_{m}\nu_{m}\incl Od(\Pi_{N})$ {\rm of length} $l+1 \, (0 \leq l\leq N-1-j_0)$ {\rm such that} 
 \benu
 \item $2\leq j_0<j_1<\cdots <j_l=N-1,$
 \item 
 $\kap_l=\pi \spand \sig\preceq\kap_m \, (m\leq l)$.
 \item $\nu_{l}\in Od(\Pi_{N})$,
 \beqn\label{eq:stboundN-1}
 \sig=\pi \Rarw \nu_{l}\leq\alp
 \eeqn
 {\rm cf. Lemma \ref{lem:N4aro},}\\
 {\rm and}
 \beqn\label{eq:stbound}
 m<l \Rarw \nu_m<\kap_{m}^{+}
 \eeqn
 \item $\tau_0=\sig, \tau_m\in\{\pi\}\cup\cald^{Q}, \sig\preceq \tau_m\, (m\leq l)$
  {\rm and}
 \beqn\label{eq:5pred}
 \tau_{l}=\pi \Rarw \sig=\pi
 \eeqn
 {\rm Cf. Lemmata \ref{lem:5.3}.\ref{lem:5.3.5} and \ref{lem:5.4}.\ref{lem:5.4.9}.}
 \eenu
{\rm Put} $\rho:=d^q_\sig\alp\in \cald^{Q}_{\sig}\incl Od$. {\rm For this} $\rho$ {\rm define}
 \benu
 \item $in_j(\rho)=st_j(\rho)rg_j(\rho)$ {\rm (a pair) and} $pd_j(\rho)$: {\rm Given} $j$ {\rm with} $2\leq j<N$, {\rm put} $m=\min\{m\leq l:j\leq j_m\}$.
  \benu
  \item $pd_j(\rho)=\tau_m$.
  \item  $\exi m\leq l(j=j_m)$: {\rm Then} $st_j(\rho)=\nu_m, \, rg_j(\rho)=\kap_m$.
  \item {\rm Otherwise:} $in_j(\rho)=in_j(pd_j(\rho))=in_j(\tau_m)$. {\rm If} $in_j(\tau_m)=\emptyset$, {\rm then set} $st_j(\rho)\uarw, rg_j(\rho)\uarw$.
  \eenu
 \item $In(\rho)=\{j_m:m\leq l\}$.
 \eenu
{\rm Then}
$\rho \in Od(\Pi_{N})$ {\rm if the following conditions are fulfilled besides (\ref{eq:Odmu}) in Definition \ref{df:opi}.\ref{caldQ}:}
\bdes
\item[$(\cald.1)$] {\rm Assume} $i\in In(\rho)$. {\rm Put} $\kap=rg_i(\rho)$. {\rm Then} \label{cnd:stn}
\item[$(\cald.11)$] $in_i(rg_i(\rho))=in_i(pd_{i+1}(\rho))$, $rg_i(\rho)\preceq_{i}pd_{i+1}(\rho)$ {\rm and} 
 $pd_i(\rho)\neq pd_{i+1}(\rho)$ {\rm if} $i<N-1$. \\
{\rm Also} $pd_i(\rho)\preceq_{i}rg_i(\rho)$ {\rm for any} $i$.\label{cnd:stn.1}
\\
{\rm Cf. Lemma \ref{lem:5.4}.\ref{lem:5.4.9}.}
\item[$(\cald.12)$]{\rm One of the following holds:} \label{cnd:stn.2}
\item[$(\cald.12.1)$] $rg_i(\rho)=pd_i(\rho)
  \spand \calb_{>\kap}(st_i(\rho))<b(\alp_{1})$ {\rm with} $\rho\preceq\alp_{1}\in\cald_{\kap}$. \label{cnd:stn.21}
  {\rm Cf. Lemma \ref{lem:5.4}.\ref{lem:5.4.3}.}
\item[$(\cald.12.2)$]  $rg_i(\rho)=rg_i(pd_i(\rho)) \spand st_i(\rho)<st_i(pd_i(\rho))$.
 \label{cnd:stn.22}
\item[$(\cald.12.3)$]  $rg_i(pd_i(\rho))\prec_i\kap \spand 
\fal \tau(rg_i(pd_i(\rho))\preceq _i\tau\prec _i\kap\rarw rg_i(\tau)\preceq _i\kap) \spand 
st_i(\rho)<st_i(\sig_1)$ {\rm with}
 \[\sig_1=\min\{\sig_1:rg_i(\sig_1)=\kap\spand pd_i(\rho)\prec _i\sig_1\prec _i\kap\}\]
 {\rm and such a} $\sig_1$ {\rm exists.} \label{cnd:stn.23}
 
 {\rm Cf. Lemma \ref{lem:5.4}.\ref{lem:5.4.8}.}
\item [$(\cald.2)$]
   \beqn \label{cnd:Kst}
   \fal\kap\leq rg_i(\rho)(K_{\kap}st_i(\rho)<\rho)
   \eeqn
    {\rm for} $i\in In(\rho)$.
\edes
\edf

 Also note that $\alp\prec_2\bet\Lrarw \alp\prec \bet$ for $\alp,\bet\in\cald^{Q}$.

In this subsection $X,Y,\ldots$ ranges over subsets of $Od(\Pi_{N})$. Ordinal diagrams in $Od(\Pi_{N})$ are denoted $\alp,\bet,\gam,\ldots$, while $\sig,\tau,\ldots$ denote o.d.'s in the set $(R\cap Od(\Pi_{N}))\cup\{\infty\}$.

\blem\label{lem:5ast3}
For $\alp,\bet\in Od(\Pi_{N})\cap\cald^{Q}$
\[\alp\prec_{N-1}\bet\Rarw st_{N-1}(\alp)<st_{N-1}(\bet)\]
\elem
\bprf This follows from the condition $(\cald.12.2)$. Note that for any $\alp\in\cald^{Q}$, $N-1\in In(\alp)\spand rg_{N-1}(\alp)=\pi$.
\eprf

We establish elementary facts on the relations $\prec_{i}$.

\blem\label{lem:5.3}
\benu
\item The finite set $\{\tau:\sig\prec_{i}\tau\}$ is linearly ordered by $\prec_{i}$.
\label{lem:5.3.1}
\item $\rho \prec_{i}pd_{i+1}(\rho)$, i.e., $\prec_{i+1}\incl\prec _i$. Also for $i<j$, $\rho \prec _{i}pd_{j}(\rho)$, i.e., $\prec_{j}\incl\prec_i$.
\label{lem:5.3.2}
\item 
$[i,j)\cap In(\rho)=\emptyset \spand i<j \Rarw \rho\prec_{j}pd_i(\rho)=pd_j(\rho)$.
\label{lem:5.3.3}
\item $in_i(\rho)=in_i(pd_i(\rho))\: \Lrarw \: i\not\in In(\rho)$. 
\label{lem:5.3.4}
\item For each $\eta\in\cald^{Q}$ and $i\in[2,N-1]$, $\max\{\eta_{\pi}<\pi: \eta\preceq_{i}\eta_{\pi}\}$ is the diagram $\eta_{\pi}$ such that $\eta\preceq\eta_{\pi}\in\cald_{\pi}$. 
Therefore $\alp\prec\bet\in\cald_{\pi} \Rarw \alp\prec_{i}\bet$ for any $i\in[2,N-1]$.
\label{lem:5.3.5}
\item Given $\gam,\kap$ so that $\exi\sig(\gam\preceq_{i}\sig\spand rg_{i}(\sig)=\kap)$, put
\\
$\sig=\max\{\sig:\gam\preceq_{i}\sig\spand rg_{i}(\sig)=\kap\}$. Then $i\in In(\sig)\spand \kap=pd_{i}(\sig)$. 
\label{lem:5.3.6}
\eenu
\elem
\bprf 
\\
\ref{lem:5.3}.\ref{lem:5.3.2}. This follows from the condition $(\cald.11)$, $pd_i(\rho)\preceq_{i}pd_{i+1}(\rho)$.
\\
\ref{lem:5.3}.\ref{lem:5.3.4}. By the definition we have the direction $[\Larw]$. For $[\Rarw]$ assume $i\in In(\rho)$. Then by the condition $(\cald.12)$ we have  $in_i(\rho)\neq in_i(pd_i(\rho))$.
\\
\ref{lem:5.3}.\ref{lem:5.3.5}. This is seen from the condition (\ref{eq:5pred}) and $pd_i(\rho)\preceq_{i}pd_{i+1}(\rho)$.
\\
\ref{lem:5.3}.\ref{lem:5.3.6}. By the maximality of $\sig$ we have $i\in In(\sig)$. In the condition $(\cald.12)$, \\
the latter two subcases $(\cald.12.2)$, $\tau=rg_i(\sig)=rg_i(pd_i(\sig))$ and \\
$(\cald.12.3)$, $\exi\sig_1(rg_i(\sig_1)=\tau\spand pd_i(\sig)\prec_i\sig_1)$ are not the cases again by the maximality of $\sig$. Hence the first subcase $(\cald.12.1)$, $\tau=pd_i(\sig)$ must occur.
\eprf

\bdf\label{df:alppi}{\rm Cf. Lemma \ref{lem:5.3}.\ref{lem:5.3.5}.)}\\
{\rm For each} $\eta\in\cald^{Q}$, $\eta_{\pi}$ {\rm denotes  the diagram} $\eta_{\pi}$ {\rm such that}
 $\eta\preceq\eta_{\pi}\in\cald_{\pi}$. 
\edf

\blem Assume $\kap=rg_i(\rho)\darw$.
\label{lem:5.4}
\benu
\item $\rho\prec_{i}rg_i(\rho)$.
\label{lem:5.4.1}
\item $st_{i}(\rho)<\kap^{+}$.
\label{lem:5.4.2}
\item $\rho\prec_i\tau \spand rg_i(\tau)\darw=\kap \spand i\in In(\rho) \Rarw 
st_i(\rho)<st_i(\tau)$.
\label{lem:5.4.5}
\item $pd_{i}(\rho)=\pi \Rarw \rho\in\cald_{\pi}\spand In(\rho)=\{N-1\}$.
Therefore $i<N-1 \Rarw rg_{i}(\rho)\leq pd_{i+1}(\rho)<\pi$.
\label{lem:5.4.9}
\item For a $j>i$, if $i\in In(\rho)$ and $rg_{j}(\rho)\darw$, then $rg_{i}(\rho)\preceq_{i}rg_{j}(\rho)$.
\label{lem:5.4.10}
\item $st_i(\rho)\leq\max\{b(\alp_{1}),\pi\} \spand \calb_{>\kap}(st_i(\rho))<b(\alp_{1})$ {\rm with} 
$\rho\preceq \alp_1\in\cald_\kap$. 
In fact we have $i<N-1 \Rarw st_{i}(\rho)<\pi$ and $st_{N-1}(\rho)\leq b(\alp_{1})$.
\label{lem:5.4.3}
\item $\rho\prec_i\sig\prec_i\tau \spand in_i(\rho)=in_i(\tau) \Rarw 
in_i(\rho)=in_i(\sig)$.
\label{lem:5.4.6}
\item
$N-1>i\in In(\rho) \spand rg_i(\rho)\preceq_{i}\del\prec_{i}pd_{i+1}(\rho) \Rarw i\not\in In(\del)$.
Therefore $rg_i(\rho)\preceq_{i+1}pd_{i+1}(\rho)$.
\label{lem:5.4.7}
\item $\rho\prec_i\tau<rg_i(\rho) \Rarw rg_i(\tau)\preceq_{i}rg_i(\rho)$.
\label{lem:5.4.8}
\eenu
\elem
{\bf Proof} by induction on $\ell\rho$.
\\
\ref{lem:5.4}.\ref{lem:5.4.1}. If $i\in In(\rho)$, then $\rho\prec_{i}pd_i(\rho)\preceq_{i}rg_i(\rho)$ by the condition $(\cald.11)$. Otherwise $pd_i(\rho)\prec_{i}rg_i(pd_i(\rho))= rg_i(\rho)$ by IH.
\\
\ref{lem:5.4}.\ref{lem:5.4.2}. This follows from the condition (\ref{eq:stbound}) in Definition \ref{df:piN} and the convention $\pi^{+}=\infty$ for $\pi=rg_{N-1}(\rho)$.
\\
\ref{lem:5.4}.\ref{lem:5.4.5}. By Lemma \ref{lem:5.4}.\ref{lem:5.4.1} we have 
$\rho\prec_i\tau\prec_{i}rg_i(\tau)=\kap=rg_i(\rho)$. Thus by the condition $(\cald.12)$ one of the following cases occur:
\bdes
\item[$(\cald.12.2)$] $rg_i(\rho)=rg_{i}(pd_i(\rho))\spand pd_{i}(\rho)\preceq_{i}\tau$: Then by IH $st_i(\rho)<st_i(pd_i(\rho))\leq st_i(\tau)$,

\hspace*{-10mm}
or 
\item[$(\cald.12.3)$] $rg_i(pd_i(\rho))\prec_i\kap$: Then $rg_i(pd_i(\rho))\neq \kap=rg_i(\tau)$ and hence 
$pd_i(\rho)\prec_i\tau$. Put $\sig_1=\min\{\sig_1:rg_i(\sig_1)=\kap\spand pd_i(\rho)\prec_i\sig_1\prec_i\kap\}$. Then $\sig_1\preceq_{i}\tau$. Therefore by IH $st_i(\rho)<st_i(\sig_1)\leq st_i(\tau)$.
\edes
\ref{lem:5.4}.\ref{lem:5.4.7}. 
This is seen from Lemmata \ref{lem:5.4}.\ref{lem:5.4.6} and \ref{lem:5.3}.\ref{lem:5.3.4}.
\\
\ref{lem:5.4}.\ref{lem:5.4.9}. This is seen from the conditions (\ref{eq:5pred}) and $(\cald.11)$.
\\
\ref{lem:5.4}.\ref{lem:5.4.10}. This is seen from the definition $(\cald.11)$ and Lemma \ref{lem:5.3}.\ref{lem:5.3.2}.
\\
\ref{lem:5.4}.\ref{lem:5.4.3}. By $\kap=rg_i(\rho)\darw$ we have 
$\exi\sig(\rho\preceq_{i}\sig\spand i\in In(\sig)\spand rg_{i}(\sig)=\kap)$. Let $\sig$ denote the maximal $\sig$ such that
$\rho\preceq_{i}\sig\spand i\in In(\sig)\spand rg_{i}(\sig)=\kap$. Then by Lemma \ref{lem:5.3}.\ref{lem:5.3.6} 
$\kap=pd_{i}(\sig)=rg_{i}(\sig)$. Hence by Lemmata \ref{lem:5.4}.\ref{lem:5.4.5} and \ref{lem:5.4}.\ref{lem:5.4.2}
 we have $st_i(\rho)\leq st_{i}(\sig)<\kap^{+}$. 
If $i<N-1$, then $\kap=rg_{i}(\sig)<\pi$ by Lemma \ref{lem:5.4}.\ref{lem:5.4.9}.
Otherwise $\kap=rg_{i}(\sig)=\pi\spand \sig=\alp_{1}\in\cald_{\pi}$ by 
  Lemma \ref{lem:5.3}.\ref{lem:5.3.6}. Hence by the condition (\ref{eq:stboundN-1}) in Definition \ref{df:piN} we have 
  $st_{i}(\sig)\leq b(\alp_{1})$. Therefore $st_i(\rho)\leq\max\{b(\alp_{1}),\pi\}$.

It remains to show $\calb_{>\kap}(st_i(\rho))<b(\alp_{1})$. Lemma \ref{lem:Od3} with \\
$st_i(\rho)\leq st_{i}(\sig)<\kap^{+}$ yields $\fal\tau>\kap[\calb_{\tau}(st_i(\rho))\leq\calb_{\tau}(st_{i}(\sig))]$, and hence $\calb_{>\kap}(st_i(\rho))\leq\calb_{>\kap}(st_{i}(\sig))$.
On the other hand we have $\calb_{>\kap}(st_{i}(\sig))<b(\alp_{1})$ by the condition $(\cald.12.1)$ in Definition \ref{df:piN}. Consequently $\calb_{>\kap}(st_i(\rho))<b(\alp_{1})$.
\\
\ref{lem:5.4}.\ref{lem:5.4.6}. This follows from Lemma \ref{lem:5.4}.\ref{lem:5.4.5} and IH.
\\
\ref{lem:5.4}.\ref{lem:5.4.8}. 
First note that if $\rho\prec_i\tau<rg_i(\rho)$, then $\tau\prec_{i}rg_i(\rho)$ by 
Lemmata \ref{lem:5.4}.\ref{lem:5.4.1} and  \ref{lem:5.3}.\ref{lem:5.3.1}.

By IH and the condition $(\cald.12)$ we can assume $i\in In(\rho)$ and the case $(\cald.12.3)$, 
$rg_i(pd_i(\rho))\prec_i\kap$ occurs. Thus also assume $pd_i(\rho)\prec_i\tau$. 
If $rg_i(pd_i(\rho))\preceq _{i}\tau\prec_i\kap$, then the condition $(\cald.12.3)$ yields $rg_i(\tau)\preceq_i\kap$. 
So assume $pd_i(\rho)\prec_{i}\tau\prec_{i}rg_i(pd_i(\rho))$. 
Then by IH we have $rg_i(\tau)\preceq_{i} rg_i(pd_i(\rho))\prec_i\kap$.
\eprf

\blem\label{lem:N4aro}
For $\rho\in\cald\cap Od(\Pi_{N})$, $Q(\rho)\leq\max\{b(\rho),\pi\}$.
\elem
\bprf
By Definition \ref{df:piN} and (\ref{eq:stbound}) we have $In(\rho)\cup\{pd_{i}(\rho), rg_{i}(\rho):i\in In(\rho)\}\cup\{st_{i}(\rho): i\in In(\rho)\spand i<N-1\}\leq\pi$. On the other hand we have $st_{N-1}(\rho)\leq b(\alp_{1})\leq b(\rho)$ by Lemmata \ref{lem:5.4}.\ref{lem:5.4.3} and \ref{lem:Npi11exist} for the diagram $\alp_{1}$ with $\rho\preceq\alp_{1}\in\cald_{\pi}$.
\eprf

\subsection{A finer analysis of the relations $\prec_{i}$}\label{subsec:od5fine}
In this subsection we give a finer analysis of the relation $\alp\prec_{i}\bet$. This is needed in Sections \ref{subsec:wfpiNid} and \ref{sec:5awf}.

First for each $\eta\in\cald^{Q}$ and each $i\in[2,N-1)=\{i\in\ome:2\leq i<N-1\}$ define a length $lh_{i}(\eta)$ and a sequence $\{\eta^{n}_{i}:n<lh_{i}(\eta)\}\incl\{\del<\pi:\eta\preceq\del\}$ of subdiagrams of $\eta$. The sequence decomposes the sequence $\{\del<\pi: \eta\preceq_{i+1}\del\}$.

\bdf\label{df:5res}Length $lh_{i}(\eta)$ and a sequence $\{\eta^{n}_{i}:n<lh_{i}(\eta)\}$ of subdiagrams of $\eta\in\cald^{Q}$
\\
{\rm {\bf Case \ref{df:5res}.1}}. $\neg\exi\del(\eta\preceq_{i}\del\spand i\in In(\del))${\rm : Then put} $lh_{i}(\eta)=1$ {\rm and} $\eta_{i}^{0}:=\eta_{\pi}${\rm , cf. Definition \ref{df:alppi}. 
Namely} $\eta_{i}^{0}$ {\rm denotes the maximal diagram such that} $\eta\preceq_{i+1}\eta_{i}^{0}<\pi${\rm .}
\\
{\rm {\bf Case \ref{df:5res}.2}}. $\exi\del(\eta\preceq_{i}\del\spand i\in In(\del))${\rm : Then} $\eta_{i}^{0}$ {\rm is defined to be the minimal diagram such that} $\eta\preceq_{i}\eta_{i}^{0}\spand i\in In(\eta_{i}^{0})${\rm .}

{\rm Suppose that} $\eta_{i}^{n}$ {\rm is defined so that} $i\in In(\eta_{i}^{n})$.
\\
{\rm {\bf Subcase \ref{df:5res}.2.1}}. $\exi\gam(rg_{i}(\eta_{i}^{n})\preceq_{i}\gam\spand i\in In(\gam))${\rm : Then} $\eta_{i}^{n+1}$ {\rm is defined to be the minimal diagram such that} $rg_{i}(\eta_{i}^{n})\preceq_{i}\eta_{i}^{n+1}\spand i\in In(\eta_{i}^{n+1})$.
\\
{\rm {\bf Subcase \ref{df:5res}.2.2}}{\rm . Otherwise: Then} $lh_{i}(\eta)=n+2$ {\rm and define} $\eta_{i}^{n+1}$ {\rm to be the maximal diagram such that} $\eta_{i}^{n}\preceq_{i+1}\eta_{i}^{n+1}<\pi${\rm , i.e., the diagram such that} $\eta\preceq\eta_{i}^{n+1}\in\cald_{\pi}$.
\edf

\blem\label{lem:5ap12}
\[\eta_{i}^{n}=\gam_{i}^{m} \Rarw 
\fal k<lh_{i}(\eta)-n=lh_{i}(\gam)-m\{\eta_{i}^{n+k}=\gam_{i}^{m+k}\}\]
and
\[rg_{i}(\eta_{i}^{n})=rg_{i}(\gam_{i}^{m}) \Rarw 
\fal k<lh_{i}(\eta)-n=lh_{i}(\gam)-m\{k>0\rarw \eta_{i}^{n+k}=\gam_{i}^{m+k}\}.\]
\elem 
\bprf
The first assertion is clear. Assume $rg_{i}(\eta_{i}^{n})=rg_{i}(\gam_{i}^{m})$. Then $\eta_{i}^{n+1}=\gam_{i}^{m+1}$, and hence the second assertion follows from the first one.
\eprf

\blem\label{lem:5Si-2} For $i<N-1$, 
\[ 
\fal\del[rg_i(\eta_{i}^{n})\preceq_{i}\del\prec_{i}pd_{i+1}(\eta_{i}^{n}) \Rarw i\not\in In(\del)] 
\spand \eta\preceq_{i+1}\eta_{i}^{0} \spand \fal n<lh_{i}(\eta)-1[\eta_{i}^{n}\prec_{i+1}\eta_{i}^{n+1}].\]
\elem
\bprf
First we show $\eta\preceq_{i+1}\eta_{i}^{0}$. By the definition and Lemma \ref{lem:5.3}.\ref{lem:5.3.5} 
we can assume that $\eta_{i}^{0}$ is the minimal diagram 
such that $\eta\preceq_{i}\eta_{i}^{0}\spand i\in In(\eta_{i}^{0})$, i.e, {\bf Case \ref{df:5res}.2}. 
Then $i\not\in In(\del)$ for any $\del$ with $\eta\preceq_{i}\del\prec_{i}\eta_{i}^{{0}}$, and hence the assertion follows from 
Lemma \ref{lem:5.3}.\ref{lem:5.3.3}.

Next we show $\fal\del[rg_i(\eta_{i}^{n})\preceq_{i}\del\prec_{i}pd_{i+1}(\eta_{i}^{n}) \Rarw i\not\in In(\del)]$. 
By the condition $(\cald.11)$, 
$in_i(rg_i(\eta_{i}^{n}))=in_i(pd_{i+1}(\eta_{i}^{n}))\spand rg_i(\eta_{i}^{n})\preceq_{i}pd_{i+1}(\eta_{i}^{n})$, 
Lemma \ref{lem:5.4}.\ref{lem:5.4.7} we have 
\[\fal\del[rg_i(\eta_{i}^{n})\preceq_{i}\del\prec_{i}pd_{i+1}(\eta_{i}^{n})
 \Rarw i\not\in In(\del)].\]

Finally we show $\eta_{i}^{n}\prec_{i+1}\eta_{i}^{n+1}$.
We can assume, by Lemma \ref{lem:5.3}.\ref{lem:5.3.5}, that $\eta_{i}^{n+1}$ is the minimal diagram such that 
$rg_{i}(\eta_{i}^{n})\preceq_{i}\eta_{i}^{n+1}\spand i\in In(\eta_{i}^{n+1})$, i.e., {\bf Subcase \ref{df:5res}.2.1}. 
We have by the definition that $\eta_{i}^{n+1}$ is the diagram such that 
$rg_i(\eta_{i}^{n})\preceq_{i}\eta_{i}^{n+1}$ and
\[\fal\del[rg_i(\eta_{i}^{n})\preceq_{i}\del\prec_{i}\eta_{i}^{n+1} \Rarw i\not\in In(\del)] \spand 
i\in In(\eta_{i}^{n+1}).\]
Therefore $pd_{i+1}(\eta_{i}^{n})\preceq_{i}\eta_{i}^{n+1}$ and 
$\fal\del[pd_{i+1}(\eta_{i}^{n})\preceq_{i}\del\prec_{i}\eta_{i}^{n+1}
 \Rarw i\not\in In(\del)]$. Hence $\eta_{i}^{n}\prec_{i+1}pd_{i+1}(\eta_{i}^{n})\preceq_{i+1}\eta_{i}^{n+1}$.
\eprf

\blem\label{lem:5Si-1}
\[\gam\prec_{i}\eta\prec_{i}\kap=rg_{i}(\gam)\darw \, \Rarw 
\exi m<lh_{i}(\eta)-1[\kap=rg_{i}(\eta_{i}^{m})].\]
\elem
{\bf Proof} by induction on $\ell\eta$.
Put 
\[\sig=\max\{\sig:\gam\preceq_{i}\sig\prec_{i}\kap\spand rg_{i}(\sig)=\kap\}.\]
Then by Lemma \ref{lem:5.3}.\ref{lem:5.3.6} we have $i\in In(\sig)\spand \kap=pd_{i}(\sig)$. Hence $\eta\preceq_{i}\sig$. If $\eta=\sig$, then $\kap=rg_{i}(\eta_{i}^{0})$ with $\eta_{i}^{0}=\eta$. Assume $\eta\prec_{i}\sig$. If $i\not\in In(\eta)$, then we have $\gam\prec_{i}pd_{i}(\eta)\preceq_{i}\sig\prec_{i}\kap$. IH with $(pd_{i}(\eta))_{i}^{m}=\eta_{i}^{m}$ yields the assertion. Suppose $i\in In(\eta)$. By Lemma \ref{lem:5.4}.\ref{lem:5.4.8} we have $rg_{i}(\eta)\preceq_{i}\kap$. If $rg_{i}(\eta)=\kap$, then we are done by $\eta_{i}^{0}=\eta$. Suppose $rg_{i}(\eta)\prec_{i}\kap$. Then we have $\eta_{i}^{1}\preceq_{i}\sig$ by the definition, and hence $\gam\prec_{i}\eta_{i}^{1}\prec_{i}\kap$. IH with $\eta_{i}^{m}=(\eta_{i}^{1})_{i}^{m-1}$ yields the assertion.
\eprf

\blem\label{lem:3.23.1}
Assume $\eta=pd_{i}(\gam)$ for an $i<N-1$. Then one of the following holds, cf. Subsection \ref{sec:prl5a.2}:
\bdes
\item[Case \ref{lem:3.23.1}.1]
$\eta=pd_{i+1}(\gam) \spand \gam_{i}^{0}=\eta_{i}^{0}$. Hence
\[
\eta=pd_{i}(\gam)=pd_{i+1}(\gam) \spand lh_{i}(\gam)=lh_{i}(\eta)\spand \fal m<lh_{i}(\gam)[\gam_{i}^{m}=\eta_{i}^{m}]
\]
\item[Case \ref{lem:3.23.1}.2]
$rg_{i}(\gam)=pd_{i}(\gam)=\eta \spand \gam_{i}^{0}=\gam \spand \gam_{i}^{1}=\eta_{i}^{0}$. Hence
\[
rg_{i}(\gam)=pd_{i}(\gam)=\eta \spand \gam_{i}^{0}=\gam \spand \fal m<lh_{i}(\eta)=lh_{i}(\gam)-1[\eta_{i}^{m}=\gam_{i}^{1+m}]
\]
\item[Case \ref{lem:3.23.1}.3]
$\eta=pd_{i}(\gam)\prec_{i}rg_{i}(\gam) \spand \gam_{i}^{0}=\gam$ and $rg_{i}(\eta_{i}^{m})=rg_{i}(\gam)\spand st_{i}(\eta_{i}^{m})>st_{i}(\gam)$ for an $m<lh_{i}(\eta)-1$. Hence
\beqnarrs
&& \eta=pd_{i}(\gam)\prec_{i}rg_{i}(\gam) \spand \gam_{i}^{0}=\gam \spand  \\
&& \exi m[0<m\leq lh_{i}(\eta)-1 \spand rg_{i}(\eta_{i}^{m-1})=rg_{i}(\gam)\spand st_{i}(\eta_{i}^{m-1})>st_{i}(\gam)\spand \\
&& \fal k<lh_{i}(\eta)-m+1=lh_{i}(\gam)(k>0 \rarw \eta_{i}^{m-1+k}=\gam_{i}^{k})]
\eeqnarrs
\edes
\elem
\bprf
Assume $\eta=pd_{i}(\gam)$ for an $i<N-1$. 

First consider the case $i\not\in In(\gam)$. Then by the definition $\eta=pd_{i+1}(\gam) \spand \gam_{i}^{0}=\eta_{i}^{0}$ holds. Hence by Lemma \ref{lem:5ap12} {\bf Case \ref{lem:3.23.1}.1} holds. 

In what follws suppose $i\in In(\gam)$. Then $\gam_{i}^{0}=\gam$ and $\eta=pd_{i}(\gam)\preceq_{i}rg_{i}(\gam)$ by Lemma \ref{lem:5.4}.\ref{lem:5.4.1}, i.e., by the condition $(\cald.11)$.
Second consider the case $rg_{i}(\gam)=pd_{i}(\gam)=\eta$. Then by the definition we have $\gam_{i}^{1}=\eta_{i}^{0}$.
Therefore by Lemma \ref{lem:5ap12} {\bf Case \ref{lem:3.23.1}.2} holds. 

Finally consider the case $\eta=pd_{i}(\gam)\prec_{i}rg_{i}(\gam)$. Then by Lemma \ref{lem:5Si-1} we have $rg_{i}(\eta_{i}^{m})=rg_{i}(\gam)\spand st_{i}(\eta_{i}^{m})>st_{i}(\gam)$ for an $m<lh_{i}(\eta)-1$.
Consequently by Lemma \ref{lem:5ap12} {\bf Case \ref{lem:3.23.1}.3} holds. 
\eprf

\bdf\label{df:free}
\benu
\item\label{df:free.1}
{\rm For} $2\leq i<N-1$ {\rm define}
\[
\alp\lhd_{i}\bet :\Lrarw \alp,\bet\in\cald^{Q} \spand i\in In(\alp) \spand \alp\prec_{i}\bet\prec_{i}rg_{i}(\alp).
\]
\item\label{df:free.2}
{\rm For} $2\leq i<N-1$ {\rm define}
\[
\alp\lhd^{i}\bet :\Lrarw \alp,\bet\in\cald^{Q} \spand i\in In(\alp) \spand 
rg_{i}(\alp)\preceq_{i}\bet\prec_{i}pd_{i+1}(\alp).
\]
\eenu
\edf

The following lemma is seen from Lemma \ref{lem:3.23.1}.

\blem\label{lem:3.23.1+}
Assume $\gam\prec_{i}\eta$ for an $i<N-1$. Then one of the following holds
\bdes
\item[Case \ref{lem:3.23.1+}.1] 
$\gam_{i}^{0}=\eta_{i}^{0}\spand \gam\prec_{i+1}\eta$.
\item[Case \ref{lem:3.23.1+}.2]
$\exi n\in(0,lh_{i}(\gam))[\gam_{i}^{n}=\eta_{i}^{0}]$, and $\gam_{i}^{n-1}\lhd^{i}\eta$.
\item[Case \ref{lem:3.23.1+}.3]
$\exi n\in[0,lh_{i}(\gam)-1)\exi m\in[0,lh_{i}(\eta)-1)[\gam_{i}^{n}\prec_{i}\eta\preceq_{i}\eta_{i}^{m}\spand rg_{i}(\gam_{i}^{n})=rg_{i}(\eta_{i}^{m})]$, and $\gam_{i}^{n}\lhd_{i}\eta$.
\edes
\elem

\blem\label{lem:rgpilhd} 
$\alp\in\cald_{\pi}\spand i<N-1 \Rarw \gam\not\!\!\lhd_{i}\alp$.
\elem
\bprf
This follows from Lemma \ref{lem:5.4}.\ref{lem:5.4.9}.
\eprf

\subsection{Decomposing ordinal diagrams}\label{subsec:daggareast}

In this subsection we introduce decompositions $\alp(s)$ of ordinal diagrams $\alp$, where 
$s$ denotes a function in ${}^{[i,k)}2\, (2\leq i\leq k\leq N-3)$.
In the next section \ref{subsec:wfpiNid} we define a suitable $\Pi^{0}_{N-1}$-operator $\Gam_{N}$ 
through the decompositions.

\bdf\label{df:ppd}
{\rm For an o.d.'s} $\eta\in\cald^{Q}, i<N-1$ {\rm define as follows:}
\benu
\item
\[
ppd_{i}(\eta):=\left\{
\begin{array}{ll}
rg_{i}(\eta), & \mbox{{\rm if }} i\in In(\eta) \\
pd_{i}(\eta), & \mbox{{\rm if }} i\not\in In(\eta)
\end{array}
\right.
\]
\item
\[
\gam\prec^{p}_{i}\alp :\Lrarw \gam\prec_{i}\alp \spand \neg\exi n<lh_{i}(\gam)[\gam^{n}_{i}\lhd_{i}\alp]
\]
{\rm and} $\gam\preceq^{p}_{i}\alp$ {\rm denotes the reflexive closure of the relation} $\gam\prec^{p}_{i}\alp$.
\eenu
\edf

\blem\label{lem:precp}
\benu
\item\label{lem:precp.1}
$\gam\prec^{p}_{i}\alp$ is the transitive closure of the relation $\{(\alp,\bet): \bet=ppd_{i}(\alp)\}$,
and hence $\prec_{i}^{p}$ is transitive and $\prec_{i+1}\incl\prec^{p}_{i}\incl\prec_{i}$.
\item\label{lem:precp.2}
$\alp\prec_{i}\bet \Rarw \exi\gam[\alp\preceq^{p}_{i}\gam \spand (ppd_{i}(\gam)=\bet \lor 
\gam\lhd_{i}\bet)]$.
\item\label{lem:precp.3}
$\alp\preceq_{i}^{p}\bet \Rarw \alp\preceq_{i+1}\bet_{i}^{0}$.
\eenu
\elem
\bprf
\\
\ref{lem:precp}.\ref{lem:precp.1}. 
This is seen from Lemmata \ref{lem:5Si-2}, \ref{lem:3.23.1+} and \ref{lem:5.4}.\ref{lem:5.4.1}.
\\
\ref{lem:precp}.\ref{lem:precp.2}. This is seen from Lemma \ref{lem:precp}.\ref{lem:precp.1}.
\\
\ref{lem:precp}.\ref{lem:precp.3}.
By induction on $\ell\alp$ and $\alp\preceq_{i+1}\alp_{i}^{0}$ we can assume that $\bet=ppd_{i}(\alp)$.

If $i\not\in In(\alp)$, then $\bet=pd_{i}(\alp)=pd_{i+1}(\alp)$, and hence 
$\alp\prec_{i+1}\bet\preceq_{i+1}\bet_{i}^{0}$.

Suppose $i\in In(\alp)$ and $\bet=rg_{i}(\alp)$. Then $\alp=\alp_{i}^{0} \spand \bet_{i}^{0}=\alp_{i}^{1}$.
Hence Lemma \ref{lem:5Si-2} yields $\alp\prec_{i+1}\bet_{i}^{0}$.
\eprf

\bdf\label{df:barst}
{\rm For} $\alp\in\cald^{Q}$ {\rm let} $\ovl{st}_{N-1}(\alp)$ {\rm denote the pair}
\[
\ovl{st}_{N-1}(\alp):=
\langle st_{N-1}(\alp_{N-2}^{0}), st_{N-1}(\alp)\rangle
\]
\edf

\blem\label{lem:piNbarstb}
For any $\alp\in\cald^{Q}$ and $\alp\preceq\alp_{\pi}\in\cald_{\pi}$,
$st_{N-1}(\alp)\leq st_{N-1}(\alp_{N-2}^{0})\leq b(\alp_{\pi})$.
\elem
\bprf
By Lemma \ref{lem:5Si-2} we have $\alp\preceq_{N-1}\alp_{N-2}^{0}$. Hence Lemma \ref{lem:5ast3} yields 
$st_{N-1}(\alp)\leq st_{N-1}(\alp_{N-2}^{0})$.
$st_{N-1}(\alp_{N-2}^{0})\leq b(\alp_{\pi})$ is seen from \ref{lem:5.4}.\ref{lem:5.4.3}.
\eprf

\blem\label{lem:barst}
$\gam\prec^{p}_{N-2}\alp
\Rarw \ovl{st}_{N-1}(\gam)<_{lex}\ovl{st}_{N-1}(\alp)$.
\elem
\bprf
Assume $\gam\prec^{p}_{N-2}\alp$.
Then by Lemma \ref{lem:3.23.1+} one of the following holds:
\\
{\bf Case \ref{lem:3.23.1+}.1} $\gam_{N-2}^{0}=\alp_{N-2}^{0}\spand \gam\prec_{N-1}\alp$: 
Then $st_{N-1}(\gam_{N-2}^{0})=st_{N-1}(\alp_{N-2}^{0})$. 
Lemma \ref{lem:5ast3} yields $st_{N-1}(\gam)<st_{N-1}(\alp)$. Hence $\ovl{st}_{N-1}(\gam)<_{lex}\ovl{st}_{N-1}(\alp)$.
\\
{\bf Case \ref{lem:3.23.1+}.2} $\exi n\in(0,lh_{N-2}(\gam))[\gam_{N-2}^{n}=\alp_{N-2}^{0}]$:  
By $n>0$ we have $\gam_{N-2}^{0}\prec_{N-1}\gam_{N-2}^{n}$ by Lemma \ref{lem:5Si-2}, and hence 
$st_{N-1}(\gam_{N-2}^{0})<st_{N-1}(\gam_{N-2}^{n})=st_{N-1}(\alp_{N-2}^{0})$ by Lemma \ref{lem:5ast3}. 
Therefore $\ovl{st}_{N-1}(\gam)<_{lex}\ovl{st}_{N-1}(\alp)$.
\\
{\bf Case \ref{lem:3.23.1+}.3} $\exi n\in[0,lh_{N-2}(\gam)-1)\exi m\in[0,lh_{N-2}(\alp)-1)[\gam_{N-2}^{n}\prec_{N-2}\alp\preceq_{N-2}\alp_{N-2}^{m} \spand rg_{N-2}(\gam_{N-2}^{n})=rg_{N-2}(\alp_{N-2}^{m})]$: Then $\gam^{n}_{N-2}\lhd_{N-2}\alp$. Hence this is not the case.
\eprf

\bdf\label{df:sequence}
$I:=I_{N}:=\bigcup\{{}^{[i,k)}2: 2\leq i\leq k\leq N-3\}$
 {\rm denotes the set of functions from the set} $[i,k)=\{j\in\ome: i\leq j<k\}$ 
{\rm to} $2=\{0,1\}$.
\benu
\item
{\rm For} $s\in{}^{[i,k)}2$ {\rm let}
\beqnarrs
d(s) & := & i \\
\ell(s) & := &  k \\
\#s & := & \#\{j\in[i,k): s(j)=1\}
\eeqnarrs

{\rm Note that there are} $(N-4)$ {\rm empty functions in}  $I${\rm .}
{\rm Each element in}  $I$ {\rm is, by definition, a triple of a function} $s$ {\rm and} $d(s), \ell(s)${\rm .
An empty function is, then, a triple} $(\emptyset, i,i)\, (2\leq i\leq N-3)$.

{\rm If} $\# s=0${\rm , then} $s$ {\rm is said to be} null.

\item {\rm For} $s\in I$, $s|i$ {\rm denotes the function in} $I$ {\rm such that}
$d(s|i)=d(s)$, $\ell(s|i)=\min\{\max(i,d(s)),\ell(s)\}$ {\rm and}
$(s|i)(j):= s(j)$ {\rm for} $d(s|i)\leq j<\ell(s|i)${\rm .}

\item {\rm For} $s,t\in I$, $s<_{lex}t$ {\rm denotes the lexicographic ordering induced by} $0<1${\rm :}
\[
s<_{lex}t :\Lrarw \exi i\in[d(s),\ell(s))\cap[d(t),\ell(t))\{s|i=t|i \spand s(i)=0<1=t(i)\}.
\]
{\rm Note that} $s<_{lex}t  \Rarw d(s)=d(t)$.

\item
{\rm For} $s,t\in I$ {\rm with} $\ell(s)=d(t)${\rm ,} $u=s*t$ {\rm denotes the concatenated sequence.
Namely} $d(u)=d(s), \ell(u)=\ell(t)$ {\rm and} $u(i)=s(i)$ {\rm for} $i\in[d(s),\ell(s))$, 
$u(i)=t(i)$ {\rm for} $i\in[d(t),\ell(t))$.

\item
$s\in I$ {\rm is said to be} unitary {\rm if}
\[
\fal i\in[d(s),\ell(s))\{s(i)=1 \Rarw i=\ell(s)-1\} \spand \fal j\in[d(s),\ell(s)-1)\{s(j)=0\}.
\]

\item
{\rm Each} $s\in I$ {\rm is decomposed uniquely to the concatenated sequence of longest unitary
components} $s_{i}$, $s=s_{0}*\cdots *s_{k}\, (k\geq 0)$ {\rm such that}
$\{d(s)\}\cup \{i+1\in(d(s),\ell(s)]: s(i)=1\}\cup\{\ell(s)\}=\{\ell_{-1}<\ell_{0}<\cdots <\ell_{k}\}\, (k\geq 0)$ {\rm and each}
$s_{i}\in{}^{[\ell_{i-1},\ell_{i})}2$ {\rm is a subseries of} $s$.

$s=s_{0}*\cdots *s_{k}$ {\rm is said to be the} unitary decomposition {\rm of} $s$.

\item
{\rm If} $d(s)=2${\rm , then} $s$ {\rm is said to be} initial.

$I(2):=\{s\in I: d(s)=2\}$ {\rm denotes the set of initial sequences in} $I$.

\item
$I(2,N-2):=\{s\in I: d(s)=2\spand \ell(s)=N-2\}$.

\item
$t\incl_{e}s$ {\rm [}$t\subset_{e}s${\rm ] designates that} $t$ {\rm is an [a proper]} initial segment 
{\rm of} $s$
{\rm , i.e.,} $\exi i\leq \ell(s)\{t=s|i\}$ {\rm [}$\exi i< \ell(s)\{t=s|i\}${\rm ], resp.}

\item {\rm For} $s,t\in I$ {\rm with} $d(s)=d(t)$, 
$s\cap t$ {\rm denotes a sequence in} $I$ {\rm defined as follows:}
\[
(s\cap t)(i)=j :\Lrarw \fal k<i[s(k)=t(k)] \spand s(i)=t(i)=j.
\]
$d(s\cap t)=d(s)$ {\rm  and} $\ell(s\cap t)=\min(\{i: s(i)\neq t(i)\}\cup\{\ell(s),\ell(t)\})$.

\eenu
\edf

\bdf\label{df:ast}
{\rm Let} $\alp\in\cald^{Q}${\rm . We define ordinal diagrams} $\alp(s)\, (s\in I)$ {\rm as follows:}
\benu
\item \label{df:ast.0}
{\rm First define} $\alp(s)$ {\rm for unitary} $s$ {\rm by induction on} $\ell\alp$ {\rm as follows:}
 \benu
 \item $\alp(s):=\alp$ {\rm if} $s$ {\rm is null}.
 
 \item {\rm Otherwise} $s(\ell(s)-1)=1\spand \fal i\in[d(s),\ell(s)-1)\{s(i)=0\}$.
   \bdes
   \item[Case 1] {\rm If}
   \beqn\label{eq:()}
   \fal i\in[d(s),\ell(s)-1)\{\alp_{\ell(s)-1}^{0}<\alp_{i}^{0} \spand 
   rg_{\ell(s)-1}(\alp)\leq\alp_{i}^{0}\}
   \eeqn
   {\rm then put}
   \[
   \alp(s):=\min\{\del: \alp_{\ell(s)-1}^{0}\preceq_{\ell(s)}^{p}\del<\alp^{\prime}(s)\},
   \]
   {\rm where}
   \[
   \alp^{\prime}=\min(\{\alp_{i}^{0}: i\in[d(s),\ell(s)-1)\}\cup\{\alp_{\ell(s)-1}^{1}\}).
   \]

   \item[Case 2] {\rm Otherwise: Then put, cf. Definition \ref{df:alppi},}
   \[
   \alp(s):=\alp_{\pi}.
   \]
   \edes
  \eenu
 
 \item
 $\alp(s)$ {\rm is defined through the unitary decomposition} $s=s_{0}*s_{1}*\cdots * s_{k}$ {\rm as follows}

\[
\alp(s):=(\cdots (\alp(s_{0}))(s_{1})\cdots)(s_{k}).
\]
\eenu
\edf

\bdf\label{df:ast[]}
{\rm Let} $\alp,\bet\in\cald^{Q}${\rm .}
\benu
\item

{\rm Let} $2\leq k\leq N-3$.
{\rm Define} $s[k;\alp,\bet]\in {}^{[k,N-3]}2$ {\rm recursively. Suppose that} $s=(s[k;\alp,\bet])|i$ 
{\rm has been defined for an} $i$ {\rm with} $k\leq i\leq N-3${\rm . Then}
$s[k;\alp,\bet](i)\in\{0,1\}$ {\rm is defined as follows:}
\beqn\label{eq:s[]}
s[k;\alp,\bet](i):=\left\{
        \begin{array}{ll}
         1 & \mbox{{\rm if }} \exi\gam[\alp(s)\preceq_{i}\gam\prec_{i}\bet(s) \spand i\in In(\gam)] \\
        0 & \mbox{{\rm otheriwse}}
        \end{array}
        \right.
\eeqn

\item
\[
t\incl s[\alp,\bet] :\Lrarw \fal i\in[d(t),\ell(t))\{t(i)=s[d(t);\alp,\bet](i)\}.
\]
\eenu
\edf

\blem\label{lem:astup}
\benu
\item \label{lem:astup0}
$\alp(s|i)=\bet(s|i) \Rarw \alp(s)=\bet(s)$.

\item \label{lem:astup1}
Let $s$ be unitary with $s(\ell(s)-1)=1$.
Assume $\alp\prec_{d(s)}\bet$ and $\fal i\in[d(s),\ell(s))[\alp_{i}^{0}\geq\bet]$.
Then $\alp(s)=\bet(s)$.

\item \label{lem:astup1.3}
Let $s$ be unitary with $s(\ell(s)-1)=1$. 
Then $\alp\preceq_{\ell(s)}\alp_{\ell(s)-1}^{0}\preceq^{p}_{\ell(s)}\alp(s)$.

\item \label{lem:astup1.4}
Let $\alp\prec\del\prec\gam$.
For an initial $t$, assume $t\not\incl_{e}s[2;\alp,\del]$, $t\incl_{e}s[2;\alp,\gam]$ and 
$\del\prec_{t|(\ell(t)-1)}\gam$.

Then $t(i)=1$ for $i=\min\{i: t(i)\neq s[2;\alp,\del](i)\}$, and hence
$$s[2;\alp,\del]<_{lex}t.$$

\item \label{lem:astup1.5}
$u=s*t \spand \ell(s)\in In(\alp(s)) \spand rg_{\ell(u)}(\alp(u))\darw \Rarw 
rg_{\ell(s)}(\alp(s))\preceq rg_{\ell(u)}(\alp(u))$.

\eenu
\elem
\bprf
\\
\ref{lem:astup}.\ref{lem:astup0}.
Let $s=(s|i)*t$. If $t$ is null, then $\alp(s)=\alp(s|i)=\bet(s|i)=\bet(s)$.
Otherwise let $s=s_{0}*\cdots*s_{k}$ be the unitary decomposition of $s$, and
$s_{j}=s'*t'$ with $s|i=s_{0}*\cdots*s_{j-1}*s'$ and $t=t'*s_{j+1}*\cdots * s_{k}$.
Then $\alp(s|i)=\alp(s_{0}*\cdots*s_{j-1})$ and similarly for $\bet(s|i)$.
By Definition \ref{df:ast}.\ref{df:ast.0} we have $\alp(s_{0}*\cdots*s_{j})=\bet(s_{0}*\cdots*s_{j})$,
and we see inductively $\alp(s)=\bet(s)$.
\\
\ref{lem:astup}.\ref{lem:astup1}.
We have $\fal i\in[d(s),\ell(s))[\alp\prec_{i}\bet]$ by Lemma \ref{lem:5.3}.\ref{lem:5.3.3}.
Also we have $\fal i\in[d(s),\ell(s))[\alp_{i}^{0}=\bet_{i}^{0}$ and $\alp_{\ell(s)-1}^{1}=\bet_{\ell(s)-1}^{1}$.
By Definition \ref{df:ast}.\ref{df:ast.0} we have $\alp(s)=\bet(s)$.
\\
\ref{lem:astup}.\ref{lem:astup1.3}.
This is seen from Lemmata \ref{lem:5.3}.\ref{lem:5.3.5}, \ref{lem:5.4}.\ref{lem:5.4.9}  and \ref{lem:5Si-2}.
\\
\ref{lem:astup}.\ref{lem:astup1.4}.
Suppose $t(i)=0$. Then $s[2;\alp,\del](i)=1$. 
By the minimality of $i$ we have $(s[2;\alp,\gam])|i=t|i=(s[2;\alp,\del])|i$, and hence $\alp(t|i)\preceq_{i}\eta\prec_{i}\del(t|i)$
for an $\eta$ by the definition (\ref{eq:s[]}).
On the other hand we have $\del(t|i)\prec_{i}\gam(t|i)$ by $\del\prec_{t|(\ell(t)-1)}\gam$.
Therefore $t(i)=s[2;\alp,\gam](i)=1$. A contradiction.
\\
\ref{lem:astup}.\ref{lem:astup1.5}.
If $\alp(u)=\alp(s)$, then the assertion follows from the assumption $\ell(s)\in In(\alp(s))$
and Lemma \ref{lem:5.4}.\ref{lem:5.4.10}.

Suppose $\alp(u)>\alp(s)$. 
Then $t$  is not null.
Let $t_{1}$ be the longest unitary subseries of $t$ such that for some $t_{0}, t_{2}$,
$t=t_{0}*t_{1}*t_{2}$, $\alp(s)=\alp(s*t_{0})$ and $\alp(s)<\alp(s*t_{0}*t_{1})$.

$\alp(s*t_{0}*t_{1})$ is defined by the {\bf Case 1} in Definition \ref{df:ast}.
Otherwise we would have $\alp(u)=\alp(s*t_{0}*t_{1})=\alp_{\pi}$, and $rg_{\ell(u)}(\alp(u))\uarw$.
Therefore we have $\alp(s)=\alp(s*t_{0})\prec_{\ell(t_{1})}\alp(s*t_{0}*t_{1})$
with $\ell(s)<\ell(t_{1})$.
Hence $rg_{\ell(s)}(\alp(s))\preceq\alp(s*t_{0}*t_{1})\preceq\alp(u)\prec rg_{\ell(u)}(\alp(u))$.
\eprf

\bdf\label{df:precs}
{\rm Let} $\alp\in\cald^{Q}$ {\rm and} $s\in I$.
\benu
\item\label{df:precs1}
\[
\alp\prec_{s}\bet :\Lrarw s\incl s[\alp,\bet] \spand
 \fal t\incl_{e} s[\alp(t)\prec_{\ell(t)}^{p}\bet(t)].
\]

\item\label{df:precs2}
\[
\alp\prec^{-}_{s}\bet :\Lrarw s\incl s[\alp,\bet] \spand 
\fal t\incl_{e}s [\# t<\#  (u*s) \Rarw \alp(t)\prec_{\ell(t)}^{p}\bet(t)].
\]

\item\label{df:precs3}
\[
\alp\lhd_{s}\bet :\Lrarw s\incl s[\alp,\bet] \spand 
(\ell(s)>d(s) \Rarw \alp\prec_{s|(\ell(s)-1)}\bet) \spand \alp(s)\lhd_{\ell(s)}\bet(s),
\]
{\rm cf. Definition \ref{df:free}.\ref{df:free.1}.}

\item\label{df:precs4}
$\alp\lhd_{s}^{+}\bet$ {\rm iff} $\alp\lhd_{s}\bet$ {\rm and for each} $i\in[d(s),\ell(s))$
{\rm if} $s(i)=1 \spand i\not\in In(\alp(s|i))${\rm , then there are sequences} $\{\alp_{k}\}_{k\leq K}$ 
{\rm of diagrams and a sequence} $\{v_{k}\}_{k<K}\incl I$
{\rm such that}
\beqn\label{eq:C+11}
 \alp_{K}=\alp(s|i) \spand \alp_{0}=(\alp_{K})_{i}^{0} 
\spand
\fal k<K\{\alp_{k+1}\lhd_{v_{k}}\alp_{k}  \spand d(v_{k})=i\}
\eeqn
\eenu
\edf

\blem\label{lem:astup11}
$
\alp\prec_{t}\bet \spand t<_{lex}u \Rarw \alp(u)=\bet(u)
$.
\elem
\bprf
This is seen from Lemma \ref{lem:astup}.\ref{lem:astup1}.
\eprf

\blem\label{clm:lemB}
Suppose that $u=t \lor u<_{lex}t \lor t<_{lex}u$.
\benu
\item\label{clm:lemB1}
Suppose 
$\del\lhd_{u}\gam\lhd_{t}\bet$.
Then  $\del\lhd_{s}\bet$ for $s=\max_{<_{lex}}\{u,t\}$.

\item\label{clm:lemB3}
Suppose that
$\del\lhd_{u}^{+}\gam\lhd_{t}^{+}\bet$.
Then  $\del\lhd_{s}^{+}\bet$  for $s=\max_{<_{lex}}\{u,t\}$.
\eenu
\elem
\bprf
\\
\ref{clm:lemB}.\ref{clm:lemB1}.
This is seen from Lemma \ref{lem:astup11}.
\\
\ref{clm:lemB}.\ref{clm:lemB3}.
By Lemma \ref{clm:lemB}.\ref{clm:lemB1} we have $\del\lhd_{s}\bet$.

If either $u=t$ or $t<_{lex}u$, then $\del\lhd_{u}^{+}\bet$ is seen from the definition.

In what follows assume $u<_{lex}t$.

Let $j$ denote the number such that $u|j=t|j \spand u(j)=0<1=t(j)$.
Suppose $t(i)=1 \spand i\not\in In(\del(t|i))$. 
We have to show that there are sequences 
enjoying the condition (\ref{eq:C+11}).

If $i<j$, then $\del\lhd_{u}^{+}\gam$ yields the assertion.
If $i>j$, then $\del(t|i)=\gam(t|i)$ by Lemma \ref{lem:astup11}, and hence $\gam\lhd_{t}^{+}\bet$ yields the assertion.
Finally assume $i=j$.
We have $(\del(u|i))_{i}^{0}=(\gam(u|i))_{i}^{0}$.
If $i\not\in In(\gam(t|i))$, then pick sequences $\{\alp_{k}\}_{k\leq K}$ and $\{v_{k}\}_{k<K}$
for $\gam(t|i)=\alp_{K}$ and $i$.
Otherwise set $K=0$.
Now let $\alp_{K+1}=\del(t|i)$ and $(t|i)*v_{K}=u$.
Sequences $\{\alp_{k}\}_{k\leq K+1}$ and $\{v_{k}\}_{k<K+1}$ are desired one for $\del, t$ and $i$.
\eprf

\blem\label{lem:E}
\benu
\item\label{lem:E.0}
Let $s$ be unitary with $s\incl s[\alp,\bet]$.
Assume $[\ell(s)>d(s) \Rarw \alp\prec_{s|(\ell(s)-1)}\bet]$ and 
$\alp(s)\prec^{p}_{\ell(s)}\gam^{\prime}\prec_{\ell(s)}\bet(s)$.
Then there exists a $\gam$ such that $\gam(s)=\gam^{\prime}$, 
$\alp\prec_{s}\gam$, $[\ell(s)>d(s) \Rarw \gam\prec_{s|(\ell(s)-1)}\bet]$ and
\beqn\label{eq:lemE}
\fal j\in[d(s),\ell(s))\{s(j)=1 \lrarw j\in In(\gam(s|j))\}
\eeqn
\item\label{lem:E.1}
Assume $\alp\prec_{s}\bet \spand \alp(s)\prec^{p}_{\ell(s)}\gam^{\prime}\prec^{p}_{\ell(s)}\bet(s)$.
Then there exists a $\gam$ such that $\gam(s)=\gam^{\prime}$ and
$\alp\prec_{s}\gam\prec_{s}\bet$.
\item\label{lem:E.2}
Assume $[\ell(s)>d(s) \Rarw \alp\prec_{s|(\ell(s)-1)}\bet]$ and 
$\alp(s)\prec^{p}_{\ell(s)}\gam^{\prime}\lhd_{\ell(s)}\bet(s)$.
Then there exists a $\gam$ such that $\gam(s)=\gam^{\prime}$, 
$\alp\prec_{s}\gam$ and $\gam\lhd_{s}\bet$.
\eenu
\elem
\bprf 
\\
\ref{lem:E}.\ref{lem:E.0} by induction on $\ell\alp$.
Suppose $s$ is  unitary.
If $s$ is  null, then $\gam=\gam^{\prime}$ works.

Suppose $s$ is not null, and put $i=\ell(s)-1$.
By IH we can assume that $\gam^{\prime}=ppd_{i+1}(\alp(s))$.
Then $\fal j\in[d(s),i)\fal\del\{\alp\preceq_{j}\del\prec_{j}\bet  \Rarw j\not\in In(\del)\}$ and 
$\alp\prec^{p}_{i}\bet$ with $\exi\del\{\alp\preceq_{i}\del\prec_{i}\bet  \spand i\in In(\del)\}$, i.e.,
$rg_{i}(\alp_{i}^{0})\preceq_{i}\bet\leq\alp_{j}^{0}$ for any  $j\in[d(s),i)$.
Therefore the condition (\ref{eq:()}) in Definition \ref{df:ast}.\ref{df:ast.0} holds.
Hence
\[
   \alp(s)=\min\{\del: \alp_{i}^{0}\preceq_{i+1}^{p}\del<\alp^{\prime}(s)\},
 \]
where
\[
   \alp^{\prime}=\min(\{\alp_{j}^{0}: j\in[d(s),i)\}\cup\{\alp_{i}^{1}\}).
\]
If $\alp^{\prime}(s)=\alp_{\pi}$, then we would have 
$\gam^{\prime}=ppd_{i+1}(\alp(s))=\alp_{\pi}\not\prec_{i+1}\bet(s)$.
Hence $\alp^{\prime}(s)<\alp_{\pi}$, and  
$\alp^{\prime}=\alp_{i}^{1}$.
This means that
\beqn\label{eq:lemE1}
\alp\prec_{s|i}\alp^{\prime}
\eeqn

We  have  $\alp_{i}^{0}\prec_{i+1}\alp_{i}^{1}=\alp^{\prime}\preceq_{i+1}\alp^{\prime}(s)$
by Lemma \ref{lem:5Si-2}, and hence either $ppd_{i+1}(\alp(s))=\alp^{\prime}(s)$ or
$\alp(s)\lhd_{i+1}\alp^{\prime}(s)$ by Lemma \ref{lem:precp}.\ref{lem:precp.2}.
If $\gam^{\prime}=ppd_{i+1}(\alp(s))=\alp^{\prime}(s)$, then $\gam=\alp^{\prime}$ works:
If $\bet\leq\alp^{\prime}$, then we would have by (\ref{eq:lemE1})  and 
Lemma \ref{lem:astup11},
$\bet(s)\in\{(\alp_{i}^{m}(s):m\leq 1\}\leq\alp^{\prime}(s)=\gam^{\prime}$.
Hence $\alp^{\prime}<\bet$, which yields $\gam=\alp_{i}^{0}\prec_{s|(\ell(s)-1)}\bet$.
On the other hand we have, for (\ref{eq:lemE}), $i\in In(\alp_{i}^{1})$.
Otherwise we would have $\bet>\alp_{i}^{1}=\alp_{\pi}$.

Suppose $\alp(s)\lhd_{i+1}\alp^{\prime}(s)$.
We have $\alp_{i}^{0}\prec^{p}_{i}rg_{i+1}(\alp(s))=ppd_{i+1}(\alp(s))=\gam^{\prime}$ by 
Lemma \ref{lem:precp}.\ref{lem:precp.1}.
Let $\del=\max\{\del: \alp_{i}^{1}\preceq_{s|i}\del \spand \alp(s)\lhd_{i+1}\del(s)\}$.
Then $rg_{i+1}(\del(s))\preceq_{i+1}rg_{i+1}(\alp(s))=\gam^{\prime}<\bet(s)\leq\alp_{\pi}$
 by Lemma \ref{lem:5.4}.\ref{lem:5.4.8},
and hence $\del=\alp_{i}^{k}$ for a  $k$ by the maximality of $\del$ and $\del\preceq_{s|i}\del_{i}^{0}$.

We have
$ \del^{\prime}=\min(\{\del_{j}^{0}: j\in[d(s),i)\}\cup\{\alp_{i}^{k+1}\})=\alp_{i}^{k+1}$
and
$\del(s)=\min\{\eta: \alp_{i}^{k}\preceq_{i+1}^{p}\eta<(\alp_{i}^{k+1})(s)\}$.
Therefore 
\beqn\label{eq:lemE2}
\alp\prec_{s|i}\alp_{i}^{1}\prec_{s|i}\alp_{i}^{k+1}
\eeqn

We have either $(\alp_{i}^{k})(s)\lhd_{i+1}(\alp_{i}^{k+1})(s)$ or
$ppd_{i+1}((\alp_{i}^{k})(s))=(\alp_{i}^{k+1})(s)$.
By the maximality of  $\del$, we have $ppd_{i+1}((\alp_{i}^{k})(s))=(\alp_{i}^{k+1})(s)$,
and $\gam^{\prime}=rg_{i+1}(\alp(s))=(\alp_{i}^{k+1})(s)$.
Thus $\gam=\alp_{i}^{k+1}$ works:
Suppose $\bet\leq\alp_{i}^{k+1}$.
Then by (\ref{eq:lemE2}) and Lemma \ref{lem:astup11} we would have 
$\bet(s)\in\{(\alp_{i}^{m})(s):m\leq k+1\}\leq(\alp_{i}^{k+1})(s)=\gam^{\prime}$.
Hence $\alp_{i}^{k+1}<\bet$, which yields $\gam=\alp_{i}^{k+1}\prec_{s|(\ell(s)-1)}\bet$.
On the other hand we have, for (\ref{eq:lemE}), $i\in In(\alp_{i}^{k+1})$.
Otherwise we would have $\bet>\alp_{i}^{k+1}=\alp_{\pi}$.
\\
\ref{lem:E}.\ref{lem:E.1} by induction on $\ell(s)$.
The unitary case follows from Lemma \ref{lem:E}.\ref{lem:E.0}.

Next let  $s=s_{0}*\cdots*s_{k-1}*s_{k}$ be the  unitary decomposition of $s$, and
$t=s_{0}*\cdots*s_{k-1}$. Then 
$\alp\prec_{s}\bet \Lrarw \alp\prec_{t}\bet \spand \alp(t)\prec_{s_{k}}\bet(t)$.
Then $\alp(s)=(\alp(t))(s_{k})$ and similarly for $\bet(s)$.
By the unitary case pick a  $\del$ such that 
$\alp(t)\prec_{s_{k}}\del\prec_{s_{k}}\bet(t) \spand \del(s_{k})=\gam^{\prime}$.
By  IH pick a $\gam$ such that $\alp\prec_{t}\gam\prec_{t}\bet \spand \gam(t)=\del$.
\\
\ref{lem:E}.\ref{lem:E.2}.
The unitary case follows again from Lemma \ref{lem:E}.\ref{lem:E.0}.

Next let  $s=s_{0}*\cdots*s_{k-1}*s_{k}$ be the  unitary decomposition of $s$, and
$t=s_{0}*\cdots*s_{k-1}$. 
Then $t\incl_{e}s|(\ell(s)-1)$, and, by Lemma  \ref{lem:E}.\ref{lem:E.1}, 
 it suffices to find a $\del$ such that
$\alp(t)\prec_{s_{k}}\del\lhd_{s_{k}}\bet(t)$ and $\del(s_{k})=\gam^{\prime}$.
This is seen from Lemma \ref{lem:E}.\ref{lem:E.0}.
\eprf

\blem\label{lem:B}
Assume $\alp\prec\bet$ and $\alp\not\prec_{s}\bet$ for an $s\incl_{e}s[2;\alp,\bet]$.
Then there exist a $\gam$ and a $t$ such that $t\incl_{e}s$
and $\alp\preceq_{t}\gam\lhd_{t}^{+}\bet$.

Moreover let $\gam$ denote the minimal diagram for which the above conditions hold 
for some $t$.
Then there  are no $\del$ and no  $u$ such that $u<_{lex}t$
and
$\alp\preceq_{t}\del\lhd_{u}^{+}\gam$.
\elem
\bprf
We can assume that $\ell(s)>2 \Rarw \alp\prec_{s^{-}}\bet$ for $s^{-}=s|i$ with  $i=\ell(s)-1$, and $\alp(s)\not\prec^{p}_{\ell(s)}\bet(s)$.
We  show $t=s$ works for a $\gam$.

First consider the case when $\ell(s)=2$, i.e., $s$ is empty.
Then $\alp\prec\bet$ while  $\alp\not\prec^{p}_{2}\bet$.
This means that $\alp\preceq_{2}^{p}\gam\lhd_{2}\bet$ for a $\gam$.
In what follows assume $\ell(s)>2$.

Second consider the case when $s(i)=0$.
Then $\alp(s)=\alp(s^{-})\prec^{p}_{i}\bet(s)$ while $\alp(s)\not\prec^{p}_{i+1}\bet(s)$.
On the other side $s(i)=0$ yields $\alp(s)\prec_{i+1}\bet(s)$.
Hence by Lemma \ref{lem:precp}.\ref{lem:precp.2} there exists a $\gam^{\prime}$ such that
$\alp(s)\preceq^{p}_{i+1}\gam^{\prime}\lhd_{i+1}\bet(s)$.
By Lemma \ref{lem:E}.\ref{lem:E.2} pick a $\gam$ so that
$\gam(s)=\gam^{\prime}$ and $\alp\preceq_{s}\gam\lhd_{s}^{+}\bet$.

Third consider the case when $s(i)=1$.
Let $s=s_{0}*\cdots s_{k-1}*v=u*v$ be the unitary decomposition of $s$.
Then $\alp(u)=\alp(s^{-})\prec^{p}_{i}\bet(u)$ while $\alp(s)\not\prec^{p}_{i+1}\bet(s)$.
From $\alp\prec_{s^{-}}\bet$ and $s(i)=1$, we see that
$(\alp(u))_{i}^{0}<\bet(u)$ and
$\alp(s)=\min\{\del: (\alp(u))_{i}^{0}\preceq_{i+1}^{p}\del<(\alp(u))_{i}^{1})(v)\}$.

If $\bet(s)=\alp_{\pi}$, then we would have $\alp(s)\prec^{p}_{i+1}\bet(s)$.
Hence $\bet(s)<\alp_{\pi}$.

We have $(\alp(u))_{i}^{0}\prec_{i+1}(\alp(u))_{i}^{1}\preceq_{i+1}(\bet(u))_{i}^{0}\preceq_{i+1}^{p}\bet(s)$.
Therefore $\alp(s)\prec_{i+1}$.
Hence by Lemma \ref{lem:precp}.\ref{lem:precp.2} there exists a $\gam^{\prime}$ such that
$\alp(s)\preceq^{p}_{i+1}\gam^{\prime}\lhd_{i+1}\bet(s)$.
By Lemma \ref{lem:E}.\ref{lem:E.2} pick a $\gam$ so that
$\gam(s)=\gam^{\prime}$ and $\alp\preceq_{s}\gam\lhd_{s}^{+}\bet$.

Lemma \ref{clm:lemB} with the minimality of $\gam$ yields the last assertion of the lemma.
\eprf

\blem\label{lem:B'}
Assume $\alp\prec_{s}\gam$ and $\alp\preceq_{u}\eta\lhd_{u}\del\lhd_{t}\gam$ for a $u\incl_{e}s\cap t$.
Then $\# u<\# t$ if $t<_{lex}s \, \lor \, t\incl_{e}s$.
\elem
\bprf
If $u=t\incl_{e}s$, then $\alp(t)\preceq_{\ell(t)}\eta(t)\lhd_{\ell(t)}\del(t)\lhd_{\ell(t)}\gam(t)$, and hence
$\alp(t)\not\prec_{\ell(t)}\gam(t)$.
This is not the case.

Assume $u\subset_{e}t \spand u\subset_{e}s \spand \#u=\#t$. Then $t=u*v$ for a non empty null $v$.
We have $\alp\preceq_{u}\eta\lhd_{u}\del\prec_{u}\gam$ and $\alp\prec_{u}\gam$.
Hence 
\[
\eta(u)\prec_{\ell(u)}\del(u)\prec_{\ell(u)}rg_{\ell(u)}(\eta(u))\preceq_{\ell(u)}\gam(u).
\]
Therefore $\exi\xi[\del(u)\preceq_{\ell(u)}\xi\prec_{\ell(u)}\gam(u) \spand \ell(u)\in In(\xi)]$
by Lemma \ref{lem:5.3}.\ref{lem:5.3.6}.
This contradicts $t(\ell(u))=0$ and $t\incl s[\del,\gam]$.
\eprf

\blem\label{lem:C+}
Let $\alp\lhd_{s}^{+}\bet$.
\benu
\item\label{lem:C+1}
Assume $s(i)=1$ and $i\not\in In(\alp(s|i))$ for an $i\in[d(s),\ell(s))$.
Then there exist a sequence $\{\alp_{k}\}_{k\leq K}\, (K>0)$ of diagrams and 
a sequence $\{v_{k}\}_{k<K}\incl I$ enjoying (\ref{eq:C+11}) and the following condition (Cf. Definition \ref{df:precs}.\ref{df:precs4}.).
\[
 \mbox{ Each } v_{k} \mbox{ is null, and } \fal k<K-1\{v_{k}\subset_{e}v_{k+1}\}.
\]
\item\label{lem:C+2}
Assume $\# u\geq\# s$ and  $u<_{lex}s$.
Then $\alp(u)=\alp_{\pi}$ and $\lnot\exi \gam[\gam\lhd_{u}\alp]$.
\eenu
\elem
\bprf
\\
\ref{lem:C+}.\ref{lem:C+1}.
By the Definition \ref{df:precs}.\ref{df:precs4} let $K>0$ denote the {\it least\/} number
for which there exist sequences
 $\{\alp_{k}\}_{k\leq K}$ and $\{v_{k}\}_{k<K}$ enjoying (\ref{eq:C+11}).

Note that $v_{k}(i)=0$ since $v_{k}\incl s[\alp_{k+1},\alp_{k}]$ and 
$s[d(v_{k});\alp_{k+1},\alp_{k}](i)=0$ by 
$(\alp_{k+1})_{i}^{0}=\alp_{0}$.

\bclm\label{clm:lemC'1}
\benu
\item\label{clm:lemC'1.1}
Each $v_{k}$ is null.
\item\label{clm:lemC'1.2}
$\fal k\in(0,K)\{v_{k-1}\subset_{e}v_{k}\}$.
\eenu
\eclm
{\bf Proof} of Claim \ref{clm:lemC'1} by simultaneous induction on $k$.

Supposing  that $v_{k}$ is not null,
let $m_{k}=\min\{m\geq i: v_{k}(m_{k})=1\}$.
$m_{k}>i$ by $v_{k}(i)=0$.

Now $m_{0}>i\in In(\alp_{0})$ yields 
$\alp_{0}(v_{0}|(m_{0}+1))=(\alp_{0})_{\pi}=\alp_{0}(v_{0})\not\prec rg_{\ell(v_{0})}(\alp_{1}(v_{0}))$.
This is not the case, and $v_{0}$ is null.

Suppose $v_{k-1}$ is null, while $v_{k}$ is not null.
Then $v_{k-1}\subset_{e}v_{k}$ and $m_{k}\geq\ell(v_{k-1})$.
Since $v_{k-1}$ is null, $\alp_{k}(v_{k-1})=\alp_{k}$.

Hence $m_{k}>\ell(v_{k-1})\in In(\alp_{k})$ yields 
$\alp_{k}(v_{k}|(m_{k}+1))=(\alp_{k})_{\pi}=\alp_{k}(v_{k})\not\prec rg_{\ell(v_{k})}(\alp_{k+1}(v_{k}))$.
This is not the case.

Suppose $m_{k}=\ell(v_{k-1})$, and let $\ell(v_{k-2})=i$ if $k=1$.
Then by IH we have $In(\alp_{k-1})\ni\ell(v_{k-2})<\ell(v_{k-1})=m_{k}$, and
$(\alp_{k})^{0}_{\ell(v_{k-2})}=\alp_{k-1}=\alp_{k-1}(v_{k-1})\prec rg_{\ell(v_{k-1})}(\alp_{k}(v_{k-1}))=rg_{m_{k}}(\alp_{k})$.
Therefore $\alp_{k}(v_{k}|(m_{k}+1))=(\alp_{k})_{\pi}$.
Again this is not the case.
We have shown Claim \ref{clm:lemC'1}.\ref{clm:lemC'1.1}.

Suppose both $v_{k-1}$ and $v_{k}$ are null, and $v_{k-1}\not\subset_{e}v_{k}$.
Then we shown that sequences can be shortened contradicting the minimality of the number $K$.

Suppose either $v_{k-1}<_{lex}v_{k}$ or $v_{k}<_{lex}v_{k-1}$. Then $\alp_{k+1}\lhd_{v}\alp_{k-1}$
for $v=\max_{<_{lex}}\{v_{k-1},v_{k}\}$ by Lemma \ref{clm:lemB}.\ref{clm:lemB1}.

Suppose $v_{k}\subset_{e}v_{k-1}$.
Then $\alp_{k}=\alp_{k}(v_{k})\prec_{\ell(v_{k})}
rg_{\ell(v_{k})}(\alp_{k+1}(v_{k}))=rg_{\ell(v_{k})}(\alp_{k+1})$ and $v_{k-1}(\ell(v_{k}))=0$.
The latter means that there is no $\del$ such that 
$\alp_{k}\preceq_{\ell(v_{k})}\del\prec_{\ell(v_{k})}\alp_{k-1}$.
Therefore from Lemma \ref{lem:5.3}.\ref{lem:5.3.6} we see $\alp_{k+1}\lhd_{v_{k}}\alp_{k-1}$.
\\
\ref{lem:C+}.\ref{lem:C+2}.
Let $i$ denote the number such that $u(i)=0<1=s(i) \spand u|i=s|i$.
By $\#u\geq\#s$, let $j=\min\{j>i: u(j)=1\}$.
We show $\alp(u|(j+1))=\alp_{\pi}$. Then $\alp(u)=\alp(u|(j+1))=\alp_{\pi}$ follows.

Suppose $i\in In(\alp(s|i))$. Then the condition (\ref{eq:()}) in Definition \ref{df:ast} is broken.
Specifically $\alp(u|i)=(\alp(u|i))_{i}^{0}\leq(\alp(u|i))_{j}^{0}$,
and hence $\alp(u|(j+1))=\alp_{\pi}$.

In what follows assume $i\not\in In(\alp(s|i))$, and let 
$\{\alp_{k}\}_{k\leq K}$ and $\{v_{k}\}_{k<K}$ denote sequences in Lemma \ref{lem:C+}.\ref{lem:C+1}
 for $i$.

$j>\ell(v_{K-1})$
Let $\ell_{k}=\ell(v_{k})$ for $k\in[0,K)$ and $\ell_{-1}=i$.
Then we have $\ell_{k}>\ell_{k-1}$ for any $k\in[0,K)$.

First consider the case when $j=\ell_{k}$ for a $k\in[0,K)$.
Then $\ell_{k}>\ell_{k-1}\in In(\alp_{k})$ and 
$(\alp_{k+1})_{\ell_{k-1}}^{0}=\alp_{k}\prec rg_{\ell_{k}}(\alp_{k+1})$.
Hence $\alp(u|(\ell_{k}+1))=\alp_{\pi}$.

Otherwise there exists a maximal $k\in[-1,K)$ such that 
$\ell_{k}<j$.
Then $(\alp_{K})_{\ell_{k}}^{0}=\alp_{k+1}\leq(\alp_{K})_{j}^{0}$, and hence
$\alp(u|(j+1))=\alp_{\pi}$.
We are done.
\eprf

\section{Wellfoundedness proof for $Od(\Pi_{N})$ by means of inductive definitions}
\label{subsec:wfpiNid}

In this section we work in the theory $\Pi_{N-1}\mbox{-Fix}$ and show the

\bth\label{th:id5wfeachpiN}
For each $\alp<d_{\Ome}\veps_{\pi+1}$,i.e., each $\alp\in Od(\Pi_{N})|\Ome$, $\Pi^{0}_{N-1}\mbox{{\rm -Fix}}$ proves that $(Od(\Pi_{N})|\alp,<)$ is a well ordering.
\eth

\subsection{Operators $\calg_{i}$}\label{subsec:calgi}

Define operators $\calg_{i}\, (1\leq i<N-1)$ on $Od(\Pi_{N})$ recursively from the operator $\calg(X)$ in 
Definition \ref{df:G()}.

\bdf\label{df:oprtNG}{\rm Define inductively} $\Pi^{0}_{i}${\rm -operators} 
$\calg_{i}\, (1\leq i<N-1)$ {\rm as follows.}
\benu
\item\label{df:oprtNG.1}
\[
\calg_{1}(X):=\calg(X).
\]
\item\label{df:oprtNG.<}
{\rm For} $1<i\leq N-1$
\[
\calg_{<i}(X) :=\bigcap\{\calg_{j}(X): 1\leq j<i\}.
\]
\item\label{df:oprtNG.s}
{\rm For} $s\in I(2)$
\[
\alp\in\calg_{s}(X) :\Lrarw \fal\gam\lhd^{+}_{s}\alp[\gam\in\calg_{<2+\# s}(X) \rarw \gam\leq X|\alp]
.\]
\item\label{df:oprtNG.G}
{\rm For} $2\leq i<N-1$
\[
\calg_{i}:=\bigcap\{\calg_{s}(X): s\in I(2) \spand 2+\# s=i\}
.\]
\eenu
\edf

Now let us define an operator $\Gam_{N}$ on $Od(\Pi_{N})$ from these operators.

\bdf\label{df:oprtN}
\beqnarrs
&& \alp\in\Gam_{N}(X) :\Lrarw \alp<\pi \spand \alp\in\calg_{<N-1}(X) \spand \\
&& [\alp\in SR \rarw \fal\gam\in\cald_{\alp}(\gam\in\calg(X)\rarw\gam\in X)] \spand  \\
&& \fal\gam\in\calg_{<N-1}(X)[\gam\prec_{s[2;\gam,\alp]}\alp
  \Rarw \gam\leq X|\alp]  
\eeqnarrs
\edf

Let us examine the complexity of these operators. By induction on $i$ we see that 
$\calg_{i}$ is a $\Pi^{0}_{i}$-operator, and hence  
$\Gam_{N}$ is $\Pi^{0}_{N-1}$.
In this subsection we work in $\Pi^{0}_{N-1}\mbox{-Fix}$, and write $\Gam$ for $\Gam_{N}$, 
$|\alp|$ for $|\alp|_{\Gam_{N}}$, resp. 

We see easily that $\Gam=\Gam_{N}$ enjoys the hypotheses {\bf ($\Gam$.0)}, {\bf ($\Gam$.1)} and {\bf ($\Gam$.5)} 
in Subsection \ref{subsec:ensureoperators}. 
Furthermore {\bf ($\Gam$.3)} and {\bf ($\Gam$.4)} follow from the facts: 
if $\alp\not\in R^{\prime}$ or $\alp\in SR$, then $\alp\in \calg_{<N-1}(X)$ for any $X$.

\blem\label{lem:A1.2piN}
All of $\calg_{i}$ and $\Gam=\Gam_{N}$ are persistent and enjoys $\calw\incl\Gam(\calw)$.
\elem
\bprf
By Lemma \ref{lem:CXG} $\calg_{1}=\calg$ is persistent. 
Hence by induction on $i$ we see, from Lemma \ref{lem:3.2}.\ref{lem:3.2.2}, 
all of $\calg_{i}$ are persitent, too. 
Therefore so is $\Gam$. $\calw\incl\Gam(\calw)$ follows from Lemma \ref{lem:A1.2}.
\eprf

\subsection{Adequacy of  the operator $\Gam$}\label{subsec:adequate}

We next show that $\Gam$ enjoys the hypothesis {\bf ($\Gam$.2)}.
First we show the following lemmata.

\blem\label{lem:adq0}
Assume $\alp,\del\in\calg_{1}(X)$. Then $\alp<\del \spand \alp\not\prec\del \Rarw \alp\leq X|\del$. 
\elem
\bprf
Let $\eta$ denote the diagram such that $\alp\preceq\eta \spand \eta<\del<pd_{2}(\eta)$.
We have $\alp\leq\eta\in\calc^{\del}(X)|\del\incl X|\del$ by Lemma \ref{lem:CX4}.\ref{lem:CX4.2} 
and $\del\in\calg_{1}(X)$.
\eprf

\blem\label{lem:adq1+}
Let $2\leq i< N-1$ and $X=\Gam^{x}$ for an ordinal $x$. 
Then

\[
\alp\prec\bet \spand \alp\in\calg_{<i}(X),\bet\in\calg_{<i+1}(X) \Rarw
 \alp\leq X|\bet \lor \alp\prec_{s_{0}}\bet,
\]
where $s_{0}$  denotes  the longest initial segment of $s[2;\alp,\bet]$ such that $\# s_{0}\leq i-2$.
\elem
{\bf Proof}  by induction on $i$.
Suppose $\alp\not\prec_{s_{0}}\bet$.
By Lemma \ref{lem:B}
pick a {\it minimal\/} $\gam$ and an $s$ such that 
$s\incl_{e}s_{0}$, and
$$\alp\preceq_{s}\gam\lhd^{+}_{s}\bet \spand 2+\# s\leq i.$$

We show 
$$\gam\in\calg_{<2+\# s}(X).$$
Then by  $\bet\in\calg_{<i+1}(X)\incl\calg_{s}(X)$ we conclude $\alp\leq\gam\leq X|\bet$.

Assume $\alp\not\leq X|\bet$.
By Lemma \ref{lem:CX6+} we have for any $\del$
\[
\alp\preceq\del \Rarw \del\in\calg_{1}(X).
\]
Thus $\gam\in\calg_{1}(X)$. 

Let $t\in I$ denote an initial sequence such that $\# t<\#  s$, and suppose
$\calg_{<2+\#  t}(X)\ni\del\lhd^{+}_{t}\gam$.
Then $t\incl_{e}s[2;\del,\gam]$ by the definition, and $\# t<\#  s\leq\#  s_{0}\leq i-2$.
We  have to show $\del\leq X|\gam$.

Put
\[
\alp_{0}=\min\{\alp,\del\}, \: \alp_{1}=\max\{\alp,\del\} \mbox{ and } v=s[2;\alp_{0},\alp_{1}].
\]
It suffices to show $\alp_{0}\leq X|\gam\incl X|\bet$.

First consider the case when $\alp_{0}\not\preceq\alp_{1}$.
By Lemma \ref{lem:adq0} we have $\alp_{0}\leq X|\alp_{1}\incl X|\gam$.

In what follows assume $\alp_{0}\preceq\alp_{1}$.

Second consider the case when $s<_{lex}t$.
Then by Lemma \ref{clm:lemB}.\ref{clm:lemB3} we have $\del\lhd^{+}_{t}\bet$.
$\bet\in\calg_{<i+1}(X)\incl\calg_{2+\# t}(X)$ yields $\del\leq X|\bet$.
By the  assumption $\alp\not\leq X|\bet$ this means that $\del\leq X|\alp\incl X|\gam$.

In what follows suppose $s\not<_{lex}t$.
Then we have either $t\subset_{e}s$ or $t<_{lex}s$ by $\# t<\#  s$.
Let $w=s\cap t$ and $i=\ell(w)$.
Then $s|i=t|i$ and $t<_{lex}s \Rarw t(i)=0<1=s(i)$.

We show the following claim.
Claim \ref{clm:adq1+}.\ref{clm:adq1+3} yields the lemma by $X|\alp_{1}\incl X|\gam$.

\bclm\label{clm:adq1+}
\benu
\item\label{clm:adq1+0}
$t<_{lex}s \Rarw \del(w)\not\preceq_{i}\alp(w)$.
Therefore $\alp\neq\del$.
\item\label{clm:adq1+1}
$w\incl_{e}v$.
\item\label{clm:adq1+2}
$\alp_{0}\not\prec_{w}\alp_{1}$.
\item\label{clm:adq1+3}
$\alp_{0}\leq X|\alp_{1}$.
\eenu
\eclm
{\bf  Proof} of Claim \ref{clm:adq1+}.
\\
\ref{clm:adq1+}.\ref{clm:adq1+0}.
Assume  $t<_{lex}s$ and $\del(w)\preceq_{i}\alp(w)$.
Then by $t|(i+1)\incl_{e}t\incl_{e}s[2;\del,\gam]$ and $t(i)=0$,
$\not\exi \eta\{\del(w)\preceq_{i}\eta\prec_{i}\gam(w)\}$ with $w=t|i=s|i$.
On the other hand we have $\exi \eta\{\alp(w)\preceq_{i}\eta\prec_{i}\gam(w)\}$
by $\alp\prec_{s}\gam$ and $s(i)=1$.
A contradiction.

Assume $\alp=\del$. Then $\del(w)=\alp(w)$ and hence $t\not<_{lex}s$, i.e., $t\subset_{e}s$.
But then we have $\alp\prec_{t}\gam \spand \alp\lhd_{t}\gam$. This is not the case.
\\
\ref{clm:adq1+}.\ref{clm:adq1+1}.
Suppose $w\not\incl_{e}v$.
We have $w\incl_{e}s[2;\alp_{0},\gam]$ and $\alp_{1}\prec_{w|(\ell(w)-1)}\gam$.
Thus by Lemma \ref{lem:astup}.\ref{lem:astup1.4} we have $v<_{lex}w$.
But then $\alp(w)=\del(w)$, and hence $t\not<_{lex}s$ by Claim \ref{clm:adq1+}.\ref{clm:adq1+0}.
Hence  we have $w=t=s|i\subset_{e}s$, and this would yield 
$\alp(s|i)\lhd_{i}\gam(s|i)$ contradicting
$\alp\preceq_{s}\gam$.
\\
\ref{clm:adq1+}.\ref{clm:adq1+2}.
Suppose $\alp_{0}\prec_{w}\alp_{1}$.

First consider the case when $\alp_{1}=\alp$.
Then $\del\prec_{w}\alp\prec_{w}\gam$ and hence $\del\prec_{w}\gam$.
If $t\subset_{e}s$, then this means that $\del\prec_{t}\gam$, contradicting $\del\lhd_{t}\gam$.
$t\not<_{lex}s$ is not the case by Claim \ref{clm:adq1+}.\ref{clm:adq1+0}.

Next consider the case when $\alp_{1}=\del$.
Then $\alp\prec_{w}\del \spand \alp\prec_{w}\gam$ and hence $\del\prec_{w}\gam$.
$t\subset_{e}s$ is not the case since this contradicts $\del\lhd_{t}\gam$.
Assume $t<_{lex}s$.
Then $\alp(u)\prec_{j}\gam(u)=\del(u)$ for any $u$ with $w\subset_{e}u\incl_{e}s$, and
$j=\ell(u)$.
Therefore $\alp\prec_{s}\del\lhd_{s}^{+}\bet$ by Lemma \ref{clm:lemB}.\ref{clm:lemB3}.
This contradicts the minimality of $\gam$.
\\
\ref{clm:adq1+}.\ref{clm:adq1+3}.
Let $j-2:=\# w\leq\# t<\#s\leq i-2$.

First consider the case when $\alp_{1}=\alp$.
Then we have $\alp_{1}=\alp\in\calg_{<i}(X)\incl\calg_{<j+1}(X)$
and $\alp_{0}=\del\in\calg_{<2+\# t}(X)\incl\calg_{<j}(X)$.
Thus IH with Claims \ref{clm:adq1+}.\ref{clm:adq1+1} and \ref{clm:adq1+}.\ref{clm:adq1+2}
yields $\alp_{0}\leq X|\alp_{1}$.

Assume $\alp_{1}=\del$. 
We have $\alp_{0}=\alp\in\calg_{<i}(X)\incl\calg_{<j}(X)$.
On the other hand we have $\del\in\calg_{<2+\# t}(X)\incl\calg_{<j}(X)$.
Therefore we can assume $\# w=\# t$.

Again IH with Claims \ref{clm:adq1+}.\ref{clm:adq1+1} and \ref{clm:adq1+}.\ref{clm:adq1+2}
yields either $\alp\leq X|\del$ or $\alp\prec_{w_{0}}\del$
for the longest initial segment of $s[2;\alp,\del]$ such that $\# w_{0}\leq j-3=\# w-1$.
Now pick a $u\incl_{e}w$ and an $\eta$ 
so that $\alp\preceq_{u}\eta\lhd_{u}\del$ by Lemma \ref{lem:B}.
Then we have $\# u<\# t=\# w$ by Lemma \ref{lem:B'}.
Hence  $u\incl_{e}w_{0}$, and $\alp\not\prec_{w_{0}}\del$.
Thus  $\alp\leq  X|\del$.
\eprf
\smlskp
{\bf Proof} of Theorem \ref{th:AM} for $\Gam=\Gam_{N}$. 
Assume $\alp,\bet\in\calw$ and $\alp<\bet$. Put $x=|\alp|, y=|\bet|$. 
We show $x<y$ by induction on the natural sum $x\#y$. Suppose $x\geq y$. 
Put $X=\Gam^{x}, Y=\Gam^{y}$.  
We show $\alp\leq Y$. 
As in \cite{Wienpi3d} we see, using IH, 
$\alp\in\calg_{1}(X)|\bet=\calg_{1}(Y)|\bet\spand \alp\in \calg_{<N-1}(Y)$, 
and we can assume $\alp\prec\bet\in\cald^{Q}$ by IH.

Then Lemma \ref{lem:adq1+} with $i=N-2$ yields either $\alp\leq Y|\bet$
or $\alp\prec_{s}\bet$ for  $s=s[2;\alp,\bet]$
since $\#s\leq N-4=i-2$.
Assume $\alp\prec_{s}\bet$. Then $\calg_{<N-1}(Y)\ni\alp\prec_{s}\bet\in\Gam_{N}(Y)$.
Consequently $\alp\leq Y|\bet$ by Definition \ref{df:oprtN}.

This completes a proof of Theorem \ref{th:AM} for $\Gam=\Gam_{N}$.
\eprf

\subsection{Proof of Lemma \ref{lem:id5wf21}}\label{subsec:lem:id5wf21}
In this subsection we prove Lemma \ref{lem:id5wf21} for $Od(\Pi_{N})$ and $\Gam$.

\bdf\label{df:VWpiN}
$\calg_{i}:=\calg_{i}(\calw)$ {\rm and} $\calg_{<1+i}:=\calg_{<1+i}(\calw)$ {\rm for} $1\leq i<N-1$.
\edf

\blem\label{lem:GWQpiN}
For any $i\leq N-1$, if $\calg_{1}\ni\alp\preceq\rho\in\cald^{Q}$ and $\tau:=rg_{i}(\rho)\darw$, 
then $st_{i}(\rho)\in\calc^{\tau}(\calw)$.

In particular, $\calg_{1}\ni\alp\preceq\rho\in\cald^{Q} \Rarw st_{N-1}(\rho)\in\calc^{\pi}(\calw)=\calw_{\pi}$ 
for $\pi=rg_{N-1}(\rho)$.
\elem
\bprf
Assume $\calg_{1}\ni\alp\preceq\rho\in\cald^{Q}_{\sig}$ and put $\nu=st_{i}(\rho)$. 
Then $\alp\in\calc^{\alp}(\calw)$, and hence $\nu\in\calc^{\alp}(\calw)$. 
On the other hand we have $\fal\kap\leq\tau[K_{\kap}\nu<\alp]$ by 
the condition (\ref{cnd:Kst}), $(\cald.2)$ in Definition \ref{df:piN} and Lemma \ref{lem:3.2}.\ref{lem:3.2.5}. 
Lemma \ref{lem:CX3} with $\calc^{\alp}(\calw)|\alp\incl\calw$ yields $\nu\in\calc^{\tau}(\calw)$.
\eprf

\bdf\label{df:lexast}
\benu
\item \label{df:ast.1}
$\ovl{\alp}$ {\rm denotes the sequence of ordinal diagrams} $\{\alp(s):s\in I(2,N-2)\}$ 
{\rm ordered by the opposite relation of}
$<_{lex}$ {\rm on} $I${\rm :}
\[
\ovl{\alp}:=\langle\alp(s):s\in I(2,N-2) \rangle=\langle\cdots,\alp(s_{n}),\alp(s_{n+1}),\cdots\rangle
\mbox{ {\rm where }} s_{n+1}<_{lex}s_{n}.
\]
\item\label{df:lexN-1}
$\gam\prec_{N-1}^{pl}\alp$ 
{\rm denotes the lexicographic ordering on the finite sequences} $\ovl{\gam}$ 
{\rm of diagrams with respect to the ordering} $\prec_{N-2}^{p}${\rm :} $\gam\prec_{N-1}^{pl}\alp$ {\rm iff}
\[
\exi s\in I(2,N-2)[\fal t\in I(2,N-2)\{s<_{lex}t \Rarw \gam(t)=\alp(t)\} \spand \gam(s)\prec_{N-2}^{p}\alp(s)].
\]

\eenu
\edf

\blem\label{lem:GWpiN}
For each $\alp\in Od(\Pi_{N})|\pi$,
\[\alp\in\calg_{<N-1} \Rarw \alp\in\calw.\]
Specifically, for each $n\in\ome$, 
\[
\fal\alp\in\cald^{Q}\fal\alp_{\pi}[\alp\preceq\alp_{\pi}\in\cald_{\pi} \spand b(\alp_{\pi})<\ome_{n}(\pi+1) 
\spand \alp\in\calg_{<N-1} \Rarw \alp\in\calw].
\] 
\elem
\bprf
Assume $b(\alp_{\pi})<\ome_{n}(\pi+1)$ for $\alp\preceq\alp_{\pi}\in\cald_{\pi}$. 
Then Lemma \ref{lem:piNbarstb} yields $st_{N-1}(\alp)\leq st_{N-1}(\alp_{N-2}^{0})\leq b(\alp_{\pi})<\ome_{n}(\pi+1)$.

Assume $\alp\in\calg_{<N-1}$. Then for any $t\in I(2,N-2)$,
$st_{N-1}(\alp(t)), st_{N-1}(\alp(t)_{N-2}^{0})\in\calw_{\pi}|\ome_{n}(\pi+1)$
by Lemma \ref{lem:GWQpiN}.
It suffices to show $\alp\in\Gam(\calw)\incl\calw$, i.e., by Definition \ref{df:oprtN},
show that for any $\gam\in\calg_{<N-1}$ with 
$s=s[2;\gam,\alp]$, if $\gam\prec_{s}\alp$, then $\gam\leq\calw|\alp$.
By Lemma \ref{lem:barst} we have 
$\ovl{st}_{N-1}(\gam(s))<_{lex}\ovl{st}_{N-1}(\alp(s))$ for 
$\ovl{st}_{N-1}(\alp)=\langle st_{N-1}(\alp_{N-2}^{0}), st_{N-1}(\alp)\rangle$.

On the other hand we have
$\fal t\in I(2,N-2)\{s<_{lex}t \Rarw \gam(t)=\alp(t)\}$, i.e.,
$\gam\prec_{N-1}^{pl}\alp$ by Lemma \ref{lem:astup11}.

Thus the lemma is seen by induction along the lexicographic ordering $\prec_{N-1}^{pl}$.
\eprf
\smlskp
{\bf Proof} of Lemma \ref{lem:id5wf21} for $Od(\Pi_{N})$. We have to show for each $n\in\ome$
\[
\fal\alp\in\calw_{\pi}|\ome_{n}(\pi+1)\fal q\incl\calw_{\pi}|\ome_{n}(\pi+1) A(\alp,q).
\] 
By main induction on $\alp\in\calw_{\pi}|\ome_{n}(\pi+1)$ with subsidiary induction on $q\incl\calw_{\pi}|\ome_{n}(\pi+1)$.
Here observe that if $\bet_{1}\in\cald$ with $b(\bet_{1})<\ome_{n}(\pi+1)$, then by Lemma \ref{lem:N4aro} we have $Q(\bet_{1})\leq\max\{b(\bet_{1}),\pi\}<\ome_{n}(\pi+1)$.

Let $\alp_{1}\in\cald_{\sig}$ with $\sig\in\calw_{\pi}$ and $\alp=b(\alp_{1})\spand q=Q(\alp_{1})$.
By Theorem \ref{th:id5wf21} we have $\alp_{1}\in\calg_{1}$. We show $\alp_{1}\in\calw$. 
By Lemma \ref{lem:GWpiN} it suffices to show $\alp_{1}\in\calg_{<N-1}$.

We show the following claim.
Claim \ref{clm:5wf21.40id}.\ref{clm:5wf21.40.3id} with $K=1$ yields $\alp_{1}\in \calg_{<N-1}$.

\bclm\label{clm:5wf21.40id}
Let $\eta\in\calg_{<2+\# s}$ for an $s\in I(2)$.

Assume that there exist sequences $\{\eta_{k}\}_{k\leq K}$ of diagrams and $\{s_{k}\}_{k<K}\incl I(2)\, (K\geq 1)$
such that $\eta_{0}=\alp_{1}$, $\eta_{K}=\eta$, $s_{K-1}=s$,
$\fal k<K[\eta_{k+1}\lhd_{s_{k}}^{+}\eta_{k}]$
and $\fal k<K-1[s_{k}\subset_{e}s_{k+1} \spand \# s_{k}\leq\# s_{k+1}]$.

Then the followings hold.
\benu
\item\label{clm:5wf21.40.2id}
$\alp_{1}\prec rg_{\ell(s)}(\eta(s))\darw$ and $st_{\ell(s)}(\eta(s))\in\calw_{\pi}|\ome_{n}(\pi+1)$.
\item\label{clm:5wf21.40.3id}
$\eta\in\calw$.
\eenu
\eclm
{\bf Proof} of Claim \ref{clm:5wf21.40id}.\ref{clm:5wf21.40.2id}.
Put $\nu=st_{\ell(s)}(\eta(s))$ and $\tau=rg_{\ell(s)}(\eta(s))$. 

First we show $\alp_{1}\prec\tau$.
We have $\eta_{k+1}(s_{k})\lhd_{\ell(s_{k})}\eta_{k}(s_{k})$, and hence 
$rg_{\ell(s_{k})}(\eta_{k}(s_{k}))\preceq rg_{\ell(s_{k})}(\eta_{k+1}(s_{k}))$.
On the other hand we have 
$rg_{\ell(s_{k})}(\eta_{k+1}(s_{k}))\preceq rg_{\ell(s_{k+1})}(\eta_{k+1}(s_{k+1}))$
by $\ell(s_{k})\in In(\eta_{k+1}(s_{k})) \spand s_{k}\incl_{e}s_{k+1}$ and Lemma \ref{lem:astup}.\ref{lem:astup1.5}.
Hence we see 
$\alp_{1}=\eta_{0}\prec rg_{\ell(s_{0})}(\eta_{1}(s_{0}))\preceq rg_{\ell(s_{K-1})}(\eta_{K}(s_{K-1}))=\tau$.

Next we show $\nu\in\calw_{\pi}|\ome_{n}(\pi+1)$.
By Lemma \ref{lem:5.4}.\ref{lem:5.4.3} and $i<N-1$ we have $\nu<\pi$. 

Lemma \ref{lem:GWQpiN} with $\eta\in\calg_{<2+\# s}\incl\calg_{1}$ yields

\beqn\label{eq:5wf21.43id}
\nu\in\calc^{\tau}(\calw)
\eeqn

By Lemmata \ref{lem:5.4}.\ref{lem:5.4.3} and \ref{lem:Npi11exist} with $\alp_{1}\prec \tau$ we have

\beqn\label{eq:5wf21.41id}
\calb_{>\tau}(\nu)<b(\alp_{1})=\alp 
\eeqn

Now Lemma \ref{lem:id3wf20} together with $\mbox{MIH}(\alp)$, (\ref{eq:5wf21.43id}) and (\ref{eq:5wf21.41id}) yields 
$\nu\in\calc^{\pi}(\calw)=\calw_{\pi}$. 
 This shows Claim \ref{clm:5wf21.40id}.\ref{clm:5wf21.40.2id}.
\smlskp
{\bf Proof} of Claim \ref{clm:5wf21.40id}.\ref{clm:5wf21.40.3id}.

For each $s\in I(2)$ let
\[
E(s):=\sum\{3^{N-3-i}\cdot(1-s(i)): 2\leq i<\ell(s)\}+\sum\{3^{N-3-i}\cdot 2: \ell(s)\leq i\leq N-3\}.
\]
Observe that for $s, t\in I(2)$
\[
s\subset_{e}t \, \lor \, s<_{lex}t \Rarw E(s)>E(t).
\]
We show the Claim \ref{clm:5wf21.40id}.\ref{clm:5wf21.40.3id} by a triple induction:
by main induction on $E(s)$
with subsidiary induction on the length $\ell(rg_{\ell(s)}(\eta(s)))$ of the diagram $rg_{\ell(s)}(\eta(s))$
with sub-subsidiary induction on $st_{\ell(s)}(\eta(s))\in\calw_{\pi}|\ome_{n}(\pi+1)$.

We have $b(\eta_{\pi})=b((\alp_{1})_{\pi})=\alp<\ome_{n}(\pi+1)$.
By Lemma \ref{lem:GWpiN}
it suffices to show $\eta\in\calg_{t}$ for any $t\in I(2)$ with $\# t\geq\# s$.

Suppose $\calg_{<2+\# t}\ni\del\lhd_{t}^{+}\eta$.
We show $\del\in\calw$.

We can assume $t\not<_{lex}s$ by Lemma \ref{lem:C+}.\ref{lem:C+2}.
Then by $\# t\geq\#s$ we have one of the three cases $s\subset_{e}t$, $s<_{lex}t$ and $s=t$.

First consider the case when $s\subset_{e}t$.
Extend the sequences $\{\eta_{k}\}, \{s_{k}\}$ by one, i.e., $\del, t$.
MIH with $E(s)>E(t)$ yields $\del\in\calw$.

Second consider the case when $s_{K-1}=s<_{lex}t$. 
Since, in general, we have $u\subset_{e}s<_{lex}t \Rarw u\subset_{e}t \lor u<_{lex}t$,
let $K'=\max(\{k<K-1:s_{k}\subset_{e}t\}\cup\{0\})$.
Consider the sequences $\{\eta_{k}\}_{k\leq K'}\cup\{\del\}$ and $\{s_{k}\}_{k<K'}\cup\{t\}$.
Then we see $\del\lhd_{t}^{+}\eta_{K'}$ from Lemma \ref{clm:lemB}.\ref{clm:lemB3}.
Thus MIH with $E(s)>E(t)$ yields $\del\in\calw$.

Finally consider the case when $s=t$.
Then $\del\lhd_{s}^{+}\eta=\eta_{K}\lhd_{s}^{+}\eta_{K-1}$ yields $\del\lhd_{s}^{+}\eta_{K-1}$
by Lemma \ref{clm:lemB}.\ref{clm:lemB3}.
On the other hand we have $rg_{\ell(s)}(\eta(s))\preceq rg_{\ell(s)}(\del(s))$.
If $rg_{\ell(s)}(\eta(s))\prec rg_{\ell(s)}(\del(s))$, then 
$\ell(rg_{\ell(s)}(\eta(s)))>\ell(rg_{\ell(s)}(\del(s)))$.
Otherwise we have $st_{\ell(s)}(\eta(s))>st_{\ell(s)}(\del(s))$ 
by Lemma \ref{lem:5.4}.\ref{lem:5.4.5}.
Considering the sequences $\{\eta_{k}\}_{k<K}\cup\{\del\}$ and $\{s_{k}\}_{k<K}$,
SIH or SSIH yields $\del\in\calw$.

This shows Claim \ref{clm:5wf21.40id}.\ref{clm:5wf21.40.3id}, and
completes a proof of Lemma \ref{lem:id5wf21} for $\calw$.
\eprf

Lemma \ref{lem:id5wf21} yields Lemma \ref{th:id4wf22}: $\alp_{1}\in\calw_{\pi}$ for each $\alp_{1}\in Od(\Pi_{N})$ as in \cite{Wienpi3d}.

Consequently Lemma \ref{lem:WWOme} yields Theorem \ref{th:id5wfeachpiN}.

\section{Wellfoundedness proof by distinguished classes}\label{sec:5awf}
In this section we work in the set theory $\mbox{KP}\Pi_{N}$ for $\Pi_{N}$-reflecting universes and show the following theorem.

\bth\label{th:5wfeach}
For each $\alp\in Od(\Pi_{N})|\Ome$, $\mbox{\rm{KP}}\Pi_{N}$ proves that $(Od(\Pi_{N})|\alp,<)$ is a well ordering.
\eth

In \cite{hndodpiN} a system $(O(\Pi_{N}),<)$ was shown to be wellfounded. 
Our wellfoundedness proof in \cite{hndodpiN} is formalizable in the second order arithmetic 
$\Pi_{1}^{1}\mbox{-CA}_0+\Sig^{1-}_2\mbox{-CA}$, that is to say, we have assumed that the largest distinguished class 
$\calw_{D}$ defined by a $\Sig^{1-}_{2}$-formula exists as a set.
\smlskp
Our wellfoundedness proof is an extension of one for $O(\Pi_{3})$ in \cite{odpi3}, and is essntially the same as given in 
\cite{odpiN}. Lemmata stated without proofs are either easy or similarly seen as in \cite{odMahlo} and \cite{odpi3}.

 $X,Y,\ldots$ range over {\em subsets} of $Od(\Pi_{N})$. While $\calx,\caly,\ldots$ range over {\em classes}.

\subsection{Distinguished classes}\label{subsec:5awf.1}
In this subsection distinguished classes are defined and elementary facts on these classes are established.

Definitions concerning the distinguished class are modified by requiring that, for any distinguished class $\calx$,
$\alp\in\calx \Rarw \alp\in V^{*}(\calx)$ for a class $V^{*}(\calx)$ defined below.

\bdf\label{df:5wfuv}
\benu
\item\label{df:5wfuv.slhd}
{\rm For} $2\leq i<N-1$ {\rm define, cf. Lemma \ref{lem:5.4}.\ref{lem:5.4.5},}
\[
\alp\lhd^{s}_{i}\bet :\Lrarw 
\alp\prec_i\bet \spand \kap=rg_i(\alp)\darw=rg_i(\bet)\darw\spand i\in In(\alp).\]
\item $\alp^{-}:=\max\{\sig\in R: \sig\leq\alp\}$.
\item\label{df:5wfuv.Vsi}
{\rm For} $2\leq i\leq N-1$,
\[\bet\in U_{i}(X;\alp): \Lrarw \bet\in\cald^{Q} \Rarw
\bigcup\{K_{\sig}\nu:\nu=st_{i}(\bet), \sig\leq rg_{i}(\bet)\}\incl X|\alp^{-}.
\]
\item {\rm For} $2\leq i<N-1$,
$V^{s}_{i}(X;\del)$ {\rm denotes the} wellfounded part {\rm of the relation} $\{(\alp,\bet): \alp\in U_{i}(X;\del) \spand \alp\lhd^{s}_{i}\bet\}$.
{\rm Thus:}
\[
\alp\in V^{s}_{i}(X;\del) \Lrarw 
\fal\bet\in U_{i}(X;\del)[\bet\lhd^{s}_{i}\alp \rarw \bet\in V^{s}_{i}(X;\del)].
\]
\item {\rm For} $2\leq i\leq N-1$,
$U_{\geq i}(X;\del):=\bigcap_{i\leq j\leq N-1}U_{j}(X;\del)$.
\item {\rm For} $2\leq i<N-1$,
$H^{s}_{\geq i}(X;\del):=\bigcap_{i\leq j< N-1}H^{s}_{j}(X;\del)$.
\item {\rm For} $2\leq i<N-1$,
\[
\alp\in H^{s}_{i}(X;\del) :\Lrarw  
\fal\bet\in U_{i}(X;\del)[\bet\lhd^{s}_{i}\alp \rarw \bet\in V^{*}_{i+1}(X;\del)\cap V^{s}_{i}(X;\del)].
\]
\item {\rm For} $2\leq i\leq N-1$,
\[
\alp\in V^{*}_{i}(X;\del) :\Lrarw \fal\bet\{\alp\preceq_{i}\bet<\pi \Rarw \bet\in H^{s}_{\geq i}(X;\del)\}.
\]
{\rm Thus} $V^{*}_{N-1}(X;\del)=Od(\Pi_{N})$.
\label{df:5wfuv.V*}
\item {\rm For} $2\leq i\leq N-1$,
$\alp\in U^{*}_{i}(X;\del) :\Lrarw \fal\bet\{\alp\preceq_{i}\bet<\pi \Rarw \bet\in U_{\geq i}(X;\del)\}$.
\item {\rm For} $2\leq i\leq N-1$,
$UV^{*}_{i}(X;\del):=U^{*}_{i}(X;\del)\cap V^{*}_{i}(X;\del)$.
\item {\rm For} $2\leq i\leq N-1$, $\alp\in U_{i}(X) :\Lrarw \alp\in U_{i}(X;\alp)$. 
\\
$\alp\in V^{s}_{i}(X)$, $\alp\in H^{s}_{\geq i}(X)$, $\alp\in V^{*}_{i}(X)$, $\alp\in U^{*}_{i}(X)$ {\rm and} $\alp\in UV^{*}_{i}(X)$ 
{\rm are defined similarly by diagonalizations.}
\item
\[
\alp\in V^{*}(X) :\Lrarw \alp\in V^{*}_{2}(X) \spand \calc^{\alp}(X)|\alp\incl V^{*}_{2}(X).
\]
{\rm Cf. Lemmata \ref{lem:wf5.3.3-1} and \ref{lem:wf5.3}.\ref{lem:wf5.3.3} for the added condition} 
$\calc^{\alp}(X)|\alp\incl V^{*}_{2}(X)$.
\label{df:5wfuv.V*d}
\item $V\calc^{\alp}(X):=V^{*}(X) \cap \calc^{\alp}(X)$.
\label{df:3wfdtg.6}
\eenu
\edf

\blem\label{lem:5uv} 
\benu
\item\label{lem:5uv.0}
For $i<N-1$, $V^{*}_{i}(X;\del) \incl V^{*}_{i+1}(X;\del)\cap H^{s}_{\geq i}(X;\del)$ 
and $H^{s}_{i}(X;\del)\incl V^{s}_{i}(X;\del)$.

$T(X;\del)\ni\alp\prec_{i}\bet \Rarw \bet\in T(X;\del)$ for $i\leq N-1$ and $T\in\{V^{*}_{i},UV^{*}_{i}\}$.
\item \label{lem:5uv.1}
$\bigcup\{K_{\sig}\nu:\nu=st_{i}(\del), i\leq N-1, \sig\leq rg_{i}(\del)\}<\gam$ for $\cald^{Q}\ni\gam\preceq\del$.
\item \label{lem:5uv.3-1}
Assume $\caly\incl\calx \spand \alp\leq\bet$. Then $T(\caly;\alp)\incl T(\calx;\bet)$ for $T\in\{U_{i},U^{*}_{i}\}$ and  
$T(\caly;\alp)\supseteq T(\calx;\bet)$ for $T\in\{V^{s}_{i},H^{s}_{i},V^{*}_{i}\}$, and hence
$V^{*}_{2}(\caly)\supseteq V^{*}_{2}(\calx)$,
$V^{*}(\caly)\supseteq V^{*}(\calx)$.
\item \label{lem:5uv.2}
For any classes $\calx,\caly\incl Od(\Pi_{N})$ enjoying the condition {\bf (A)} in Definition \ref{df:CXcond}, i.e., $\fal\alp\in \calx[\alp\in\calc^{\alp}(\calx)]$ the following holds:
\[
\calx|\alp=\caly|\alp \Rarw \fal\bet<\alp^{+}\{V\calc^{\bet}(\calx)|\bet^{+}=V\calc^{\bet}(\caly)|\bet^{+}\}.
\]
\eenu
\elem
\bprf
\\
\ref{lem:5uv}.\ref{lem:5uv.0}. $V^{*}_{i}(X;\del) \incl V^{*}_{i+1}(X;\del)$ follows from Lemma \ref{lem:5.3}.\ref{lem:5.3.2}.
\\
\ref{lem:5uv}.\ref{lem:5uv.1}. 
Put $Y:=\bigcup\{K_{\sig}\nu:\nu=st_{i}(\del), i\leq N-1, \sig\leq rg_{i}(\del)\}$.
Then $Y<\del=\del^{-}$ follows from the condition (\ref{cnd:Kst}) 
of $(\cald.2)$ in Definition \ref{df:piN}.
Since $\ell Y<\ell \del$, we have $K_{\del}Y=\emptyset$, 
and hence $Y<\gam$ follows from Lemma \ref{lem:3.2}.\ref{lem:3.2.5}.
\\
\ref{lem:5uv}.\ref{lem:5uv.3-1}.
$V^{s}_{i}(\caly;\alp)\supseteq V^{s}_{i}(\calx;\bet)$ is seen from $U_{i}(\caly;\alp)\incl U_{i}(\calx;\bet)$.
Using this and $V^{*}_{N-1}(\calx;\del)=Od(\Pi_{N})$ for any $\calx$ and $\del$, we see $T(\caly;\alp)\supseteq T(\calx;\bet)$ for $T\in\{H^{s}_{i},V^{*}_{i}\}$ by induction on $N-i$.
\\
\ref{lem:5uv}.\ref{lem:5uv.2}. 
Let $\bet<\alp^{+}$, and $\calx, \caly$ enjoy the condition {\bf (A)} with $\calx|\alp=\caly|\alp$.
We have $\bet^{-}\leq\alp$, and hence $U_{i}(\calx;\bet)=U_{i}(\caly;\bet)$ for any $i$.
Therefore $H^{s}_{i}(\calx;\bet)=H^{s}_{i}(\caly;\bet)$. 
This means that $V^{*}_{2}(\calx)|\bet^{+}=V^{*}_{2}(\caly)|\bet^{+}$. 
Finally consider $V\calc^{\alp}(X)$.
By Lemmata \ref{lem:CX2}.\ref{lem:CX2.4} and \ref{lem:CX1} we have $\calc^{\bet}(\calx)=\calc^{\alp}(\calx)=\calc^{\alp}(\caly)=\calc^{\bet}(\caly)$.
Hence $V^{*}(\calx)|\bet^{+}=V^{*}(\caly)|\bet^{+}$ and $V\calc^{\bet}(\calx)|\bet^{+}=V\calc^{\bet}(\caly)|\bet^{+}$ for any $\bet<\alp^{+}$.
\eprf

\bremark\label{rem:alp}
{\rm Lemma \ref{lem:5uv}.\ref{lem:5uv.2} is needed for us to ensure the one of the most basic properties of distinguished classes} $\calx, \caly${\rm :}
\[
\alp\leq\calx \spand \alp\leq\caly \Rarw \calx|\alp^{+}=\caly|\alp^{+}.
\]

{\rm Here is the reason why we have restricted the sets} $U_{i}(X;\del)$ {\rm and} $V^{s}_{i}(X;\del)$ {\rm to} $\del^{-}$.
{\rm If the restriction is absent, say} $V^{s}_{i}(X):=V_{i}^{s}(X;\pi)${\rm , we would have}
$V^{*}_{2}(\calx)|\alp^{+}\neq V^{*}_{2}(\caly)|\alp^{+}$ {\rm even if} $\calx|\alp=\caly|\alp$
{\rm since for some} $\gam>\alp$, $\alp\prec\bet\spand \gam\prec\bet$ {\rm and hence we could have, e.g.,} 
$\gam\in U_{i}(\calx)\spand \gam\not\in U_{i}(\caly)$.

{\rm On the other side the upward requirement} $\bet\in H^{s}_{\geq i}(X;\del)$ {\rm for any} $\bet$ 
{\rm with} $\alp\preceq_{i}\bet<\pi$ {\rm in Definition \ref{df:5wfuv}.\ref{df:5wfuv.V*} of} $V^{*}_{i}(X;\del)$ {\rm will be used in Claim \ref{clm:5etaMh.b0} and Theorem \ref{th:pi11exist}.}
\eremark

\bdf \label{df:3wfdtg}{\rm For} $X,Y\subseteq Od(\Pi_{N})$ {\rm and} $\alp\in Od(\Pi_{N})$,
\benu
\item $D[X] :\Lrarw X<\pi \spand
\fal\alp(\alp\leq X\rarw WV\calc^{\alp}(X)|\alp^+= X|\alp^+)$.
\\
{\rm A class} $\calx$ {\rm is said to be a} distinguished class {\rm if} $D[\calx]${\rm . A} distinguished set {\rm is a set which is a distinguished class.}
\label{df:3wfdtg.8}
\item $\calw_{D}:=\bigcup\{X|\pi:D[X]\}$.
\label{df:3wfdtg.9}
\eenu
\edf

Observe that $V\calc^{\alp}(X)$ is $\Pi^{1}_{1}$ in $X$, $WV\calc^{\alp}(X)$ is $\Pi^{1}_{1}$ in $\Pi^{1}_{1}(X)$ and hence $D[X]$ is $\Del^{1}_{2}$ in $X$. Thus $\calw_{D}$ is $\Sig^1_2$ and hence a proper class. $\calw_{D}$ would exist as a set if we assume $\Sig^{1-}_{2}\mbox{-CA}$.

Obviously any distinguished class $\calx$ enjoys the condition {\bf (A)} $\fal\alp\in \calx[\alp\in\calc^{\alp}(\calx)]$ in Definition \ref{df:CXcond}.

\blem\label{lem:3.11.6}
$D[X] \spand \alp\in X \Rarw \fal\bet[\alp\in V^{*}(X)\cap\calc^{\bet}(X)]$.
\elem
\bprf
Assume $D[X] \spand \alp\in X$. Then $\alp\in X|\alp^{+}=WV\calc^{\alp}(X)\incl V^{*}(X)\cap\calc^{\alp}(X)$.
Hence $\alp\in\calc^{\bet}(X)$ for any $\bet\leq\alp$ by Lemma \ref{lem:CX2}.\ref{lem:CX2.3}. Moreover for $\bet>\alp$ we have $\alp\in X|\bet\incl\calc^{\bet}(X)$.
\eprf

In the following Lemmata \ref{lem:3wf6-1} and \ref{lem:3wf6}, ${\cal I}$ denotes a class family of distinguished classes, i.e., $\fal X\in{\cal I}\, D[X]$. Set
\[\calw_{{\cal I}}=\bigcup\{X\in{\cal I}:D[X]\}.\]

The following Lemmata \ref{lem:3wf6-1}, \ref{lem:3wf6}, \ref{lem:3acalwti} and \ref{lem:3wf9} are seen as in \cite{odMahlo} from Lemma \ref{lem:5uv}.\ref{lem:5uv.2} and \ref{lem:3.11.6}.

\blem\label{lem:3wf6-1}
$D[X] \spand X\in{\cal I} \spand \alp\leq X \Rarw \calw_{{\cal I}}|\alp^{+}=X|\alp^{+}$.
\elem

\blem\label{lem:3wf6}
$D[\calw_{{\cal I}}]$. In paricular $D[\calw_{D}]$.
\elem

\blem\label{lem:3acalwti}
$TI[\calw_{D}]$.
\elem

\blem\label{lem:3wf9}Suppose $D[X]$, 
$
\alp\in\calg(X)\cap V^{*}(X)$
and
\beqn\label{eq:3wf9hyp.2} 
\fal\bet(X<\bet \spand \bet^+<\alp^+ \Rarw 
WV\calc^{\bet}(X)|\bet^+\incl X)
\eeqn
 Then
 $\alp\in WV\calc^{\alp}(X)|\alp^+ \spand D[WV\calc^{\alp}(X)|\alp^+]$.
\elem

\blem\label{lem:wf5.3.3-1}
Let $X$ be a distinguished set. Assume $\gam\in V^{*}(X)$, $X|\gam\incl\calg(X)$ and 
$\alp<\gam \spand \fal\sig\leq\gam[K_{\sig}\alp\incl X]$.
Then $\calc^{\alp}(X)|\alp\incl V^{*}_{2}(X)$.
\elem
\bprf
We show by induction on $\ell\bet$ that
\[
\bet\in\calc^{\alp}(X)|\alp \Rarw \bet\in\calc^{\gam}(X)|\gam.
\]
On the other hand we have $\calc^{\gam}(X)|\gam\incl V^{*}_{2}(X)$ by $\gam\in V^{*}(X)$. 
Hence $\calc^{\alp}(X)|\alp\incl V^{*}_{2}(X)$ follows.

If $\bet\in X$, then $\bet\in\calc^{\gam}(X)$ follows from Lemma \ref{lem:3.11.6}.
If $\bet\not\in\cald$, then $\bet\in\calc^{\gam}(X)$ is seen from IH.
Assume $\bet\in\cald_{\sig}$ with a $\sig>\alp$.

First consider the case $\gam<\sig$. 
Then $\fal\kap\leq\gam[K_{\kap}\bet<\bet<\alp]$ by Lemma \ref{lem:3.2}.\ref{lem:3.2.10}. 
Lemma \ref{lem:CX3}.\ref{lem:CX3.1} with IH yields $\bet\in\calc^{\gam}(X)|\gam$. 

Finally assume $\alp<\sig\leq\gam$. Then pick a $\del\in K_{\sig}\alp\incl X|(\alp+1)$ such that 
$\bet\leq\del\in X|\gam$ by Lemmata \ref{lem:3.2}.\ref{lem:3.2.9} and \ref{lem:3.2}.\ref{lem:3.2.5}.
We claim $\bet\in X$. Then $\bet\in\calc^{\gam}(X)$ follows from Lemma \ref{lem:3.11.6} again.
Assume $\bet<\del$.
We have $\del\in\calg(X)$, and hence $\bet\in\calc^{\alp}(X)|\del\incl\calc^{\del}(X)|\del\incl X$ by Lemma \ref{lem:CX2}.\ref{lem:CX2.3}. 
We are done.
\eprf

\blem\label{lem:wf5.3}
Let $X$ be a distinguished set. Assume $\gam\in X$.
\benu
\item\label{lem:wf5.3.0}
$\fal\bet\in X|\gam\fal\sig[\alp\in\calc^{\bet}(X) \Rarw K_{\sig}\alp\incl\calc^{\bet}(X)]$.
\item\label{lem:wf5.3.1}
$\alp\in X|\gam \Rarw \fal\sig(K_{\sig}\alp\incl X)$.
\item\label{lem:wf5.3.2}
$\alp\in\calc^{\gam}(X) \Rarw \fal\sig\leq\gam(K_{\sig}\alp\incl X)$.
\item\label{lem:wf5.3.3}
$\alp\in\calc^{\gam}(X)|\gam \Rarw \alp\in X$.
\eenu
Therefore 
\[
D[X] \Rarw X\incl\calg(X).\]
\elem
\bprf
By main induction on $\gam\in X$ with subsidiary induction on $\ell\alp$ we show these simultaneously.
\\
\ref{lem:wf5.3}.\ref{lem:wf5.3.0}.
Let $\bet\in X|\gam$ and $\alp\in\calc^{\bet}(X)$. 
If $\alp\in X|\bet$, then $K_{\sig}\alp\incl X$ by MIH on Lemma \ref{lem:wf5.3}.\ref{lem:wf5.3.1}. 
Hence $K_{\sig}\alp\incl\calc^{\bet}(X)$ by Lemma \ref{lem:3.11.6}.
Otherwise the assertion follows from SIH.
For example if $\alp\in\cald_{\rho}$ with $\rho>\bet$, then $\{\rho\}\cup c(\alp)\incl\calc^{\bet}(X)$. 
SIH yields $K_{\sig}\alp\incl \{\alp\}\cup K_{\sig}(\{\rho\}\cup c(\alp))\incl\calc^{\bet}(X)$.
\\
\ref{lem:wf5.3}.\ref{lem:wf5.3.1}.
Assume $\alp\in X|\gam$. 
Then $\alp\in\calc^{\alp}(X)$, and hence $K_{\sig}\alp\incl\calc^{\alp}(X)$ by Lemma \ref{lem:wf5.3}.\ref{lem:wf5.3.0}.
We can assume $K_{\sig}\alp<\alp$ by Lemma \ref{lem:3.2}.\ref{lem:3.2.9}.
Then $K_{\sig}\alp\incl\calc^{\alp}(X)|\alp\incl X$ by MIH on Lemma \ref{lem:wf5.3}.\ref{lem:wf5.3.3}.
\\
\ref{lem:wf5.3}.\ref{lem:wf5.3.2}.
Assume $\alp\in\calc^{\gam}(X)$ and $\sig\leq\gam$.
If $\alp\in X|\gam$, then Lemma \ref{lem:wf5.3}.\ref{lem:wf5.3.1} yields $K_{\sig}\alp\incl X$.
Otherwise the assertion follows from SIH.
For example if $\alp\in\cald_{\rho}$ with $\rho>\gam\geq\sig$, then $\alp\not\prec\sig$ and 
$\{\rho\}\cup c(\alp)\incl\calc^{\gam}(X)$. 
SIH yields $K_{\sig}\alp\incl K_{\sig}(\{\rho\}\cup c(\alp))\incl X$.
\\
\ref{lem:wf5.3}.\ref{lem:wf5.3.3}.
Assume $\alp\in\calc^{\gam}(X)|\gam$. We show $\alp\in X$.
We have $WV\calc^{\gam}(X)|\gam^{+}=X|\gam^{+} \spand \gam\in V^{*}(X)$ by $\gam\in X$. 
Thus it suffices to show $\alp\in V^{*}(X)$. By $\gam\in V^{*}(X)$ we have $\alp\in V^{*}_{2}(X)$. 
On the other hand we have $\fal\sig\leq\gam[K_{\sig}\alp\incl X]$ by Lemma \ref{lem:wf5.3}.\ref{lem:wf5.3.2}, and
$X|\gam\incl\calg(X)$ by MIH on Lemma \ref{lem:wf5.3}.\ref{lem:wf5.3.3}.
Consequently Lemma \ref{lem:wf5.3.3-1} yields $\calc^{\alp}(X)|\alp\incl V^{*}_{2}(X)$, and hence $\alp\in V^{*}(X)$.
We are done.
\eprf

Thus $\tht[X]:\Lrarw D[X]$ enjoys these hypotheses {\bf ($\tht$.i)}($i\leq 2$) in Subsection \ref{subsec:proopr}
Therefore by Lemmata \ref{lem:id7wf8}.\ref{lem:id7wf8.3}, \ref{lem:id7wf8}.\ref{lem:id7wf8.3.5}, \ref{lem:KC} and \ref{cor:A1.3}.\ref{cor:A1.3.1} we have
the conditions {\bf (K)} $\fal\tau[\alp\in\calx \Rarw K_{\tau}\alp\incl\calx]$, 
{\bf (KC)} $\fal\alp\fal\bet\fal\sig[\alp\in\calc^{\bet}(\calx)\spand\sig\leq\bet \Rarw K_{\sig}\alp\incl \calx]$,
$\fal\bet\fal\tau[\alp\in\calc^{\bet}(X) \Rarw K_{\tau}\alp\incl\calc^{\bet}(X)]$ 
and $\calx\incl\calg(\calx)$ for any $\calx\in\{X: D[X]\}\cup\{\calw_{D}\}$.

\blem\label{lem:3wf16}
Let $X$ be a distinguished set, and suppose 
\beqn\label{eq:3wf16hyp.1}
\eta\in\calg(X)\cap V^{*}(X)
\eeqn
and 
\beqn\label{eq:3wf16hyp.2}
\fal\gam\prec\eta(\gam\in\calg(X)\cap V^{*}(X)  \rarw \gam\in X)
\eeqn
 Then
\[\eta\in WV\calc^{\eta}(X)|\eta^{+}\spand D[WV\calc^{\eta}(X)|\eta^{+}].\]
\elem
\bprf 
By Lemma \ref{lem:3wf9} and the hypothesis (\ref{eq:3wf16hyp.1}) it suffices to show
\beqn
\renewcommand{\theequation}{\ref{eq:3wf9hyp.2}} 
\fal\bet(X<\bet \spand \bet^+<\alp^+ \Rarw 
WV\calc^{\bet}(X)|\bet^+\incl X)
\eeqn
\addtocounter{equation}{-1}
Assume $X<\bet\spand \bet^+<\eta^+$. We have to show  $WV\calc^{\bet}(X)|\bet^+\incl X$. 
We prove this by induction on $\gam\in WV\calc^{\bet}(X)|\bet^+$. 
Suppose $\gam\in V\calc^{\bet}(X)|\bet^{+}$ and $V\calc^{\bet}(X)|\gam\incl X$. We show $\gam\in X$.

We show first 
\[\gam\in\calg(X).\]
 First $\gam\in V\calc^{\gam}(X)$ by $\gam\in V\calc^{\bet}(X)|\bet^{+}$. 
 Second we show the following claim by induction on $\ell\alp$:
\bclm\label{clm:3wf16}
$\alp\in\calc^{\gam}(X)|\gam \Rarw  \alp\in X$.
\eclm
{\bf Proof} of Claim \ref{clm:3wf16}. 
We have $\alp\in V^{*}_{2}(X)$ by $\gam\in V^{*}(X)$, and $\alp\in V^{*}(X)$ by Lemmata \ref{lem:wf5.3.3-1}, \ref{lem:wf5.3}.\ref{lem:wf5.3.3} and {\bf (KC)} for $X$. Therefore we can assume $\gam^{+}\leq\bet$ by $V\calc^{\bet}(X)|\gam\incl X$.

First consider the case $\alp\not\in\cald$. Then Lemma \ref{lem:CX3}.\ref{lem:CX3.1} with IH yields $\alp\in\calc^{\bet}(X)$.
Hence $\alp\in V\calc^{\bet}(X)|\gam\incl X$.

Therefore we can assume $\alp\in\cald_{\rho}$ for some $\rho>\gam$. Then $\{\rho\}\cup c(\alp)\incl\calc^{\gam}(X)$.
\\
{\bf Case1}. $\bet<\rho$: Then $\fal\kap\leq\bet[K_{\kap}(\{\rho\}\cup c(\alp))=K_{\kap}\alp<\alp<\gam]$.
Hence Lemma \ref{lem:CX3}.\ref{lem:CX3.1} with IH yields $\alp\in\calc^{\bet}(X)$
 and $\alp\in V\calc^{\bet}(X)|\gam\incl X$. 
\\
{\bf Case2}. $\bet\geq\rho$: We have $\cald_{\rho}\ni\alp<\gam<\rho\leq\bet$. 
Pick a $\del\in K_{\rho}\gam$ such that $\alp\leq\del\leq\gam$ by Lemma \ref{lem:3.2}.\ref{lem:3.2.5}.
$\gam\in \calc^{\bet}(X)$ with {\bf (KC)} yields $V\calc^{\del}(X)\ni\alp\leq\del\in X$. 
Therefore $\alp\in WV\calc^{\del}(X)|\del^{+}=X|\del^{+}$.

Thus Claim \ref{clm:3wf16} was shown.
\eprf

Hence we have $\gam\in\calg(X)\cap V^{*}(X)$.
We have $\gam<\bet^{+}\leq\eta \spand \gam\in\calc^{\gam}(X)$.
If $\gam\prec\eta$, then the hypothesis (\ref{eq:3wf16hyp.2}) yields $\gam\in X$.
In what follows assume $\gam\not\prec\eta$.

If $\fal\kap\leq\eta[K_{\kap}\gam<\gam]$, then Lemma \ref{lem:CX3}.\ref{lem:CX3.2} yields 
$\gam\in\calc^{\eta}(X)|\eta\incl X$ by $\eta\in\calg(X)$.

Suppose $\exi\kap\leq\eta[K_{\kap}\gam=\{\gam\}]$. This means that $\gam\in\cald$ and
$\exi\kap<\eta[\gam\prec\kap]$ by $\gam\not\prec\eta$. Let $\tau$ denote the maximal such one.
Then $\tau\in\calc^{\gam}(X)\spand \gam<\tau<\eta\spand \fal\kap\leq\eta[K_{\kap}\tau<\gam]$ by 
Lemmata \ref{lem:CX4}.\ref{lem:CX4.0}, \ref{lem:3.2}.\ref{lem:3.2.5} and \ref{lem:3.2}.\ref{lem:3.2.10}.
Lemma \ref{lem:CX3}.\ref{lem:CX3.2} yields 
$\tau\in\calc^{\eta}(X)|\eta\incl X$ by $\eta\in\calg(X)$.
Therefore $\tau\in X<\bet$.
$\gam\in\calc^{\bet}(X)$ with {\bf (KC)} yields $\{\gam\}=K_{\tau}\gam\incl X$.
We are done.
\eprf

Thus $\tht[X]:\Lrarw D[X]$ enjoys these hypotheses {\bf ($\tht$.i)}($i\leq 4$) in Subsection \ref{subsec:proopr} (demonstrably in the set theory KPM+V=L of recursively Mahlo universes with the axiom of constructibility).
Here note that we have $\alp\in V^{*}(X)$ for any $X$ and any $\alp\not\in \cald^{Q}$.

\subsection{Mahlo universes}\label{subsec:5awf.2}
In this subsection we introduce several classes of Mahlo universes and establish key facts on these classes. This is a crux in showing $Od(\Pi_{N})$ to be wellfounded without assuming the existence of the maximal distinguished class $\calw_{D}$.

From Theorem 2.4 in p.315 of \cite{Richter-Aczel74} we know that there exists a $\Pi_{3}$-sentence $ad$ such that $z$ is admissible iff $(z;\in)\models ad$. Put
\[lmtad :\Lrarw \fal x\exi y(x\in y\spand ad^{y})\]
Observe that $lmtad$ is a $\Pi_{2}$-sentence. Let $Lmtad$ denote the class of limits of admissible sets in a whole universe.

\bdf\label{df:3auni}
\benu
\item {\rm By a} universe {\rm we mean either a whole universe L with} $(\mbox{{\rm L}};\in)\models \mbox{{\rm KP}}\Pi_{N}$ {\rm or a transitive set} $Q\in \mbox{{\rm L}}$ {\rm in a whole universe L such that} $\ome\in Q${\rm . Universes are denoted} $P,Q,\ldots$
\item {\rm For a universe} $P$ {\rm and a set-theoretic sentence} $\vphi$, $P\models\vphi :\Lrarw (P;\in)\models\vphi$.
\item {\rm A universe} $P$ {\rm is said to be a} limit universe {\rm if} $lmtad^{P}$ {\rm holds, i.e.,} $P$ {\rm is a limit of admissible sets.}
\item {\rm For a universe} $P$, $\Del_{0}(\Del_{1})$ in $P$ {\rm denotes the class of predicates which are} $\Del_{0}$ {\rm in some} $\Del_{1}$ {\rm predicates on} $P${\rm .}
\item {\rm L denotes a whole universe L with} $(\mbox{{\rm L}};\in)\models \mbox{{\rm KP}}\Pi_{N}$.
\eenu
\edf

\blem\label{lem:3ahier} Let $P$ be a limit universe and $X\in {\cal P}(\ome)\cap P$. 
\benu
\item $V^{*}(X)$ and $WV\calc^{\alp}(X)$ are $\Del_{1}$ and $D[X]$ is $\Del_{0}(\Del_{1})$.
\label{lem:3ahier.1}
\item $V^{*}(X)=\{\alp: P\models\alp\in V^{*}(X)\}$, $WV\calc^{\alp}(X)=\{\alp:P\models \alp\in WV\calc^{\alp}(X)\}$ 
and $D[X]\Lrarw P\models D[X]$.
\label{lem:3ahier.2}
\eenu
\elem

\bdf\label{df:3awp}
{\rm For a limit universe} $P$ {\rm set}
\[\calw^{P}=\bigcup\{X\in P:D[X]\}=\bigcup\{X\in P:P\models D[X]\}.\]
\edf

Thus $\calw^{\mbox{\footnotesize L}}=\calw_{D}$ for the whole universe L.

\blem\label{lem:ausinP} For any limit universe $P$
\[D[\calw^{P}].\]
\elem
\bprf
This follows from Lemma \ref{lem:3wf6}.
\eprf

\blem\label{lem:5wuv}
For limit universes $P,Q$,
\[Q\in P\Rarw \calw^{Q}\incl\calw^{P}\spand \calw^{Q}\in P.\]
\elem

\blem\label{lem:3afin}
For any limit universe $P$
\[\bet\in \calc^{\alp}(\calw^{P}) \lrarw 
\exi X\in P\{D[X]\spand \bet\in \calc^{\alp}(X)\}.\]
\elem
\bprf
By the monotonicity of $\calc^{\alp}(X)$ we have the direction $[\larw]$.

The converse direction $[\rarw]$ is seen by induction on $\ell\bet$ using the fact 
\\
$\{X_{i}\}_{i<n}\incl P \spand \fal i<n D[X_{i}] \Rarw \bigcup_{i<n}X_{i}\in P \spand D[\bigcup_{i<n}X_{i}]$.
\eprf

The following lemma is seen as in Lemma \ref{lem:3afin} using Lemma \ref{lem:5uv}.\ref{lem:5uv.3-1}.

\blem\label{lem:3afinU}
For any limit universe $P$
\[\alp\in U^{*}_{i}(\calw^{P};\del) \rarw 
\exi X\in P\{D[X] \spand \alp\in U^{*}_{i}(X;\del)\}.\]
\elem

Some preparatory definitions are introduced. We say that a class $\calx\incl Lmtad$ is a $\Pi_{n}${\it -class} for $n\geq 2$ if there exists a set-theoretic $\Pi_{n}$-formula $F(\bar{a})$ with parameters $\bar{a}$ such that for any set $P$ with $\bar{a}\incl P$
\[P\in\calx\Lrarw (P;\in)\models F(\bar{a})\wedge lmtad.\]
Thus $P\in\calx$ is a $\Del_{0}$-formula. 
For a whole universe L, $\mbox{L}\in\calx$ denotes the formula $F(\bar{a})$.

By a $\Pi^{1}_{0}${\it -class} we mean a $\Pi_{n}$-class for some $n\geq 2$.

Referring \cite{Richter-Aczel74}, pp.322-327 let $\Pi_{i}(a)\, (i>0)$ denote a universal $\Pi_{i}$-formula uniformly on admissibles. Set
\[P\in M_{i}(\calx) :\Lrarw P\in Lmtad \spand \fal b\in P[P\models\Pi_{i}(b)\rarw \exi Q\in\calx\cap P(Q\models\Pi_{i}(b))].\]

Observe that $M_{i}(\calx)$ is a $\Pi_{i+1}$-class if $\calx$ is $\Pi^{1}_{0}$-class.

\blem\label{lem:4acalg} Let $\calx$ be a $\Pi^{1}_{0}$-class such that $\calx\incl Lmtad$.
Suppose $P\in M_{2}(\calx)$ and $\alp\in\calg(\calw^{P})$. Then there exists a universe $Q\in \calx$ such that $\alp\in\calg(\calw^{Q})$. 
\elem
\bprf Suppose $P\in M_{2}(\calx)$ and $\alp\in\calg(\calw^{P})$. 

First by $\alp\in\calc^{\alp}(\calw^{P})$ and Lemma \ref{lem:3afin} pick a distinguished set $X_{0}\in P$ so that 
$\alp\in \calc^{\alp}(X_{0})$

Next writing $\calc^{\alp}(\calw^{P})|\alp\incl\calw^{P}$ analytically we have
\[\fal\bet<\alp[\bet\in \calc^{\alp}(\calw^{P}) \Rarw \exi Y\in P(D[Y]\spand \bet\in Y)]\]
Again by Lemma \ref{lem:3afin} we have
\[\bet\in \calc^{\alp}(\calw^{P}) \lrarw \exi X\in P\{D[X]\spand \bet\in \calc^{\alp}(X)\}.\] 
Thus we have
\[
\fal\bet<\alp\fal X\in P\exi Y\in P[(D[X]\spand \bet\in \calc^{\alp}(X)) \Rarw (D[Y]\spand \bet\in Y)].
\]
By Lemma \ref{lem:3ahier}.\ref{lem:3ahier.2} we have $D[X]\lrarw P\models D[X]$ for any $X\in P$. Hence by Lemma \ref{lem:3ahier}.\ref{lem:3ahier.1} the following $\Pi_{2}$-predicate holds in the universe $P\in M_{2}(\calx)$:
\beqn\label{eq:3acalg}
\fal\bet<\alp\fal X\exi Y[(D[X]\spand \bet\in \calc^{\alp}(X)) \Rarw (D[Y]\spand \bet\in Y)]
\eeqn

Now pick a universe $Q\in P\cap\calx$ such that $X_{0}\in Q$, $Q\models (\ref{eq:3acalg})$.
Tracing the above argument backwards in the limit universe $Q$ we have $\calc^{\alp}(\calw^{Q})|\alp\incl\calw^{Q}$ and 
$X_{0}\incl\calw^{Q}=\bigcup\{X\in Q: Q\models D[X]\}\in P$. 
Thus by Lemma \ref{lem:3afin} we have $\alp\in \calc^{\alp}(\calw^{Q})$. Hence $\alp\in\calg(\calw^{Q})$.
\eprf
\smlskp
In the following key Definition \ref{df:5etaMh} we define $\Pi_{N}$-classes $M[\eta,i;\alp]\, (2\leq i\leq N-1)$ and $\Pi_{i+1}$-classes $M(\eta,i;\alp)\, (2\leq i<N-1)$ by induction on $N-i$ for $\eta\in\cald^{Q}$ and $\alp\in Od(\Pi_{N})|\pi$, 
cf. Subsections \ref{sec:prl5a.1}, \ref{sec:prl5a.2} and Remark \ref{rem:alp}.

\bdf\label{df:5etaMh} {\rm Let} $2\leq i\leq N-1$.
\benu
\item\label{df:5etaMh.1}
 $P\in M[\eta,N-1;\alp]$ {\rm iff} $P\in Lmtad$ {\rm and}
\[P\in\bigcap\{M_{N-1}(M[\gam,N-1;\alp]) : U^{*}_{N-1}(\calw^{P};\alp)\ni\gam\prec_{N-1}\eta\}.\]
\item\label{df:5etaMh.2}
 {\rm For} $2\leq i<N-1$, $P\in M(\eta,i;\alp)$ {\rm iff} $P\in Lmtad \spand i\in In(\eta)$ {\rm and}
\[P\in\bigcap\{M_{i}(M[\gam,i;\alp]) : U^{*}_{i}(\calw^{P};\alp)\ni\gam\lhd^{s}_{i}\eta\}.\]
\item\label{df:5etaMh.3}
 {\rm For} $2\leq i<N-1$, 
\[P\in M[\eta,i;\alp] :\Lrarw P\in M[\eta,i+1;\alp]\cap M_{i}(\calx(\eta_{i}^{0},i;\alp))\]
{\rm where} $\calx(\eta,i;\alp)$ {\rm denotes the class}
\[\calx(\eta,i;\alp):=M[\eta,i+1;\alp]\cap \bigcap\{M(\eta_{i}^{m},i;\alp) : m<lh_{i}(\eta)-1\}.\]
{\rm Note that}
\[\calx(\eta_{i}^{0},i;\alp)=M[\eta_{i}^{0},i+1;\alp]\cap \bigcap\{M(\eta_{i}^{m},i;\alp) : m<lh_{i}(\eta)-1\}.\]
\item\label{df:5etaMh.4}
$M[\eta,2]:=M[\eta,2;\eta]$.
\eenu
\edf

Let us examine how these classes are defined for a fixed $\alp$. Note that $U^{*}_{i}(\calw^{P};\alp)$ is $\Sig_{1}$ on $P$.
First $M[\eta,N-1;\alp]$ is defined as a $\Pi_{N}$-class using a recursion lemma, cf. \cite{Richter-Aczel74}, pp.322-327. Namely there exists a primitive recursive function $g^{[N-1]}(\eta;\alp)$ such that 
\[P\in M[\eta,N-1;\alp] \Lrarw P\models\Pi_{N}(g^{[N-1]}(\eta;\alp)).\]

Suppose that $\Pi_{N}$-classes $M[\gam,i+1;\alp]$ has been defined for any $\gam$ and let $g^{[i+1]}(\gam;\alp)$ be a primitive recursive function such that 
\[P\in M[\gam,i+1;\alp] \Lrarw P\models\Pi_{N}(g^{[i+1]}(\gam;\alp)).\]
Then $M(\eta,i;\alp)$ and $M[\eta,i;\alp]$ are defined as a $\Pi_{i+1}$-class and a $\Pi_{N}$-class, resp. using the function $g^{[i+1]}$ and a simultaneous recursion lemma.

We say that a universe $P$ is an $\eta${\it -Mahlo universe} if $P\in M[\eta,2]$.

Now we show that if $P$ is an $\eta$-Mahlo universe and $UV^{*}_{2}(\calw^{P})\ni\gam\prec\eta$, then $P$ is $\Pi_{2}$-reflecting on $\gam$-Mahlo universes, i.e., $P\in M_{2}(M[\gam,2])$, cf. Theorem \ref{th:5awf16}.

\blem\label{lem:5etaMh-2}If $\alp<\bet$, then $M[\eta,i;\bet]\incl M[\eta,i;\alp]$ for $2\leq i\leq N-1$.
\elem
\bprf
Assume $\alp<\bet$. Then $U^{*}_{i}(\calx;\alp)\incl U^{*}_{i}(\calx;\bet)$ for any classes $\calx$. By induction on $\in$ with subsidiary induction on $N-i$ show simultaneously 
\[P\in M[\eta,i;\bet] \Rarw P\in M[\eta,i;\alp]\, (2\leq i\leq N-1),\]
 and 
\[P\in M(\eta,i;\bet) \Rarw P\in M(\eta,i;\alp)\, (2\leq i<N-1).\]
\eprf

\blem\label{lem:5etaMh-1}
Assume $P\in M_{i+1}(M[\eta,i+1;\alp])$, and either $P\in M_{i}(\calx(\eta_{i}^{0},i;\alp))$ or $\eta=\eta_{i}^{0}\spand P\in \bigcap\{M(\eta_{i}^{m},i;\alp) : m<lh_{i}(\eta)-1\}$. Then $P\in M_{i}(M[\eta,i;\alp])$.
\elem
\bprf

First consider the case when $P\in M_{i+1}(M[\eta,i+1;\alp])\cap M_{i}(\calx(\eta_{i}^{0},i;\alp))$. Since $M_{i}(\calx(\eta_{i}^{0},i;\alp))$ is a $\Pi_{i+1}$-class, $P$ reflects it on $M[\eta,i+1;\alp]$. 
Namely we have with $M[\eta,i;\alp]=M[\eta,i+1;\alp]\cap M_{i}(\calx(\eta_{i}^{0},i;\alp))$,
$P\in M_{i+1}(M[\eta,i;\alp])\incl M_{i}(M[\eta,i;\alp])$.

Next consider the case when $P\in M_{i+1}(M[\eta,i+1;\alp])\cap\bigcap\{M(\eta_{i}^{m},i;\alp) : m<lh_{i}(\eta)-1\}$ and $\eta=\eta_{i}^{0}$. By the first case it suffices to show $P\in M_{i}(\calx(\eta_{i}^{0},i;\alp))$ with the class
\[\calx(\eta_{i}^{0},i;\alp)=M[\eta_{i}^{0},i+1;\alp]\cap \bigcap\{M(\eta_{i}^{m},i;\alp) : m<lh_{i}(\eta)-1\}.\]
Since $M(\eta_{i}^{m},i;\alp)$ are $\Pi_{i+1}$-classes, $P$ reflects these on $M[\eta_{i}^{0},i+1;\alp]=M[\eta,i+1;\alp]$. 
Namely we have $P\in M_{i+1}(\calx(\eta_{i}^{0},i;\alp))\incl M_{i}(\calx(\eta_{i}^{0},i;\alp))$. We are done.
\eprf

The following is the key lemma.

\blem\label{lem:5etaMh}
\benu
\item\label{lem:5etaMh.a}
Suppose $2\leq i\leq N-1$, $pd_{i}(\gam)=\eta$ and $\gam\in UV^{*}_{i}(\calw^{P};\alp)$. 
Then for any {\rm class} $P$
\[
P\in M[\eta,i;\alp] \Rarw P\in M_{i}(M[\gam,i;\alp]).
\] 
\item\label{lem:5etaMh.b}
 Suppose $N-1>i\in In(\gam) \spand \gam\in H^{s}_{i}(\calw^{P};\alp) \spand \eta=\gam_{i}^{1}$. 
 Then for any {\rm set} $P$
\[P\in \calx(\eta,i;\alp) \Rarw P\in M(\gam,i;\alp)\]
with $\calx(\eta,i;\alp)=M[\eta,i+1;\alp]\cap \bigcap\{M(\eta_{i}^{m},i;\alp) : m<lh_{i}(\eta)-1\}$.
\eenu
\elem
\bprf
These are shown simultaneously by induction on $N-i$. First we show Lemma \ref{lem:5etaMh}.\ref{lem:5etaMh.a} for the case $i=N-1$. Second Lemma \ref{lem:5etaMh}.\ref{lem:5etaMh.b} is shown assuming Lemma \ref{lem:5etaMh}.\ref{lem:5etaMh.a} for the case $i+1$. Finally Lemma \ref{lem:5etaMh}.\ref{lem:5etaMh.a} is proved assuming Lemma \ref{lem:5etaMh}.\ref{lem:5etaMh.b}.
\\
\ref{lem:5etaMh}.\ref{lem:5etaMh.a} for the case $i=N-1$ follows from the definition of $P\in M[\eta,N-1;\alp]$.
Namely we have
\[U^{*}_{N-1}(\calw^{P};\alp)\ni\gam\prec_{N-1}\eta \spand P\in M[\eta,N-1;\alp] \Rarw P\in M_{N-1}(M[\gam,N-1;\alp]).\]
\smallskip

In what follows assume $i<N-1$. First we show the following claim.
\bclm\label{clm:5etaMh.b0}
Assume Lemma \ref{lem:5etaMh}.\ref{lem:5etaMh.a} holds for an $i<N-1$. Then the hypothesis $pd_{i}(\gam)=\eta$ can be weakened to $\gam\prec_{i}\eta$ as follows:
Suppose $\gam\prec_{i}\eta$ and $\gam\in UV^{*}_{i}(\calw^{P};\alp)$. 
 Then
\[P\in M[\eta,i;\alp] \Rarw P\in M_{i}(M[\gam,i;\alp]).\] 
\eclm
{\bf Proof} of Claim \ref{clm:5etaMh.b0}.
This is seen from Lemmata \ref{lem:5uv}.\ref{lem:5uv.0}, \ref{lem:5uv}.\ref{lem:5uv.3-1} and \ref{lem:3afinU}, and the fact $M_{i}(M_{i}(\calx))\incl M_{i}(\calx)$.
\eprf
\\
\ref{lem:5etaMh}.\ref{lem:5etaMh.b}. Suppose $i\in In(\gam) \spand \gam\in H^{s}_{i}(\calw^{P};\alp)$.
 We show $P\in M(\gam,i;\alp)$ by subsidiary induction on $\gam\in V^{s}_{i}(\calw^{P};\alp)$ for {\it sets\/} $P$ assuming $P\in \calx(\eta,i;\alp)$ for the diagram $\eta=\gam_{i}^{1}=(rg_{i}(\gam))^{0}_{i}$. 
Note that the relation $\lhd^{s}_{i}$ is transitive.

 Assume $U^{*}_{i}(\calw^{P};\alp)\ni\del\lhd^{s}_{i}\gam$. Then $\del\in V^{s}_{i}(\calw^{P};\alp)\cap UV^{*}_{i+1}(\calw^{P};\alp)$ by $\gam\in H^{s}_{i}(\calw^{P};\alp)$, and $rg_{i}(\del)=rg_{i}(\gam)$ and hence $\del^{1}_{i}=\eta$. By Lemma \ref{lem:5ap12} and Definition \ref{df:5res} we have $\fal m[0<m<lh_{i}(\del)-1 \Rarw \del_{i}^{m}=\eta_{i}^{m-1}]$ and $\del_{i}^{0}=\del$ since $i\in In(\del)$. 

We have to show $P\in M_{i}(M[\del,i;\alp])$. By Lemma \ref{lem:5etaMh-1} and $\del_{i}^{0}=\del$ it suffices to show the following claim.

\bclm \label{clm:5etaMh.b}
\benu
\item\label{clm:5etaMh.b.1}
$P\in M_{i+1}(M[\del,i+1;\alp])$.
\item\label{clm:5etaMh.b.2}
$P\in\bigcap\{M(\del_{i}^{m},i;\alp) : m<lh_{i}(\del)-1\}$.
\eenu
\eclm
{\bf Proof} of Claim \ref{clm:5etaMh.b}.
\\
\ref{clm:5etaMh.b}.\ref{clm:5etaMh.b.1}. 
By $\del^{1}_{i}=\eta$ and Lemma \ref{lem:5Si-2} we have $\del=\del_{i}^{0}\prec_{i+1}\eta$. 
We have $P\in \calx(\eta,i;\alp)\incl M[\eta,i+1;\alp]$.
Therefore Claim \ref{clm:5etaMh.b0} and IH on Lemma \ref{lem:5etaMh}.\ref{lem:5etaMh.a} for the case $i+1$ yield Claim \ref{clm:5etaMh.b}.\ref{clm:5etaMh.b.1}.
\\
\ref{clm:5etaMh.b}.\ref{clm:5etaMh.b.2}.
First consider the case $m>0$. 
Then we have $\del_{i}^{m}=\eta_{i}^{m-1}$ and $P\in \calx(\eta,i;\alp)\incl M(\eta_{i}^{m-1},i;\alp)=M(\del_{i}^{m},i;\alp)$. 
On the other side $P\in M(\del,i;\alp)$ follows from SIH.
Thus Claim \ref{clm:5etaMh.b}.\ref{clm:5etaMh.b.2} was shown.
\eprf
\smlskp
\ref{lem:5etaMh}.\ref{lem:5etaMh.a} for the case $i<N-1$. Suppose $pd_{i}(\gam)=\eta$, $\gam\in UV^{*}_{i}(\calw^{P};\alp)$ and $P\in M[\eta,i;\alp]$. We have to show $P\in M_{i}(M[\gam,i;\alp])$.

By Lemma \ref{lem:3.23.1} one of the following cases occur:
\\
{\bf Case \ref{lem:3.23.1}.1} 
\[
\eta=pd_{i}(\gam)=pd_{i+1}(\gam) \spand lh_{i}(\gam)=lh_{i}(\eta)\spand \fal m<lh_{i}(\gam)[\gam_{i}^{m}=\eta_{i}^{m}]
\]
$\gam\in UV^{*}_{i+1}(\calw^{P};\alp)$ and $P\in M[\eta,i;\alp]\incl M[\eta,i+1;\alp]$. 
Thus by IH we have $P\in M_{i+1}(M[\gam,i+1;\alp])$, and hence by Lemma \ref{lem:5etaMh-1} it suffices to show $P\in M_{i}(\calx(\gam_{i}^{0},i;\alp))$. This follows from $\gam_{i}^{0}=\eta_{i}^{0}$ and $P\in M[\eta,i;\alp]\incl M_{i}(\calx(\eta_{i}^{0},i+1;\alp))$.

In what follows assume $i\in In(\gam)$.
\\
{\bf Case \ref{lem:3.23.1}.2} 
\[
rg_{i}(\gam)=pd_{i}(\gam)=\eta \spand \gam_{i}^{0}=\gam \spand \fal m<lh_{i}(\eta)=lh_{i}(\gam)-1[\eta_{i}^{m}=\gam_{i}^{1+m}]
\]
Then $\gam=\gam_{i}^{0}\prec_{i+1}\gam_{i}^{1}=\eta_{i}^{0}$. We have $P\in M[\eta,i;\alp]\incl M_{i}(\calx(\eta_{i}^{0},i;\alp))$ with 
\beqnarrs
\calx(\eta_{i}^{0},i;\alp) & = & \calx(\gam_{i}^{1},i;\alp) \\
& = & M[\gam_{i}^{1},i+1;\alp]\cap\bigcap\{M(\gam_{i}^{m},i;\alp): 0<m<lh_{i}(\gam)-1\}.
\eeqnarrs
Let $Q$ be any limit universe in $P$ such that $Q\in \calx(\gam_{i}^{1},i;\alp)$ and 
$\gam\in U^{*}_{i+1}(\calw^{Q};\alp)$, cf. Lemma \ref{lem:3afinU}. 
We claim that $Q\in M_{i}(M[\gam,i;\alp])$. This yields 
\[P\in M_{i}(\calx(\gam_{i}^{1},i;\alp))\incl M_{i}(M_{i}(M[\gam,i;\alp]))\incl M_{i}(M[\gam,i;\alp]),\]
 and hence we are done. By the definition we have
\beqnarr
Q & \in & M[\gam_{i}^{1},i+1;\alp] \label{eq:5etaMh.a-1} \\
Q & \in & \bigcap\{M(\gam_{i}^{m},i;\alp): 0<m<lh_{i}(\gam)-1\} \label{eq:5etaMh.a0}
\eeqnarr

By Lemma \ref{lem:5wuv} we have $\calw^{Q}\incl\calw^{P}$, and hence 
\beqn\label{eq:5etaMh.a1}
\gam\in UV^{*}_{i+1}(\calw^{Q};\alp)
\eeqn
 and 
\beqn\label{eq:5etaMh.a2}
\gam\in H^{s}_{i}(\calw^{Q};\alp)
\eeqn
On the other hand we have $\gam\prec_{i+1}\gam_{i}^{1}$. 
Therefore IH on Lemma \ref{lem:5etaMh}.\ref{lem:5etaMh.a}, i.e., Claim \ref{clm:5etaMh.b0} with 
(\ref{eq:5etaMh.a-1}) and (\ref{eq:5etaMh.a1}) yields 
\beqn\label{eq:5etaMh.a3}
Q\in M_{i+1}(M[\gam,i+1;\alp])
\eeqn
 On the other side Lemma \ref{lem:5etaMh}.\ref{lem:5etaMh.b} with the {\it set\/} $Q\in \calx(\gam_{i}^{1},i;\alp)$ and (\ref{eq:5etaMh.a2}) yields $Q\in M(\gam,i;\alp)$.
Hence we have by (\ref{eq:5etaMh.a0})
\beqn\label{eq:5etaMh.a5}
Q\in\bigcap\{M(\gam_{i}^{m},i;\alp) : m<lh_{i}(\gam)-1\}
\eeqn
Now Lemma \ref{lem:5etaMh-1} with $\gam_{i}^{0}=\gam$, (\ref{eq:5etaMh.a3}) and (\ref{eq:5etaMh.a5}) yields $Q\in M_{i}(M[\gam,i;\alp])$ as desired.
\\
{\bf Case \ref{lem:3.23.1}.3}
\beqnarrs
&& \eta=pd_{i}(\gam)\prec_{i}rg_{i}(\gam) \spand \gam_{i}^{0}=\gam \spand  \\
&& \exi m[0<m\leq lh_{i}(\eta)-1 \spand rg_{i}(\eta_{i}^{m-1})=rg_{i}(\gam)\spand st_{i}(\eta_{i}^{m-1})>st_{i}(\gam)\spand \\
&& \fal k<lh_{i}(\eta)-m+1=lh_{i}(\gam)(k>0 \rarw \eta_{i}^{m-1+k}=\gam_{i}^{k})]
\eeqnarrs
 Then we have $U^{*}_{i}(\calw^{P};\alp)\ni\gam\lhd^{s}_{i}\eta_{i}^{m-1}$ for an $m$ with $0<m\leq lh_{i}(\eta)-1$. $P\in M[\eta,i;\alp]\incl M(\eta_{i}^{m-1},i;\alp)$ yields $P\in M_{i}(M[\gam,i;\alp])$.
\eprf

Lemma \ref{lem:5etaMh}.\ref{lem:5etaMh.a}, i.e., Claim \ref{clm:5etaMh.b0} yields the following Theorem \ref{th:5awf16}.

\bth\label{th:5awf16} Let $P$ be an $\eta$-Mahlo universe. Then $P$ is $\Pi_{2}$-reflecting on $\gam$-Mahlo universes for 
$UV^{*}_{2}(\calw^{P})\ni\gam\prec\eta$: 
\[
P\in M[\eta,2]\spand UV^{*}_{2}(\calw^{P})\ni\gam\prec\eta\Rarw P\in M_{2}(M[\gam,2]).
\]
\eth
\bprf
By Lemma \ref{lem:5etaMh-2} we have $P\in M[\eta,2]=M[\eta,2;\eta]\incl M[\eta,2;\gam]$. On the other hand we have $\gam\in UV^{*}_{2}(\calw^{P};\gam)$. 
Thus Lemma \ref{lem:5etaMh}.\ref{lem:5etaMh.a}, i.e., Claim \ref{clm:5etaMh.b0} yields $P\in M_{2}(M[\gam,2;\gam])$.
\eprf

\blem\label{lem:gvinclu}
$\calg(X)\incl U^{*}_{2}(X)$.
\elem
\bprf
Assume $\gam\in\calg(X)$. Let $\del\in\cald^{Q}$ such that $\gam\preceq\del$ and $\nu=st_{i}(\del)$ for an $i\geq 2$. We have to show $Y:=\bigcup\{K_{\sig}\nu: \sig\leq rg_{i}(\del)\}\incl X|\gam$.
By Lemma \ref{lem:5uv}.\ref{lem:5uv.1} we have $Y<\gam$.

On the other hand we have $\gam\in\calc^{\gam}(X)$, and this yields $\del\in\calc^{\gam}(X)$, and hence $\nu\in\calc^{\gam}(X)$ by the definition of the set $\calc^{\gam}(X)$.
Therefore $Y\incl\calc^{\gam}(X)$ follows from Lemma \ref{lem:id7wf8}.\ref{lem:id7wf8.3.5}.
Thus we have $Y\incl\calc^{\gam}(X)|\gam\incl X$.
\eprf

\bth\label{th:5awf16.1} For any {\rm set} $P$ and $\eta$
\[\eta\in\calg(\calw^{P})\cap V^{*}(\calw^{P})\spand P\in M_{2}(M[\eta,2]) 
\Rarw \eta\in\calw^{P}\]
\eth
\bprf We show this by induction on $\in$. Suppose, as IH, the theorem holds for any limit universe $Q\in P$.
By Lemma \ref{lem:4acalg} pick a $Q\in P$ so that for $X=\calw^{Q}\in P$,
$\eta\in\calg(X) \spand Q\in M[\eta,2]$. 
By Lemma \ref{lem:5uv}.\ref{lem:5uv.3-1} we also have $\eta\in V^{*}(X)$.
 Hence we have
\beqn\renewcommand{\theequation}{\ref{eq:3wf16hyp.1}} 
\eta\in\calg(X)\cap V^{*}(X)
\eeqn
\addtocounter{equation}{-1}
On the other side Lemma \ref{lem:gvinclu} and Theorem \ref{th:5awf16} yield 
\[
\fal\gam\prec\eta\{\gam\in\calg(X)\cap V^{*}(X) \Rarw  Q\in M_{2}(M[\gam,2])\}.
\]

Now IH yields
\beqn\renewcommand{\theequation}{\ref{eq:3wf16hyp.2}}
\fal\gam\prec\eta(\gam\in\calg(X)\cap V^{*}(X)  \rarw \gam\in X)
\eeqn
\addtocounter{equation}{-1}
Therefore by Lemmata \ref{lem:ausinP} and \ref{lem:3wf16} we conclude 
\[\eta\in WV\calc^{\eta}(X)|\eta^{+}\in P\spand D[WV\calc^{\eta}(X)|\eta^{+}]\]
and hence $\eta\in\calw^{P}$.
\eprf

Now we establish the existence of $\eta$-Mahlo universes for each $\eta$. Here the $\Pi_{N}$-reflection of the whole universe enters proofs.

Recall that $\calw_{\pi}=\calc^{\pi}(\calw^{\mbox{\footnotesize L}})$.

\blem\label{lem:4aro}
Let $\mbox{{\rm L}}$ be a whole universe such that $\mbox{{\rm L}}\models \mbox{{\rm KP}}\Pi_{N}$, and $\eta\in\cald^{Q}$. 
For each $n\in\ome$, if $b(\eta)<\ome_{n}(\pi+1)$, then
\[\fal\gam\in\cald^{Q}\{\eta\preceq\gam\in U^{*}_{N-1}(\calw^{\mbox{\footnotesize L}};\eta) 
\Rarw \mbox{{\rm L}}\in M[\gam,N-1;\eta]\}.\]
\elem
\bprf 
Let $\eta\preceq\gam\in\cald^{Q}$. 

First we show
\beqn\label{eq:4aro}
\eta\preceq\gam\in U^{*}_{N-1}(\calw^{\mbox{\footnotesize L}};\eta)  \Rarw st_{N-1}(\gam)\in\calw_{\pi}|\ome_{n}(\pi+1)
\eeqn

By Lemmata \ref{lem:N4aro} and \ref{lem:Npi11exist} we have 
$st_{N-1}(\gam)\leq Q(\gam)\leq\max\{b(\gam),\pi\}\leq\max\{b(\eta),\pi\}<\ome_{n}(\pi+1)$.

On the other side $\gam\in U^{*}_{N-1}(\calw^{\mbox{\footnotesize L}};\eta)$ means that $\fal\bet[\gam\preceq_{N-1}\bet \Rarw \bet\in U_{N-1}(\calw^{\mbox{\footnotesize L}};\eta)]$.
In particular $\bigcup\{K_{\sig}\nu: \sig\leq\pi\}\incl\calw^{\mbox{\footnotesize L}}$ with $\nu=st_{N-1}(\gam)$. 
Therefore Lemma \ref{lem:CX4aro} with the condition $(\cald.2)$ (\ref{cnd:Kst}) in Definition \ref{df:piN}
 yields $\nu\in\calc^{\pi}(\calw^{\mbox{\footnotesize L}})=\calw_{\pi}$.
Thus we have shown (\ref{eq:4aro}).

We show that 
$\fal\gam\in U^{*}_{N-1}(\calw^{\mbox{\footnotesize L}};\eta)\{\eta\preceq\gam \Rarw \mbox{{\rm L}}\in M[\gam,N-1;\eta]\}$ 
by induction on $st_{N-1}(\gam)\in\calw_{\pi}$ up to each $\ome_{n}(\pi+1)$, cf. (\ref{eq:4aro}) and Lemma \ref{lem:id3wf19-1}.\ref{lem:id3wf19-1.4}. Suppose $\eta\preceq\gam\in U^{*}_{N-1}(\calw^{\mbox{\footnotesize L}};\eta)$.

$\mbox{{\rm L}}\in M[\gam,N-1;\eta]$ is equivalent to
\beqnarrs
&& \fal\del\in U^{*}_{N-1}(\calw^{\mbox{\footnotesize L}};\eta)\fal b(\del\prec_{N-1}\gam \spand \mbox{{\rm L}}\models\Pi_{N-1}(b) \Rarw \\
&& \exi Q(b\in Q\spand Q\models\Pi_{N-1}(b)\spand Q\in M[\del,N-1;\eta]).
\eeqnarrs

Suppose $\del\in U^{*}_{N-1}(\calw^{\mbox{\footnotesize L}};\eta) \spand \del\prec_{N-1}\gam \spand \mbox{{\rm L}}\models\Pi_{N-1}(b)$ for a $b$. Then 
$st_{N-1}(\del)<st_{N-1}(\gam)$ by Lemma \ref{lem:5ast3}. 

By IH, $\Pi_{N-1}(b) \land \mbox{{\rm L}}\in M[\del,N-1;\eta]$ holds in $\mbox{{\rm L}}$. Since this is a $\Pi_{N}$-formula, $\Pi_{N}$-reflection for the whole universe $\mbox{{\rm L}}$ yields $\exi Q(b\in Q\spand Q\models\Pi_{N-1}(b)\spand Q\in M[\del,N-1;\eta])$.
\eprf

\bth\label{th:pi11exist}For each $n\in\ome$
\[\fal\eta\in U^{*}_{N-1}(\calw^{\mbox{\footnotesize L}})\cap V^{*}_{2}(\calw^{\mbox{\footnotesize L}})\cap\cald^{Q}[b(\eta)<\ome_{n}(\pi+1)\Rarw \mbox{{\rm L}}\in M[\eta,2]].\]
\eth
\bprf
Assume $\eta\in U^{*}_{N-1}(\calw^{\mbox{\footnotesize L}})\cap V^{*}_{2}(\calw^{\mbox{\footnotesize L}})$ and $b(\eta)<\ome_{n}(\pi+1)$ for an $\eta\in\cald^{Q}$. We show the
\bclm\label{clm:pi11exist}
\benu
\item $\eta\preceq\gam<\pi \Rarw \mbox{{\rm L}}\in M_{i}(\calx(\gam,i;\eta))$ for $i<N-1$.
\label{clm:pi11exista}
\item $\eta\preceq\gam<\pi \Rarw \mbox{{\rm L}}\in M[\gam,i;\eta]$ for $i\leq N-1$.
\label{clm:pi11existb}
\eenu
\eclm
{\bf Proof} of Claim \ref{clm:pi11exist} by simultaneous induction on $\ell\gam$ with subsidiary induction on $N-i$. 
Suppose $\eta\preceq\gam<\pi$. 

First $\gam\in U^{*}_{N-1}(\calw^{\mbox{\footnotesize L}};\eta)$ and hence by Lemma \ref{lem:4aro} we have $\mbox{L}\in M[\gam,N-1;\eta]$, i.e., Claim \ref{clm:pi11exist}.\ref{clm:pi11existb} for the case $i=N-1$ follows. In what follows assume $i<N-1$. 
\\
\ref{clm:pi11exist}.\ref{clm:pi11exista}. 
First by SIH we have $\mbox{L}\in M[\gam,i+1;\eta]\cap\bigcap\{M[\gam_{i}^{m+1},i+1;\eta]: m<lh_{i}(\gam)-1\}$. 
By reflecting $\Pi_{N}$-classes $M[\gam,i+1;\eta]$ and $M[\gam_{i}^{m+1},i+1;\eta]$ we have 
\beqnarrs
\mbox{L} & \in & M_{N}(M[\gam,i+1;\eta]\cap\bigcap\{M[\gam_{i}^{m+1},i+1;\eta]: m<lh_{i}(\gam)-1\}) \\
& \incl & M_{i}(M[\gam,i+1;\eta]\cap\bigcap\{M[\gam_{i}^{m+1},i+1;\eta]: m<lh_{i}(\gam)-1\}).
\eeqnarrs

Let $\mbox{L}\models\Pi_{i}(b)$ for a $b$.
Pick a {\it set\/} $P$ in L such that $P\in M[\gam,i+1;\eta]\cap\bigcap\{M[\gam_{i}^{m+1},i+1;\eta]: m<lh_{i}(\gam)-1\}$,
$P\models\Pi_{i}(b)$ 
and 
$\eta\in V^{*}_{2}(\calw^{P};\eta)$ by Lemma \ref{lem:5uv}.\ref{lem:5uv.3-1}. 
We claim that
\[P\in\calx(\gam,i;\eta)=M[\gam,i+1;\eta]\cap\bigcap\{M(\gam_{i}^{m},i;\eta): m<lh_{i}(\gam)-1\}.\]
This yields $\mbox{L}\in M_{i}(\calx(\gam,i;\eta))$ as desired. 

Now we show that $P\in M(\gam_{i}^{m},i;\eta)$ by reverse induction on $m<lh_{i}(\gam)-1$. 
Suppose $P\in \bigcap\{M(\gam_{i}^{m+1+k},i;\eta): k<lh_{i}(\gam)-m-2\}$ and put $\del=\gam_{i}^{m}$. 
We have to show $P\in M(\del,i;\eta)$. 
By Lemma \ref{lem:5ap12} $\del_{i}^{1}=\gam_{i}^{m+1}$, $lh_{i}(\del_{i}^{1})=lh_{i}(\gam)-m-1$ and 
$\fal k<lh_{i}(\del_{i}^{1})[(\del_{i}^{1})_{i}^{k}=\del_{i}^{1+k}=\gam_{i}^{m+1+k}]$. 
Therefore we have $P\in\calx(\del_{i}^{1},i;\eta)=M[\gam_{i}^{m+1},i+1;\eta]\cap\bigcap\{M(\gam_{i}^{m+1+k},i;\eta): k<lh_{i}(\gam)-m-2\}$.
On the other hand we have $\del\in H^{s}_{i}(\calw^{P};\eta)$ by
$\eta\preceq\del$ and $\eta\in V^{*}_{2}(\calw^{P};\eta)$.
Consequently Lemma \ref{lem:5etaMh}.\ref{lem:5etaMh.b} yields $P\in M(\del,i;\eta)$ for the set $P$. We are done.
\\
\ref{clm:pi11exist}.\ref{clm:pi11existb}. By SIH we have $\mbox{{\rm L}}\in M[\gam,i+1;\eta]$. 
Also by Claim \ref{clm:pi11exist}.\ref{clm:pi11exista} we have $\mbox{{\rm L}}\in M_{i}(\calx(\gam_{i}^{0},i;\eta))$. 
Therefore $\mbox{{\rm L}}\in M[\gam,i+1;\eta]\cap M_{i}(\calx(\gam_{i}^{0},i;\eta))=M[\gam,i;\eta]$.

Thus we have shown Claim \ref{clm:pi11exist}.
\eprf

Claim \ref{clm:pi11exist}.\ref{clm:pi11existb} for the case $\gam=\eta\spand i=2$ yields $\mbox{{\rm L}}\in M[\eta,2]$.
\eprf

Now we conclude the following theorem.

\bth\label{th:3awf16.2}For each $n\in\ome$
\[
\fal\eta\in\cald[\eta\in\calg(\calw^{\mbox{\footnotesize L}})\cap V^{*}(\calw^{\mbox{\footnotesize L}}) \spand b(\eta)<\ome_{n}(\pi+1) \Rarw \eta\in\calw^{\mbox{\footnotesize L}}].
\]
\eth
\bprf 
Without loss of generality we can assume $\eta\in\cald^{Q}$ by Lemmata \ref{lem:3wf16} and \ref{lem:4acalg}. 
Assume $\eta\in\calg(\calw^{\mbox{\footnotesize L}})\cap V^{*}(\calw^{\mbox{\footnotesize L}}) \spand b(\eta)<\ome_{n}(\pi+1)$.  
By Lemma \ref{lem:gvinclu} and Theorem \ref{th:pi11exist} we have $\mbox{{\rm L}}\in M[\eta,2]$ and hence 
$\mbox{{\rm L}}\in M_{2}(M_{2}(M[\eta,2]))$.
By Lemma \ref{lem:4acalg}, pick a limit universe $P$ such that 
$\eta\in\calg(\calw^{P})\cap V^{*}(\calw^{P}) \spand P\in M_{2}(M[\eta,2])$. 
Then Theorem \ref{th:5awf16.1} yields $\eta\in\calw^{P}\incl\calw^{\mbox{\footnotesize L}}$.
\eprf

\subsection{Wellfoundedness proof (concluded)}\label{subsec:5awf.3}

In this subsection we show $\{\sig\}\cup c(\alp_{1})\incl\calw_{\pi}\Rarw \alp_{1}\in\calw^{\mbox{\footnotesize L}}$ for each $\alp_{1}\in\cald$. Thus a proof of Theorem \ref{th:5wfeach} is completed.

By Theorem \ref{th:3awf16.2}, cf. Lemma \ref{lem:GWpiN}, we have
for each $\alp\in Od(\Pi_{N})|\pi$,
\[\alp\in\calg\cap V^{*} \Rarw \alp\in\calw^{\mbox{\footnotesize L}}\]
for
\[
\calg:=\calg(\calw^{\mbox{\footnotesize L}}) \mbox{ and } V^{*}:=V^{*}(\calw^{\mbox{\footnotesize L}}).
\]
{\bf Proof} of Lemma \ref{lem:id5wf21} for $\calw^{\mbox{\footnotesize L}}$. We have to show for each $n\in\ome$
\[
\fal\alp\in\calw_{\pi}|\ome_{n}(\pi+1)\fal q\incl\calw_{\pi}|\ome_{n}(\pi+1) A(\alp,q).
\] 
By main induction on $\alp\in\calw_{\pi}|\ome_{n}(\pi+1)$ with subsidiary induction on $q\incl\calw_{\pi}|\ome_{n}(\pi+1)$.
Here observe that if $\bet_{1}\in\cald$ with $b(\bet_{1})<\ome_{n}(\pi+1)$, then by Lemma \ref{lem:N4aro} we have $Q(\bet_{1})\leq\max\{b(\bet_{1}),\pi\}<\ome_{n}(\pi+1)$.

Let $\alp_{1}\in\cald_{\sig}$ with $\sig\in\calw_{\pi}$ and $\alp=b(\alp_{1})\spand q=Q(\alp_{1})$.
By Theorem \ref{th:id5wf21} we have $\alp_{1}\in\calg$. We show $\alp_{1}\in\calw^{\mbox{\footnotesize L}}$. 
By Theorem \ref{th:3awf16.2} it suffices to show $\alp_{1}\in V^{*}$.
We have $\sig\in V^{*}\cup\{\pi\}$ by $\sig\in\calw_{\pi}$. 
 We show the following claim.
\bclm\label{clm:5wf21.40}
Let $2\leq i<N-1$.
\benu
\item\label{clm:5wf21.40.0}
$\calc^{\alp_{1}}(\calw^{\mbox{\footnotesize L}})|\alp_{1}\incl V^{*}_{2}(\calw^{\mbox{\footnotesize L}})$.
\item\label{clm:5wf21.40.1}
$\sig\preceq\bet<\pi \Rarw \bet\in H^{s}_{\geq 2}(\calw^{\mbox{\footnotesize L}};\alp_{1})$.
\item\label{clm:5wf21.40.2}
$\eta\preceq\alp_{1} \spand \eta\in U_{i}(\calw^{\mbox{\footnotesize L}};\alp_{1}) \spand \alp_{1}\prec rg_{i}(\eta)\darw
\Rarw \eta\in V^{s}_{i}(\calw^{\mbox{\footnotesize L}};\alp_{1})$.
\item\label{clm:5wf21.40.3}
$\eta\preceq\alp_{1} \spand \alp_{1}\prec rg_{i}(\eta)\darw
\Rarw \eta\in H^{s}_{i}(\calw^{\mbox{\footnotesize L}};\alp_{1})$.
\eenu
\eclm
{\bf Proof} of Claim \ref{clm:5wf21.40}.
\\
\ref{clm:5wf21.40}.\ref{clm:5wf21.40.0}.
$\calc^{\alp_{1}}(\calw^{\mbox{\footnotesize L}})|\alp_{1}\incl V^{*}_{2}(\calw^{\mbox{\footnotesize L}})$ follows from $\alp_{1}\in\calg(\calw^{\mbox{\footnotesize L}})$.
\\
\ref{clm:5wf21.40}.\ref{clm:5wf21.40.1}.
If $\sig\preceq\bet<\pi$, then 
$\sig\in V^{*}_{2}(\calw^{\mbox{\footnotesize L}};\sig)$ yields $\bet\in H^{s}_{\geq 2}(\calw^{\mbox{\footnotesize L}};\sig)\incl H^{s}_{\geq 2}(\calw^{\mbox{\footnotesize L}};\alp_{1})$ by Lemma \ref{lem:5uv}.\ref{lem:5uv.3-1}.
\\
\ref{clm:5wf21.40}.\ref{clm:5wf21.40.2}.
Assume $\eta\preceq\alp_{1} \spand \eta\in U_{i}(\calw^{\mbox{\footnotesize L}};\alp_{1}) \spand \alp_{1}\prec rg_{i}(\eta)\darw$ for an $i<N-1$.
We show
\beqn\label{eq:5wf21.40.2}
st_{i}(\eta)\in\calw_{\pi}|\ome_{n}(\pi+1)
\eeqn
We have $\gam\lhd^{s}_{i}\eta \Rarw st_{i}(\gam)<st_{i}(\eta)$ for 
$\sig\preceq rg_{i}(\gam)=rg_{i}(\eta)$. Therefore by induction on $st_{i}(\eta)\in\calw_{\pi}|\ome_{n}(\pi+1)$ we see 
$\eta\in V^{s}_{i}(\calw^{\mbox{\footnotesize L}};\alp_{1})$.

Put $\nu=st_{i}(\eta)$ and $\tau=rg_{i}(\eta)$. 
By Lemma  \ref{lem:5.4}.\ref{lem:5.4.3} and $i<N-1$ we have $\nu=st_{i}(\eta)<\pi$. 
Thus we have shown $st_{i}(\eta)<\ome_{n}(\pi+1)$. 

On the other hand we have $\eta\in U_{i}(\calw^{\mbox{\footnotesize L}};\alp_{1})$. Namely 
$\bigcup\{K_{\kap}\nu: \kap\leq\tau\}\incl\calw^{\mbox{\footnotesize L}}$.

Hence Lemma \ref{lem:CX4aro} with the condition $(\cald.2)$ (\ref{cnd:Kst}) in Definition \ref{df:piN} yields 

\beqn\label{eq:5wf21.43}
\nu\in\calc^{\tau}(\calw^{\mbox{\footnotesize L}})
\eeqn

By Lemmata \ref{lem:5.4}.\ref{lem:5.4.3} and \ref{lem:Npi11exist} we have

\beqn\label{eq:5wf21.41}
\calb_{>\tau}(\nu)<\alp 
\eeqn

Now Lemma \ref{lem:id3wf20} together with $\mbox{MIH}(\alp)$, (\ref{eq:5wf21.43}) and (\ref{eq:5wf21.41}) yields $\nu\in\calc^{\pi}(\calw^{\mbox{\footnotesize L}})=\calw_{\pi}$. 
 This shows (\ref{eq:5wf21.40.2}).
\\
\ref{clm:5wf21.40}.\ref{clm:5wf21.40.3} by induction on $N-i$.

Assume $\eta\preceq\alp_{1} \spand \alp_{1}\prec rg_{i}(\eta)\darw$ for an $i<N-1$.
Let $U_{i}(\calw^{\mbox{\footnotesize L}};\alp_{1})\ni\gam\lhd^{s}_{i}\eta$. 
We show $\gam\in V^{*}_{i+1}(\calw^{\mbox{\footnotesize L}};\alp_{1})\cap V^{s}_{i}(\calw^{\mbox{\footnotesize L}};\alp_{1})$.

First by Claim \ref{clm:5wf21.40}.\ref{clm:5wf21.40.2} we have $\gam\in V^{s}_{i}(\calw^{\mbox{\footnotesize L}};\alp_{1})$.

Next suppose $\gam\prec_{i+1}\bet<\pi$. 
Then $\alp_{1}\prec rg_{i}(\eta)=rg_{i}(\gam)\preceq_{i}pd_{i+1}(\gam)\preceq_{i+1}\bet$ by $(\cald.11)$ in Definition \ref{df:piN}. Hence $\sig\preceq\bet$, and Claim \ref{clm:5wf21.40}.\ref{clm:5wf21.40.1} yields 
$\bet\in H^{s}_{\geq i+1}(\calw^{\mbox{\footnotesize L}};\alp_{1})$.

Finally Lemma \ref{lem:5.4}.\ref{lem:5.4.10} 
with $i\in In(\gam)$ 
yields $\alp_{1}\prec rg_{i}(\gam)\preceq_{i}rg_{j}(\gam)$ for any $j\geq i+1$ with $rg_{j}(\gam)\darw$. 
Therefore $\gam\in H^{s}_{\geq i+1}(\calw^{\mbox{\footnotesize L}};\alp_{1})$ by IH.
Consequently 
$\gam\in V^{*}_{i+1}(\calw^{\mbox{\footnotesize L}};\alp_{1})$.

This shows Claim \ref{clm:5wf21.40}.
\eprf

Now by Claim \ref{clm:5wf21.40}.\ref{clm:5wf21.40.3} we have $\alp_{1}\in H^{s}_{\geq 2}(\calw^{\mbox{\footnotesize L}})$.
On the other hand we have $\alp_{1}\prec\bet<\pi \Rarw \bet\in H^{s}_{\geq 2}(\calw^{\mbox{\footnotesize L}};\alp_{1})$
by Claim \ref{clm:5wf21.40}.\ref{clm:5wf21.40.1}, 
and hence $\alp_{1}\in V^{*}_{2}(\calw^{\mbox{\footnotesize L}};\alp_{1})$.
Thus Claim \ref{clm:5wf21.40}.\ref{clm:5wf21.40.0} yields $\alp_{1}\in V^{*}$.

This completes a proof of Lemma \ref{lem:id5wf21} for $\calw^{\mbox{\footnotesize L}}$.
\eprf

Lemma \ref{lem:id5wf21} yields Lemma \ref{th:id4wf22}: $\alp_{1}\in\calw_{\pi}$ for each $\alp_{1}\in Od(\Pi_{N})$ as in \cite{Wienpi3d}.

Consequently Lemma \ref{lem:WWOme} yields Theorem \ref{th:5wfeach}.

\end{document}